\documentclass[12pt,twoside]{article}    
\usepackage{amssymb}
\usepackage{amsmath}
\usepackage{amsfonts}%
\usepackage{graphicx}

\usepackage{url}

\providecommand{\U}[1]{\protect\rule{.1in}{.1in}}

\newtheorem{theorem}{Theorem}[section]
\newtheorem{remark}[theorem]{Remark}
\newtheorem{lemma}[theorem]{Lemma}

\newtheorem{proposition}[theorem]{Proposition}

\numberwithin{equation}{section}
\newenvironment{proof}[1][Proof]{\textbf{#1.} }{\ \rule{0.5em}{0.5em}}

\makeatletter \@addtoreset{equation}{section} \makeatother
\begin{document}


\pagestyle{myheadings}

\markboth{\hfill {\small AbdulRahman Al-Hussein and Boulakhras Gherbal} \hfill}{\hfill
{\small FBDSDEs with Poisson Jumps} \hfill }


\thispagestyle{plain}


\begin{center}
{\large \textbf{Existence and uniqueness of the
solutions of forward-backward doubly stochastic differential equations with Poisson jumps}} \\
\vspace{0.7cm} AbdulRahman Al-Hussein$^{a,}$\footnote{This work is supported by the Science College Research Center at Qassim University, project No. 3479-cos-2018-1-14-S.},
Boulakhras Gherbal$^{b,}$\footnote{This work is supported by by the Algerian PRFU, project No. C00L03UN070120180005.},
\\
\vspace{0.2cm} {\footnotesize
{\it $^{a}$Department of Mathematics, College of Science, Qassim University, \\
 P.O.Box 6644, Buraydah 51452, Saudi Arabia \\ {\emph E-mail:} alhusseinqu@hotmail.com, hsien@qu.edu.sa \\ \smallskip
 $^{b}$Laboratory of Applied Mathematics, University of Mohamed Khider, \\ P.O.Box 145, Biskra 07000, Algeria
 \\ {\emph E-mail:} bgherbal@yahoo.fr, b.gherbal@univ-biskra.dz}}
\end{center}

\vspace{0.15cm}

\begin{abstract}
The aim of this paper is to establish the existence and uniqueness of the solution to a system of nonlinear fully coupled forward-backward doubly stochastic differential equations with Poisson jumps. Our system is Markovian in the sense that initial and terminal values depend on solutions, and are not just fixed random variables. We establish under some monotonicity conditions, the existence and uniqueness of strong solutions of such equations by using a continuation method.
\end{abstract}



\bigskip

\section{Introduction}

Forward-backward stochastic differential equations (FBSDEs in short) were
first studied by Antonelli in \cite{Ant93}, where the system of such equations is driven by Brownian motion on a small time interval. The proof there relies on the fixed point theorem. Since then FBSDEs are
encountered in stochastic optimal control problem and mathematical finance.
There are also many other methods to study FBSDEs on an arbitrarily given time interval. For example, the four-step scheme approach of Ma et al. \cite{MPY}, in which the authors proved the existence and uniqueness of solutions for fully coupled FBSDEs on an arbitrarily given time interval, where the diffusion coefficients were assumed to be nondegenerate and deterministic. Their work is based on continuation method. See also Hu and Peng \cite{HP}, Pardoux and Tang \cite{PT}, Peng and Wu \cite{PW}, and Yong \cite{Y}. There is also a numerical approach for handling some linear FBSDEs as e.g. in Delarue and Menozzi \cite{DM} and Ma et al. \cite{MSZ}; see also Ma and Yong \cite{MY}.

Recently more and more research attentions are drawn towards the optimal control problem for stochastic systems with random
jumps. The reason is clear for its applicable aspect. For example, there is compelling evidence that the dynamics of prices of financial instruments exhibit jumps that cannot be adequately captured solely by diffusion processes. Several
empirical studies demonstrate the existence of jumps in stock markets, the foreign exchange
market, and bond markets. Jumps constitute also a key feature in the description of credit risk
sensitive instruments. BSDEs with jumps process (BSDEJ) had been discussed for the first time by Tang and Li \cite{TL}, and after that Situ \cite{S} proved an existence and uniqueness result for BSDEJ with non-Lipschitz coefficients and gave a probabilistic interpretation for solutions to some associated partial differential-integral equations (PDIEs). Barles et al. in \cite{BBP} and Yin and Mao \cite{YM} discussed viscosity solutions to a system of PDIEs in terms of BSDEs with jumps.

A new class of stochastic differential equations with terminal condition, called backward doubly stochastic differential equations (BDSDEs) was introduced in 1994 by Pardoux and Peng in \cite{PP}. Precisely, they proved there the existence and uniqueness of the solutions of this kind of systems and produced also a probabilistic representation of certain quasi-linear stochastic partial differential equations (SPDEs) extending a Feynman Kac formula for linear SPDEs.

The existence and uniqueness result for the solutions of backward doubly stochastic differential equations with jumps process (BDSDEJ) with Lipschitz coefficients on a fixed time interval was achieved by Sun and Lu in \cite{SL}.

Peng and Shi, \cite{PS}, introduced fully coupled FBDSDEs and showed the existence and uniqueness of their solutions with arbitrarily fixed time duration and under some monotone conditions. The equations there work on the same space, and generalize stochastic Hamiltonian systems. This result have been extended by Zhu et al. in \cite{ZSG} to different dimensional FBDSDEs and weakened the monotone assumptions. Zhu and Shi, \cite{ZS}, introduced the notion of bridge for systems of coupled FBDSDEs. They showed that, if two FBDSDEs are linked by a bridge, then they have the same unique solvability. A probabilistic interpretation for the solutions to an associated class of quasilinear SPDEs is provided there. There is also a direct link between the optimal filtering problem and FBDSDEs as it is shown in \cite{Bao}.

In the present paper we extend the results of Peng and
Shi, \cite{PS}, to FBDSDEs with jumps, which are driven particularly by Brownian motions and Poisson process.
We establish in particular the existence and uniqueness of the solutions of the following Markovian fully coupled FBDSDEs with jumps, i.e. FBDSDEsJ (or merely FBDSDEJ):
\begin{eqnarray}\label{intr-syst}
\!\!\! \left\{
\begin{array}{ll}%
\!\! dy_{t}=b\left(  t,y_{t},Y_{t},z_{t},Z_{t},k_{t}\right)
dt+\sigma\left(  t,y_{t},Y_{t},z_{t},Z_{t},k_{t}\right)  dW_{t}%
-z_{t}d\overleftarrow{B}_{t}\\ \hspace{4cm}
+\int_{\Theta}\varphi\left(  t,y_{t},Y_{t},z_{t},Z_{t},k_{t},\rho\right)
\widetilde{N}\left(  d\rho,dt\right)\\ \\
dY_{t}=f\left(  t,y_{t},Y_{t},z_{t},Z_{t},k_{t}\right)  dt+g\left(  t,y_{t},Y_{t},z_{t},Z_{t},k_{t}\right)  d\overleftarrow{B}%
_{t}+Z_{t}dW_{t}\\ \hspace{4cm}
+\int_{\Theta}k_{t}\left(  \rho\right)  \widetilde{N}\left(  d\rho,dt\right) ,
\\ \\
y_{0}=\Psi\left(  Y_{0}\right)  ,Y_{T}=h\left(  y_{T}\right) , \; \; t \in [0,T] .
\end{array}
\right.
\end{eqnarray}
The mappings $b,\sigma,\varphi,f,g$ and $h$ are given, $(W_{t})_{t\geq0}$ and
$(B_{t})_{t\geq0}$ be two mutually independent standard Brownian motions,
defined on a probability space $(\Omega,\mathcal{F},\mathbb{P}),$
taking their values respectively in $\mathbb{R}^{d}$ and in $\mathbb{R}^{l}, \, \tilde{N}(d\rho,dt)$ is a Poisson
measure with characteristic $\Pi(d\rho)dt$.
The integral with respect to $B_{t} $ is a backward It\^{o} integral,
while the integral with respect to $W_{t}$ is a standard forward
It\^{o} integral.

Note that our equations in (\ref{intr-syst}) live in different spaces $\mathbb{R}^n$ and $\mathbb{R}^m$ and are Markovian, in which the initial value $y$ and the terminal value of $Y$ depend on state solution processes. This indeed affects the assumptions set for this study as we will see in our assumptions especially when we try to prove the existence part of the solution of (\ref{intr-syst}) which is established in Theorem~\ref{Propo 3.4}. More precisely, as stated in this theorem we shall see that three cases have to be considered: $m > n, \, m < n$ and $m = n,$ and we have to distinguish between the conditions of some of the Lipschitz constants verified by the mappings $\sigma$ and $\varphi$ appearing in (\ref{intr-syst}). This matter exists when we have different dimensions and such mappings $\sigma$ and $\varphi$ that depend on variable $z.$ We shall see in Lemma~\ref{lem:final-lemma} below, which is one of our promised vital results when $m < n$, that the Lipschitz constants of the mappings $\sigma$ and $\varphi$ with respect to the variable $z$ must be less than $1/2.$ One may would like now to see quickly the statements of Theorem~\ref{Propo 3.4} and Lemma~\ref{lem:final-lemma} and compare the latter with the other case in Lemma~\ref{Lemma: 3.5}.

Our results here are new and cover also previous studies in the same field (even without jumps). In fact, in our system (\ref{intr-syst}) we have Poisson jumps and we allow all mappings to depend on all variables $(t,y,Y,z,Z,k)$ giving more flexibility and generality besides being fully-coupled and Markovian as well.

On the other hand, we are now able to apply our work here to stochastic maximum principle for
FBDSDEJ. In fact, this work furnishes now a solid ground for studying such stochastic control problems governed by FBDSDEJ as in \cite{Al-G-relaxed} in terms of the maximum principle approach. One can see also \cite{AG}. Applications of FBDSDEs with jumps to semilinear stochastic PDEs can be developed to give a probabilistic representation for the solution of a semilinear stochastic partial
differential-integral equation in parallel to \cite{PP} and \cite{BBP}.

Let us start the paper with an introductory information concerning some notions of backward filtration and backward integration. They will help us to understand the meaning of solution of (\ref{intr-syst}), and also to derive it especially in its existence part.

\section{Introduction to backward filtration and backward integrals}\label{sec2}
Let $(\Omega,\mathcal{F},\mathbb{P}) $ be a complete probability space. Let $(
W_{t}) _{t\in[ 0,T ] }$ and $( B_{t}) _{t\in[ 0,T ] }$ be two Brownian motions taking their values in $\mathbb{R}^{d}$ and $\mathbb{R}^{l}$ respectively. Let $\eta$ be a Poisson point process taking
its values in a measurable space $( \Theta,\mathcal{B}( \Theta) ) .$ In this paper, we always assume $\Theta$ is a standard Borel space in the sense of \cite{Parth} and $\mathcal{B}( \Theta)$ is the topological $\sigma$-algebra on $\Theta$ (equivalently, $\Theta$ is a
Lusin space in the sense of \cite{Bour} and $\mathcal{B}( \Theta)$ is the totality of Borel subsets of $\Theta$). See also \cite[Section~2]{Watanabe}. For example $\Theta = \mathbb{R}$ and $\mathcal{B}( \Theta) = \mathcal{B}(\mathbb{R}) .$
We denote by $\Pi( d\rho) $ the characteristic measure of $\eta$ which is assumed to be a
$\sigma$-finite measure on $( \Theta,\mathcal{B}( \Theta) ),$ by $N (d\rho,dt)
$ the Poisson counting measure (jump measure) induced by $\eta$ with compensator $\Pi( d\rho)
dt$, and by
\[
\tilde{N}( d\rho,dt) = N ( d\rho,dt) -\Pi( d\rho) dt,
\]
the compensation of the jump measure $N(\cdot, \cdot)$ of $\eta .$ Hence $\Pi (O) = \mathbb{E} [N(O,1)]$ for $O\in \mathcal{B}(\Theta) .$

We assume that these three processes $W , B$ and $\eta$ are mutually independent.

Recall that $W,B$ are Brownian motion in $\mathbb{R}^{d}$ and $\mathbb{R}^{l}$, respectively.
Let $\mathcal{F}_{t}^{W}:=\sigma\{W_{r}\mid 0\leq r\leq t \}\vee \mathcal{N},\mathcal{F}_{t}^{\eta}:=\sigma\{\eta_{r}\mid 0\leq r\leq t \}\vee \mathcal{N},$ for all $0\leq t\leq T,$ and $\mathcal{F}_{t,T}^{B}:=\sigma\{B_{r}-B_{T}\mid t\leq r\leq T\}\vee \mathcal{N},$ where $\mathcal{N}$ is the $\mathbb{P}$-null sets in $\Omega.$ Then, $\{\mathcal{F}_{t}^{W}\mid t\geq 0\}, \{\mathcal{F}_{t}^{\eta}\mid t\geq 0\}$ are two filtrations in the sense that $\mathcal{G}_{t}\subset\mathcal{G}_{s}$, if $t\leq s,$ for $\mathcal{G}=\mathcal{F}^{W}$ or $\mathcal{G}=\mathcal{F}^{\eta}$, while $\{\mathcal{F}_{t,T}^{B}\mid 0\leq t\leq T\}$ is a \emph{backward filtration}  in the sense that $\mathcal{F}_{t,T}^{B}\supseteq \mathcal{F}_{s,T}^{B}$ if $s\leq t.$ If $\{ M_{t}\mid 0\leq t\leq T\}$ is a stochastic process over $(\Omega,\mathcal{F},\mathbb{P})$ satisfying: $M_{t}$ is $\mathcal{F}_{t,T}^{B}$-measurable $\, \forall \, 0\leq t\leq T$, we say that $M$ is $\{{\mathcal{F}_{t,T}^{B}, t\leq T}\}$-\emph{adapted}.

If $M$ is adapted and $\mathbb{E}\, [|M_{t}|]<+\infty,$ for all $0\leq t\leq T,$ we say that $M$ is a \emph{backward martingale} if $\mathbb{E}\, [M_{t}\mid \mathcal{F}_{s,T}^{B}]=M_{s}, \, \forall \, t\leq s.$

Let $\{h_{t}\mid 0\leq t\leq T\}$ be an $\{\mathcal{F}^{B}_{t,T}\mid 0\leq t\leq T\}$-adapted $\mathbb{R}^{n\times d}$-valued process, satisfying $\mathbb{E}\, [\int_{0}^{T}|h_{s}|^{2}ds]<\infty$. Recall that the backward It\^{o} integral of $h$ with respect to $B$ is defined by
\begin{eqnarray}\label{eq:2.1}
&& \int_{\alpha}^{\beta}h_{s}\overleftarrow{dB}_{s}:=\lim_{|\pi|\rightarrow 0}\overset{n}{\underset{i=1}{\mathbb{\sum}}}h(t_{i+1})(B_{t_{i+1}}-B_{t_{i}}), \; \; (\text{in} \; L^{2}(\Omega,\mathcal{F},\mathbb{P})),
\end{eqnarray}
where $\pi=\{\alpha=t_{1},t_{2},\cdots,t_{n+1}=\beta\}$ is a partition of $[\alpha,\beta]$ satisfying $$|\pi|=\rm{mesh}\, \pi :=\underset{1\leq i\leq n}{\max}(t_{i+1}-t_{i})\rightarrow 0.$$
Observe that if $h$ is constant, i.e, $h_{s}=c \; \forall \, 0\leq s\leq T,$ for some constant $c\in \mathbb{R}^{n},$ then
\begin{eqnarray}\label{eq:2.2}
&& \int_{\alpha}^{\beta}c \, \overleftarrow{dB}_{s}=c\, (B_{\beta}-B_{\alpha})=\int_{\alpha}^{\beta}c \, dB_{s}.
\end{eqnarray}
Hence, $$\int_{0}^{u}1 \, \overleftarrow{dB}_{s}=B_{u},\int_{0}^{T-t}1 \, \overleftarrow{dB}_{s}=B_{T-t},\; \int_{t}^{T}1 \, \overleftarrow{dB}_{s}=B_{T}-B_{t}.$$
It is easy to see that $M_{t}=\int_{t}^{T}h_{s}\overleftarrow{dB}_{s},0\leq t\leq T,$ is $\{\mathcal{F}_{t,T},0\leq t\leq T\}$-backward martingale.

On the other hand, letting $\breve{B}_{s}:=B_{T-s}-B_{T},0\leq s\leq T,$ shows that $\breve{B}$ is a Brownian motion as well, and for all $0\leq s\leq T$, we have
\begin{eqnarray*}
&& \mathcal{F}_{T-t,T}^{B}=\sigma\{ B_{r}-B_{T}\mid T-t\leq r\leq T\}\\
&& =\sigma\{ B_{T-s}-B_{T}\mid 0\leq s\leq t\},(r=T-s),\\
&& =\sigma\{ \breve{B}_{s}\mid 0\leq s\leq t\}=\mathcal{F}_{t}^{\breve{B}}.
\end{eqnarray*}
Similarly, if $\breve{W}_{s}:=W_{T-s}-W_{T},0\leq s\leq T$, then $\breve{W}$ is a Brownian motion, and $$\mathcal{F}_{T-t,T}^{\breve{W}}=\mathcal{F}_{t}^{W},\,\forall \,0\leq t\leq T,$$
or, in particular $$\mathcal{F}_{T-t}^{W}=\mathcal{F}_{t,T}^{\breve{W}},\,\forall \,0\leq t\leq T.$$
Consequently, if $h_{s}$ is $\mathcal{F}_{s,T}^{B}$-measurable for all $0\leq s\leq T,$ then $\breve{h}_{s}:=h_{T-s}$ is $\mathcal{F}_{s}^{\breve{B}}$-measurable for all $0\leq s\leq T,$ and if $k_{s}$ is $\mathcal{F}_{s}^{W}$-measurable for all $0\leq s\leq T,$ then $\breve{k}_{s}:=k_{T-s}$ is $\mathcal{F}_{T-s,T}^{\breve{W}}$-measurable for all $0\leq s\leq T.$

Moreover, if we fix $t\in[0,T],$ and let $\pi=\{t=t_{1},t_{2},\cdots,t_{n+1}=T\}$ be a partition of $[t,T]$ satisfying $\rm{mesh}~\pi \rightarrow 0,$ then $\widehat{\pi}=\{0=s_{1},s_{2},\cdots,s_{n}=T-t\}$ is a partition of $[0,T-t],$ where $s_{i}=T-t_{n+2-i}$ for all $1\leq i\leq n+1,$ i.e, $$s_{1}=T-t_{n+1}=0,s_{2}=T-t_{n},s_{3}=T-t_{n-1},\cdots,s_{n-1}=T-t_{3},$$ $$s_{n}=T-t_{2},s_{n+1}=T-t_{1}=T-t,$$ and $$\rm{mesh}~\widehat{\pi}=\underset{1\leq i\leq n+1}{\max}|s_{i+1}-s_{i}|=\underset{1\leq i\leq n+1}{\max}|t_{i}-t_{i+1}|=\rm{mesh}~\pi \rightarrow 0.$$

Let $h$ be as in (\ref{eq:2.1}). Then, if $\breve{h}_{s}:=h_{T-s}, \; 0\leq s\leq T-t,$
\begin{eqnarray*}
&& \hspace{-0.75cm}\int_{0}^{T-t}\breve{h}_{s}d\breve{B}_{s}=\lim_{n\rightarrow \infty}\overset{n}{\underset{i=1}{\mathbb{\sum}}}\breve{h}(s_{i})(\breve{B}_{s_{i+1}}-\breve{B}_{s_{i}})\\
&& \hspace{0.25cm}=\lim_{n\rightarrow \infty}[\breve{h}(s_{1})(\breve{B}_{s_{2}}-\breve{B}_{s_{1}})+\breve{h}(s_{2})(\breve{B}_{s_{3}}-\breve{B}_{s_{2}})\\
&& \hspace{1.25cm}+\cdots +\breve{h}(s_{n-1})(\breve{B}_{s_{n}}-\breve{B}_{s_{n-1}})+\breve{h}(s_{n})(\breve{B}_{s_{n+1}}-\breve{B}_{s_{n}})]\\
&& \hspace{0.25cm}=\lim_{n\rightarrow \infty}[\breve{h}(0)(\breve{B}_{T-t_{n}}-\breve{B}_{0})+\breve{h}(T-t_{n})(\breve{B}_{T-t_{n-1}}-\breve{B}_{T-t_{n}})\\
&& \hspace{1.25cm}+\cdots +\breve{h}(T-t_{3})(\breve{B}_{T-t_{2}}-\breve{B}_{T-t_{3}})+\breve{h}(T-t_{2})(\breve{B}_{T-t_{1}}-\breve{B}_{T-t_{2}})]\\
&& \hspace{0.25cm}=\lim_{n\rightarrow \infty}[h(T)(B_{t_{n}}-B_{T})+h(t_{n})(B_{t_{n-1}}-B_{t_{n}})\\
&& \hspace{1.25cm}+\cdots +h(t_{3})(B_{t_{2}}-B_{t_{3}})+h(t_{2})(B_{t_{1}}-B_{t_{2}})]\\
&& \hspace{0.25cm}=-\lim_{n\rightarrow \infty}[h(t_{2})(B_{t_{2}}-B_{t_{1}})+h(t_{3})(B_{t_{3}}-B_{t_{2}})\\
&& \hspace{1.25cm}+\cdots +h(t_{n})(B_{t_{n}}-B_{t_{n-1}})+h(t_{n+1})(B_{t_{n+1}}-B_{t_{n}})]\\
&& =-\int_{t}^{T}h_{s}\overleftarrow{dB}_{s}.
\end{eqnarray*}
Hence
\begin{eqnarray}\label{eq:2.3}
&& \int_{t}^{T}h_{s}\overleftarrow{dB}_{s}=-\int_{0}^{T-t}\breve{h}_{s}d\breve{B}_{s},\; 0\leq t\leq T.
\end{eqnarray}
In particular,
\begin{eqnarray}\label{eq:2.4}
&& \int_{T-u}^{T}h_{s}\overleftarrow{dB}_{s}=-\int_{0}^{u}\breve{h}_{s}d\breve{B}_{s},\; 0\leq u\leq T.
\end{eqnarray}
Therefore,
\begin{eqnarray}\label{eq:2.5}
\int_{0}^{T-t}h_{s}\overleftarrow{dB}_{s}&=& \int_{0}^{T}h_{s}\overleftarrow{dB}_{s}-\int_{T-t}^{T}h_{s}\overleftarrow{dB}_{s},\nonumber \\
&=&-\int_{0}^{T}\breve{h}_{s}d\breve{B}_{s}+\int_{0}^{t}\breve{h}_{s}d\breve{B}_{s},\nonumber \\
&=& -\int_{t}^{T}\breve{h}_{s}d\breve{B}_{s},\; 0\leq t\leq T,
\end{eqnarray}
or, in particular,
\begin{eqnarray}\label{eq:2.6}
&& \int_{0}^{u}h_{s}\overleftarrow{dB}_{s}=-\int_{T-u}^{T}\breve{h}_{s}d\breve{B}_{s},\,0\leq u\leq T.
\end{eqnarray}

These information (\ref{eq:2.1})-(\ref{eq:2.6}) hold evidently when  $h_{t}$ is $\mathcal{F}_{t}=\mathcal{F}_{t}^{W}\vee \mathcal{F}_{t,T}^{B}\vee \mathcal{F}_{t}^{\eta}\vee \mathcal{N}$-measurable for all $0\leq t\leq T,$ and $\mathbb{E}\, [\int_{0}^{T}|h_{t}|^{2}dt]<\infty,$ i.e, when $h$ is an element of the space $\mathcal{M}(0,T,\mathbb{R}^{n}),$ introduced below.

\medskip

Let us close this section by introducing the following list of notations.

\medskip

For each $t\in[ 0,T ] $ define $\mathcal{F}_{t}\triangleq\mathcal{F}%
_{t}^{W}\vee\mathcal{F}_{t,T}^{B}\vee\mathcal{F}_{t}^{\eta}$.
For a Euclidean space $E,$ let $\mathcal{M}^{2}(0,T; E) $ denote the set of jointly measurable,
processes $\left\{  X_{t},t\in\big[ 0,T  \big] \right\}
$ with values in $E$ such that $X_t$ is $\mathcal{F}_{t}$-measurable for a.e. $t \in [0,T],$ and satisfy
\[
\mathbb{E} \big[  {\displaystyle\int_{0}^{T}} \left\vert \mathcal{X}_{t}\right\vert_E ^{2}dt \big]
<\infty.
\]

Let $L_{\Pi}^{2}( E)$ be the set of $\mathcal{B}( \Theta)$-measurable mapping $k$ with values in $E$
such that
\[
||| k ||| := \big[ \int_{\Theta}\left\vert k(
\rho) \right\vert_{E}^{2}\Pi( d\rho) \big]^{\frac{1}{2}} < \infty .
\]

Denote by $\mathcal{N}_{\eta}^{2}( 0,T ;E) $ to the set of $L_{\Pi}^{2}( E)$-valued processes $\{ \mathcal{K}_t , \; t \in \big[0,T \big] \}$ that satisfy: $\mathcal{K}_t$ is $\mathcal{F}_{t}$-measurable for a.e. $t \in [0,T], $ and
\[
\mathbb{E} \big[ {\displaystyle\int_{0}^{T}} \int_{\Theta}\left\vert \mathcal{K}_{t}( \rho)
\right\vert_{E}^{2}\Pi( d\rho) dt \big] <\infty.
\]

We denote%
\begin{eqnarray*}
&& \hspace{-1cm} \mathbb{H}^{2} := \mathcal{M}^{2}\left(  0,T ;\mathbb{R}^{n}\right)
\times\mathcal{M}^{2}\left( 0,T;\mathbb{R}^{m}\right)
\times\mathcal{M}^{2}\left(0,T;\mathbb{R}^{n\times l}\right)  \\ && \hspace{2.25in}
\times \mathcal{M}^{2}\left(0,T ;\mathbb{R}^{m\times d}\right)
\times\mathcal{N}_{\eta}^{2}\left(0,T ;\mathbb{R}^{m}\right)  .
\end{eqnarray*}
Then $\mathbb{H}^{2}$ is a Hilbert space
equipped with the norm given by:
\begin{equation*}
\left\Vert \zeta_{\cdot}\right\Vert _{\mathbb{H}^{2}}^{2} \triangleq\mathbb{E}\left[ \int_{0}^{T} \big{(} \left\vert
y_{t}\right\vert ^{2} + \left\vert Y_{t}%
\right\vert ^{2} + \left\Vert z_{t}\right\Vert^{2} + \left\Vert Z_{t}\right\Vert^{2} + |||
k_{t}||| ^{2} \big{)} dt\right]  ,
\end{equation*}
where $\zeta_{\cdot}=\left(  y_{\cdot},Y_{\cdot},z_{\cdot},Z_{\cdot},k_{\cdot
}  \right) \in \mathbb{H}^{2} .$ We shall sometimes use the notation
\[ \left\vert
\zeta_{t}\right\vert^{2} := \left\vert
y_{t}\right\vert ^{2} + \left\vert Y_{t}%
\right\vert ^{2} + \left\Vert z_{t}\right\Vert^{2} + \left\Vert Z_{t}\right\Vert^{2} + |||
k_{t}||| ^{2} .\]

Similarly, let
\begin{eqnarray*}
&& \hspace{-1cm} \widehat{\mathbb{H}}^{2} := \mathcal{M}^{2}\left(  0,T ;\mathbb{R}^{m}\right)
\times\mathcal{M}^{2}\left( 0,T;\mathbb{R}^{n}\right)
\times\mathcal{M}^{2}\left(0,T;\mathbb{R}^{m\times l}\right)  \\ && \hspace{2.25in}
\times \mathcal{M}^{2}\left(0,T ;\mathbb{R}^{n\times d}\right)
\times\mathcal{N}_{\eta}^{2}\left(0,T ;\mathbb{R}^{n}\right) .
\end{eqnarray*}

\section{The existence and uniqueness theorem of FBDSDEJ}\label{sec3}
Let
\begin{eqnarray*}
&& b :[  0,T]  \times \mathbb{R}^{n}\times
\mathbb{R}^{m}\times \mathbb{R}^{n\times l}\times \mathbb{R}^{m\times d}\times L_{\Pi  }^{2}(  \mathbb{R}^{m})  \rightarrow\mathbb{R}^{n},\\
&& \sigma :[  0,T]  \times\mathbb{R}^{n}\times\mathbb{R}^{m}\times\mathbb{R}^{n\times l}\times\mathbb{R}^{m\times d}\times L_{\Pi }^{2}(  \mathbb{R}^{m})\rightarrow\mathbb{R}^{n\times d},\\
 && \varphi :[  0,T]  \times\mathbb{R}^{n}\times\mathbb{R}^{m}\times\mathbb{R}^{n\times l}\times\mathbb{R}^{m\times d}\times L_{\Pi  }^{2}(  \mathbb{R}^{m})  \times \Theta \rightarrow\mathbb{R}^{n},\\
 && f  :[  0,T]  \times\mathbb{R}^{n}\times\mathbb{R}^{m}\times\mathbb{R}^{n\times l}\times\mathbb{R}^{m\times d}\times L_{\Pi  }^{2}(  \mathbb{R}^{m}) \rightarrow
\mathbb{R}^{m},\\
&& g :[  0,T]  \times\mathbb{R}^{n}\times\mathbb{R}^{m}\times\mathbb{R}^{n\times l}\times\mathbb{R}^{m\times d}\times L_{\Pi  }^{2}(  \mathbb{R}^{m})\rightarrow\mathbb{R}^{m\times l},\\
 && \Psi:\Omega \times\mathbb{R}^{m}\rightarrow\mathbb{R}^{n},\\
 && h :\Omega \times\mathbb{R}^{n}\rightarrow\mathbb{R}^{m},
\end{eqnarray*}
be a given mappings satisfying assumptions to be given shortly.

Consider the following system of equations, which we call \emph{forward-backward doubly stochastic differential equations with jumps} (FBDSDEJ):
\begin{eqnarray}\label{eq:3.1}
\left\{
\begin{array}{ll}%
dy_{t}=b\left(  t,y_{t},Y_{t},z_{t},Z_{t},k_{t}\right)
dt+\sigma\left(  t,y_{t},Y_{t},z_{t},Z_{t},k_{t}\right)  dW_{t}%
-z_{t}d\overleftarrow{B}_{t}\\ \hspace{4cm}
+\int_{\Theta}\varphi\left(  t,y_{t},Y_{t},z_{t},Z_{t},k_{t},\rho\right)
\widetilde{N}\left(  d\rho,dt\right)\\ \\
dY_{t}=f\left(  t,y_{t},Y_{t},z_{t},Z_{t},k_{t}\right)  dt+g\left(  t,y_{t},Y_{t},z_{t},Z_{t},k_{t}\right)  d\overleftarrow{B}%
_{t}+Z_{t}dW_{t}\\ \hspace{4cm}
+\int_{\Theta}k_{t}\left(  \rho\right)  \widetilde{N}\left(  d\rho,dt\right) ,
\\ \\
y_{0}=\Psi\left(  Y_{0}\right)  ,Y_{T}=h\left(  y_{T}\right)  .
\end{array}
\right.
\end{eqnarray}

A \emph{solution} of (\ref{eq:3.1}), is a stochastic process $(y,Y,z,Z,k)$ such that for each $t\in[0,T]$ we have $a.s:$
\begin{eqnarray}\label{eq:3.2}
\left\{
\begin{array}{ll}%
y_{t}=\Psi\left(  Y_{0}\right)+\int_{0}^{t}b\left(  s,y_{s},Y_{s},z_{s},Z_{s},k_{s}\right)
ds+\int_{0}^{t}\sigma\left(  s,y_{s},Y_{s},z_{s},Z_{s},k_{s}\right)  dW_{s}%
\\ \hspace{2cm}
-\int_{0}^{t}z_{t}d\overleftarrow{B}_{s}+\int_{0}^{t}\int_{\Theta}\varphi\left(  s,y_{s},Y_{s},z_{s},Z_{s},k_{s},\rho\right)
\widetilde{N}\left(  d\rho,ds\right)\\ \\
Y_{t}=h\left(  y_{T}\right)-\int_{t}^{T}f\left(  s,y_{s},Y_{s},z_{s},Z_{s},k_{s}\right)  ds-\int_{t}^{T}g\left(  s,y_{s},Y_{s},z_{s},Z_{s},k_{s}\right)  d\overleftarrow{B}%
_{s}\\ \hspace{4cm}
-\int_{t}^{T}Z_{t}dW_{s}-\int_{t}^{T}\int_{\Theta}k_{s}\left(  \rho\right)  \widetilde{N}\left(  d\rho,ds\right).
\end{array}
\right.
\end{eqnarray}

Given an $m\times n$\ full-rank matrix $R$ let us introduce the following
notation:
$$\upsilon  =\left(  y,Y,z,Z,k\right)  ,$$
$$A\left(  t,\upsilon\right)    :=\left(  R^{\ast}f,Rb,R^{\ast}g,R\sigma,R\varphi\right)  \left(
t,\upsilon\right)  ,$$
$$\left\langle A,\upsilon\right\rangle  :=\left\langle
y,R^{\ast}f\right\rangle +\left\langle Y,Rb\right\rangle
+\left\langle z,R^{\ast}g\right\rangle +\left\langle Z,R\sigma
\right\rangle +\left\langle \left\langle k,R\varphi\right\rangle
\right\rangle ,$$
where ($^{\ast}$) denotes matrix transpose, and%
\begin{align*}
R^{\ast}g  & =\left(  R^{\ast}g_{1},\cdots,R^{\ast}g_{l}\right)
,R\sigma=\left(  R\sigma_{1},\cdots,R\sigma_{d}\right)  ,\\
\left\langle \left\langle k,R\varphi\right\rangle \right\rangle  &
=\int_{\Theta}\left\langle k_{t}\left(  \rho\right)  ,R\varphi\left(
t,\upsilon_{t},\rho\right)  \right\rangle \Pi\left(  d\rho\right)  .
\end{align*}

We set the following assumptions.

(A1) (Monotonicity condition): $\, \forall \,\upsilon=\left(
y,Y,z,Z,k\right)  ,\overline{\upsilon}=\left(  \overline{y},\overline{Y},\overline{z},\overline{Z}
,\overline{k}\right)  \in\mathbb{R}^{n+m+n\times l+m\times d}\times L^{2}_{\Pi}(\mathbb{R}^{m})$, $\, \forall \, t\in\left[  0,T\right]  $%
\begin{align*}
\left\langle A\left(  t,\upsilon\right)  -A\left(
t,\overline{\upsilon}\right)  ,\upsilon-\overline{\upsilon}\right\rangle  & \leq
-\theta_{1}\left(  \left\vert R\left(  y-\overline{y}\right)  \right\vert
^{2}+\left\Vert R\left(  z-\overline{z}\right)  \right\Vert^{2}\right)  \\
& \hspace{-1cm}-\theta_{2}\left(  \left\vert R^{\ast}\left(  Y-\overline{Y}\right)  \right\vert
^{2}+\left\Vert R^{\ast}\left(  Z-\overline{Z}\right)  \right\Vert^{2}+||| R^{\ast}\left(  k-\overline{k}\right)|||^{2}\right)  .
\end{align*}

(A2) We have
\begin{align*}
\left\langle \Psi\left(  Y\right)  -\Psi\left(  \overline{Y}\right)  ,R^{\ast
}\left(  Y-\overline{Y}\right)  \right\rangle  & \leq-\beta_{2}\left\vert R^{\ast
}\left(  Y-\overline{Y}\right)  \right\vert ^{2},\; \forall \; Y,\overline{Y}\in\mathbb{R}%
^{m},\\
\left\langle h\left(  y\right)  -h\left(  \overline{y}\right)  ,R\left(  y-\overline{y}\right)  \right\rangle  & \geq\beta_{1}\left\vert R\left(  y-\overline{y}\right)  \right\vert ^{2}, \; \forall \; y,\overline{y}\in\mathbb{R}%
^{n}.
\end{align*}
Here $\theta_{1},\theta_{2},\beta_{1},$\ and $\beta_{2}$\ are given
nonnegative constants with ${\theta_{1}+\theta_{2}>0},$ ${\beta_{1}+\beta
_{2}>0,} ~ {\theta_{1}+\beta_{2}>0,} ~ {\theta_{2}+\beta_{1}>0.}$ Moreover, we have
$\theta_{1}>0,\beta_{1}>0$\ (resp. $\theta_{2}>0,\beta_{2}>0$) when
$m>n$\ (resp. $n>m$).

\medskip

(A3) For each $\upsilon\in\mathbb{H}^{2}, \, A\left(
t,\upsilon\right)  $ is an $\mathcal{F}_{t}$-measurable vector process defined
on $\left[  0,T\right]  $ with $\widetilde{A}\left(  .,0\right)  \in\mathbb{H}%
^{2}$ and $\, \forall \, y\in\mathbb{R}^{n},h\left(  y\right)  $ is an
$\mathcal{F}_{T}$-measurable vector process with $h\left(  0\right)  \in
L^{2}\left(  \Omega,\mathcal{F}_{T},\mathbb{P};\mathbb{R}^{m}\right),$ and for each
$Y\in\mathbb{R}^{m},\Psi\left(  Y\right)  $ is an $\mathcal{F}_{0}$-measurable
vector process with $\Psi\left(  0\right)  \in L^{2}\left(  \Omega
,\mathcal{F}_{0},\mathbb{P};\mathbb{R}^{n}\right)  .$

\medskip

(A4) (Lipschitz condition): $\exists \, c>0, \gamma$ and $\gamma'>0$ such that
\begin{eqnarray*}
&& \hspace{-0.75cm} \left\vert b\left(  t,y,Y,z,Z,k\right)  -b\left(  t,\overline{y},\overline{Y},\overline{z},\overline{Z},\overline{k}\right)  \right\vert^2 \\
&& \leq c\left(  \left\vert y-\overline{y}\right\vert ^{2}+\left\vert Y-\overline{Y}\right\vert ^{2}+\left\Vert z-\overline{z}\right\Vert^{2}+\left\Vert Z-\overline{Z}\right\Vert^{2}+||| k-\overline{k}|||^{2}\right)  ,\\
&& \hspace{-0.75cm} \left\vert f\left(  t,y,Y,z,Z,k\right)  -f\left(  t,\overline{y},\overline{Y},\overline{z},\overline{Z},\overline{k}\right)  \right\vert ^{2}\\
&& \leq c\left(  \left\vert y-\overline{y}\right\vert ^{2}+\left\vert Y-\overline{Y}\right\vert ^{2}+\left\Vert z-\overline{z}\right\Vert^{2}+\left\Vert Z-\overline{Z}\right\Vert^{2}+||| k-\overline{k}|||^{2}\right)  ,\\
&& \hspace{-0.75cm} \left\vert \sigma\left(  t,y,Y,z,Z,k\right)  -\sigma\left(
t,\overline{y},\overline{Y},\overline{z},\overline{Z},\overline{k}\right)  \right\vert ^{2}\\
&& \leq c\left(  \left\vert y-\overline{y}\right\vert ^{2}+\left\vert Y-\overline{Y}\right\vert ^{2}+\left\Vert Z-\overline{Z}\right\Vert^{2}+||| k-\overline{k}||| ^{2}\right)+\gamma' \left\Vert z-\overline{z}\right\Vert^{2}  , \\
&& \hspace{-0.75cm} \left\vert g\left(  t,y,Y,z,Z,k\right)  -g\left(
t,\overline{y},\overline{Y},\overline{z},\overline{Z},\overline{k}\right)  \right\vert ^{2}\\
&& \leq c\left(  \left\vert y-\overline{y}\right\vert ^{2}+\left\vert Y-\overline{Y}\right\vert ^{2}+\left\Vert z-\overline{z}\right\Vert^{2}\right)
+\gamma\left(  \left\Vert Z-\overline{Z}\right\Vert^{2}+||| k-\overline{k}||| ^{2}\right)  ,\\
&& \hspace{-0.75cm} \underset{\rho\in\Theta}{\sup}\left\vert \varphi\left(  t,y,Y,z,Z,k,\rho
\right)  -\varphi\left(  t,\overline{y},\overline{Y},\overline{z},\overline{Z},\overline{k}
,\rho\right)  \right\vert ^{2}\\
&& \leq c\left(  \left\vert y-\overline{y}\right\vert ^{2}+\left\vert Y-\overline{Y}\right\vert ^{2}+ \left\Vert Z-\overline{Z}\right\Vert^{2}+||| k-\overline{k}|||^{2}\right)+\gamma' \left\Vert z-\overline{z}\right\Vert^{2}  ,\\
&& \hspace{-0.75cm} \left\vert \Psi\left(  Y\right)  -\Psi\left(  \overline{Y}\right)  \right\vert
\leq c\left\vert Y-\overline{Y}\right\vert ,\left\vert h\left(  y\right)  -h\left(
\overline{y}\right)  \right\vert \leq c\left\vert y-\overline{y}\right\vert ,
\end{eqnarray*}
for all argument written in the right hand side of each inequality and for all $t\in [0,T]$ whenever they appear.
\begin{remark}\label{Remark 3.1}
Assumptions (A1) and (A2) can be replaced by  the following ones with essentially the same proofs of the solution theorem and its lemmas.

(A1)$'$ $\, \forall \,\upsilon=\left(  y,Y%
,z,Z,k\right)  ,\overline{\upsilon}=\left(  \overline{y},\overline{Y}
,\overline{z},\overline{Z},\overline{k}\right)  \in\mathbb{R}^{n+m+n\times l+m\times d}\times L^{2}_{\Pi}(\mathbb{R}^{m}),$ $\, \forall \,
t\in\left[  0,T\right]  $%
\begin{align*}
& \left\langle A\left(  t,\upsilon\right)  -A\left(
t,\overline{\upsilon}\right)  ,\upsilon-\overline{\upsilon}\right\rangle \geq\theta_{1}\left(  \left\vert R\left(  y-\overline{y}\right)  \right\vert
^{2}+\left\Vert R\left(  z-\overline{z}\right)  \right\Vert^{2}\right)  \\
& \hspace{3cm} +\theta_{2}\left(  \left\vert R^{\ast}\left(  Y-\overline{Y}\right)  \right\vert
^{2}+\left\Vert R^{\ast}\left(  Z-\overline{Z}\right)  \right\Vert^{2}+||| R^{\ast}\left(  k-\overline{k}\right)||| ^{2}\right) ,
\end{align*}
and

(A2)$'$
\begin{align*}
\left\langle \Psi\left(  Y\right)  -\Psi\left(  \overline{Y}\right)  ,R^{\ast
}\left(  Y-\overline{Y}\right)  \right\rangle  & \geq\beta_{2}\left\vert R^{\ast
}\left(  Y-\overline{Y}\right)  \right\vert ^{2},\; \forall \; Y,\overline{Y}\in\mathbb{R}%
^{m},\\
\left\langle h\left(  y\right)  -h\left(  \overline{y}\right)  ,R\left(  y-\overline{y}\right)  \right\rangle  & \leq-\beta_{1}\left\vert R\left(  y-\overline{y}\right)  \right\vert ^{2},\; \forall \; y,\overline{y}\in\mathbb{R}^{n}.
\end{align*}
We shall consider only assumptions (A1)--(A4).
\end{remark}

\medskip

\begin{remark}\label{Remark 3.2}
 Note that, since $R:\mathbb{R}^{n}\rightarrow\mathbb{R}^{m}$
(a matrix $m\times n$) has a full-rank, if $m>n,$ then $R$ is injective and there exists $C_{1},C_{2}>0$ such that
\begin{align*}
C_{1}\left\vert y\right\vert _{\mathbb{R}^{n}}  & \leq\left\vert
Ry\right\vert _{\mathbb{R}^{m}}\leq C_{2}\left\vert y\right\vert
_{\mathbb{R}^{n}}, \, \forall \, y\in \mathbb{R}^{n}.
\end{align*}
In fact $C_{1}={\rm min}\{|Rx|_{\mathbb{R}^{m}}, x\in S^{n-1}\},$ which is positive since $R$ is injective, and $C_{2}={\rm max}\{|Rx|_{\mathbb{R}^{m}}, x\in S^{n-1}\}=\|R\|.$

In the case $m<n,$ we have a similar inequality for the transpose matrix $R^{\ast},$ which is injective as a mapping from $\mathbb{R}^{m}$ to $\mathbb{R}^{n},$ namely:
\begin{align*}
C_{3}\left\vert Y\right\vert _{\mathbb{R}^{m}}  & \leq\left\vert
R^{\ast}Y\right\vert _{\mathbb{R}^{n}}\leq C_{4}\left\vert Y\right\vert
_{\mathbb{R}^{m}}, \, \forall \, Y\in \mathbb{R}^{m},
\end{align*}
where $C_{3}={\rm min}\{|R^{\ast}Y|_{\mathbb{R}^{n}};Y\in S^{m-1}\}>0,$ and $C_{1}={\rm max}\{|R^{\ast}Y|_{\mathbb{R}^{n}};Y\in S^{m-1}\}=\|R^{\ast}\|.$

If $m=n,$ then ${\rm rank}(R)=m=n,$ and so $R$ is invertible, and we get the following result:
\begin{align*}
C_{5}\left\vert y\right\vert _{\mathbb{R}^{n}}  & \leq\left\vert
Ry\right\vert _{\mathbb{R}^{n}}\leq C_{6}\left\vert y\right\vert
_{\mathbb{R}^{n}}, \, \forall \, y\in \mathbb{R}^{n},
\end{align*}
where $C_{5}=\|R^{-1}\|^{-1}$ and $C_{6}=\|R\|.$
\end{remark}

Let us state on It\^{o}'s formula (see e.g. \cite{SL}), that will be used throughout the paper.
\begin{proposition}\label{Lemma 3.3}
Let $(\alpha,\widehat{\alpha})\in
\big[ \mathcal{M}^{2}(0,T;\mathbb{R}^{n}) \big]^{2},(\beta,\widehat{\beta})\in\big[ \mathcal{M}^{2}(0,T;\mathbb{R}^{n}) \big]^{2},
(\gamma,\widehat{\gamma})\in\big[ \mathcal{M}^{2}(0,T;\mathbb{R}^{n\times l}) \big]^{2},(\delta,\widehat{\delta})
\in \big[ \mathcal{M}^{2}(0,T;\mathbb{R}^{n\times d}) \big]^{2}$, and $(K,\widehat{K})\in
 \big[\mathcal{N}_{\eta}^{2}(0,T;\mathbb{R}^{n})\big]^{2} .$ Assume that
\[
\alpha_{t}=\alpha_{0}+\int_{0}^{t}\beta_{s}ds+\int_{0}^{t}\gamma_{s}%
\overleftarrow{dB}_{s}+\int_{0}^{t}\delta_{s}dW_{s}+\int_{0}^{t}\int_{\Theta}K_{s}%
(\rho)\tilde{N}(d\rho,ds),
\]
and
\[
\widehat{\alpha}_{t}=\widehat{\alpha}_{0}+\int_{0}^{t}\widehat{\beta}_{s}ds+\int_{0}^{t}\widehat{\gamma}_{s}%
\overleftarrow{dB}_{s}+\int_{0}^{t}\widehat{\delta}_{s}dW_{s}+\int_{0}^{t}\int_{\Theta}\widehat{K}_{s}(\rho)\tilde{N}(d\rho,ds),
\]
for $t \in [0,T] .$
Then, for each $t \in [0,T] ,$
\begin{eqnarray*}
&& \langle \alpha_{t},\widehat{\alpha}_{t}\rangle =\langle \alpha_{0},\widehat{\alpha}_{0}\rangle+\int_{0}^{t}\left\langle \alpha_{s},d\widehat{\alpha}_{s}\right\rangle+\int_{0}^{t}\left\langle d\alpha_{s},\widehat{\alpha}_{s}\right\rangle+\int_{0}^{t}d\left\langle \alpha,\widehat{\alpha}\right\rangle_{s} \;\; \; \mathbb{P}-a.s.,
\end{eqnarray*}
and
\begin{eqnarray*}
&& \mathbb{E}\, \big[
\langle \alpha_{t},\widehat{\alpha}_{t}\rangle \big] =  \mathbb{E}\, \big[ \langle \alpha_{0},\widehat{\alpha}_{0}\rangle \big] + \mathbb{E}\,  \big[ \int_{0}^{t}\left\langle \alpha_{s},d\widehat{\alpha}_{s}\right\rangle \big]+ \mathbb{E}\,  \big[ \int_{0}^{T}\left\langle d\alpha_{s},\widehat{\alpha}_{s}\right\rangle \big]
\\ && \hspace{0.4in} - \, \mathbb{E}\, \big[ \int_{0}^{t}\langle\gamma_{s},\widehat{\gamma}_{s}\rangle ds \big]+\mathbb{E}\, \big[ \int_{0}^{t}\langle\delta_{s},\widehat{\delta}_{s}\rangle ds \big]+\mathbb{E}\, \big[ \int_{0}^{t}\int_{\Theta}\langle K_{s},\widehat{K}_{s}\rangle \Pi(  d\rho) ds \big].
\end{eqnarray*}
\end{proposition}

\begin{theorem}[Uniqueness]\label{Propo 3.4-1}
Under assumptions (A1)--(A4) (or (A1)$'$, (A2)$'$, (A3)--(A4); cf. Remark~\ref{Remark 3.1}) FBDSDEJ~(\ref{eq:3.1}) has at most one solution
$\left(y,Y,z,Z,k\right)$  in $\mathbb{H}^{2}.$
\end{theorem}
\begin{proof} Let $\upsilon:=\left(  y,Y%
,z,Z,k\right)  $ and $\upsilon^{\prime}:=\left(  y%
^{\prime},Y^{\prime},z^{\prime},Z^{\prime},k^{\prime}\right)
$ be two solutions of (\ref{eq:3.1}). Denote $$\widehat{\upsilon}=\left(  \widehat{y},\widehat
{Y},\widehat{z},\widehat{Z},\widehat{k}\right)  =\left(  y-y^{\prime},Y-Y^{\prime
},z-z^{\prime},Z-Z^{\prime},k-k^{\prime}\right)  .$$ Applying integration of
parts (Proposition~\ref{Lemma 3.3}) to $\left\langle R\widehat{y},\widehat{Y}\right\rangle $\ on $\left[
0,T\right]  $ we get
\begin{align*}
& \mathbb{E}\, \left[  \left\langle R\widehat{y}_{T},h\left(  y_{T}\right)  -h\left(
y_{T}^{\prime}\right)  \right\rangle \right]  -\mathbb{E}\, \left[  \left\langle
R\left(  \Psi\left(  Y_{0}\right)  -\Psi\left(  Y_{0}^{\prime}\right)
\right)  ,\widehat{Y}_{0}\right\rangle \right]  \\
& \hspace{2cm} =\mathbb{E}\, \left[  \int_{0}^{T}\left\langle A\left(  t,\upsilon
_{t}\right)  -A\left(  t,\upsilon_{t}^{\prime}\right)  ,\widehat
{\upsilon}_{t}\right\rangle dt\right]  \\
& \hspace{2cm} \leq-\theta_{1}\, \mathbb{E} \, \left[  \int_{0}^{T}\left(  \left\vert R\widehat{y}%
_{t}\right\vert ^{2}+\left\Vert R\widehat{z}_{t}\right\Vert^{2}\right)
dt\right]  \\
& \hspace{2cm} -\theta_{2}\, \mathbb{E} \, \left[  \int_{0}^{T}\left(  \left\vert R^{\ast}\widehat
{Y}_{t}\right\vert ^{2}+\left\Vert R^{\ast}\widehat{Z}_{t}\right\Vert
^{2}+||| R^{\ast}\widehat{k}_{t}|||^{2}\right)  dt\right]  .
\end{align*}
Thus
\begin{eqnarray*}
&& \theta_{1}\, \mathbb{E} \, \left[  \int_{0}^{T}\left(  \left\vert R\widehat{y}%
_{t}\right\vert ^{2}+\left\Vert R\widehat{z}_{t}\right\Vert^{2}\right)
dt\right]\\
&& \hspace{0.5in}  +\theta_{2}\, \mathbb{E} \, \left[  \int_{0}^{T}\left(  \left\vert
R^{\ast}\widehat{Y}_{t}\right\vert ^{2}+\left\Vert R^{\ast}\widehat{Z}_{t}\right\Vert
^{2}+||| R^{\ast}\widehat{k}_{t}|||^{2}\right)  dt\right]  \leq0.
\end{eqnarray*}

If $m>n$\ then $\theta_{1}>0$, and so we have $\left\vert R\widehat{y}_{t}\right\vert
^{2}\equiv0$ and $\left\Vert R\widehat{z}_{t}\right\Vert^{2}\equiv0$. Thus $y_{t}=
y_{t}^{\prime}$ and $z_{t}= z^{\prime}$. In particular, $h\left(  y_{T}\right)
=h\left(  y_{T}^{\prime}\right)  $. The system (\ref{eq:3.2}) becomes
\begin{eqnarray*}
\left\{
\begin{array}{ll}%
0=\Psi\left(  Y_{0}\right)-\Psi\left(  Y'_{0}\right)+\int_{0}^{t}\widehat{b}\left(  s,Y_{s},Z_{s},k_{s}\right)
ds+\int_{0}^{t}\widehat{\sigma}\left(  s,Y_{s},Z_{s},k_{s}\right)  dW_{s}%
\\ \hspace{5cm}
+\int_{0}^{t}\int_{\Theta}\widehat{\varphi}\left(  s,Y_{s},Z_{s},k_{s},\rho\right)
\widetilde{N}\left(  d\rho,ds\right)\\ \\
\widehat{Y}_{t}=-\int_{t}^{T}\widehat{f}\left(  s,Y_{s},Z_{s},k_{s}\right)  ds-\int_{t}^{T}\widehat{g}\left(  s,Y_{s},Z_{s},k_{s}\right)  d\overleftarrow{B}%
_{s}\\ \hspace{2.2in}
-\int_{t}^{T}\widehat{Z}_{t}dW_{s}-\int_{t}^{T}\int_{\Theta}\widehat{k}_{s}\left(  \rho\right)  \widetilde{N}\left(  d\rho,ds\right),
\end{array}
\right.
\end{eqnarray*}
where
\begin{align*}
\widehat{\pi}\left(  t,\upsilon_{t}\right)
:=\pi\left(  t,\upsilon_{t}\right)  -\pi \left(  t,\upsilon
_{t}^{\prime}\right)
\end{align*}
for $\pi  :=b,\sigma,f,g,$ and
\[
\widehat{\varphi}\left(  t, \upsilon_{t},\cdot\right) :=\varphi
\left(  t, \upsilon_{t},\cdot\right)  -\varphi\left(  t, \upsilon_{t}^{\prime},\cdot\right).\]

Consequently, from the uniqueness of
solution of BDSDEJ (see Yin and Situ \cite{YS}), it follows that $Y_{t}= Y_{t}^{\prime},Z_{t}=
Z_{t}^{\prime}$\ and $k_{t}= k_{t}^{\prime}$.

If $m<n$\ then $\theta_{2}>0$, and so we have $Y_{t}= Y_{t}^{\prime}%
,Z_{t}= Z_{t}^{\prime}$\ and $k_{t}= k_{t}^{\prime}$. In particular,
$\Psi\left(  Y_{0}\right)  =\Psi\left(  Y_{0}^{\prime}\right)  ,$ and so system (\ref{eq:3.2}) becomes
\begin{eqnarray*}
\left\{
\begin{array}{ll}%
\widehat{y}_{t}=\int_{0}^{t}\widehat{b}\left(  s,y_{s},z_{s}\right)
ds+\int_{0}^{t}\widehat{\sigma}\left(  s,y_{s},z_{s}\right)  dW_{s}%
\\ \hspace{2cm}
-\int_{0}^{t} \widehat{z}_{t}d\overleftarrow{B}_{s}+\int_{0}^{t}\int_{\Theta}\widehat{\varphi}\left(  s,y_{s},z_{s},\rho\right)
\widetilde{N}\left(  d\rho,ds\right)\\ \\
0=\widehat{h}\left(  y_{T}\right)-\int_{t}^{T}\widehat{f}\left(  s,y_{s},z_{s}\right)  ds-\int_{t}^{T}\widehat{g}\left(  s,y_{s},z_{s}\right)  d\overleftarrow{B}%
_{s}.
\end{array}
\right.
\end{eqnarray*}
As done earlier, we apply the uniqueness of the solution of BDSDEJ in \cite{YS} to deduce that
$y_{t}= y_{t}^{\prime}$ and $z_{t}= z^{\prime}.$

By arguments similar to the above cases, the desired uniqueness property of solutions can be obtained
easily in the case $m=n$
\end{proof}

\begin{theorem}[Existeness]\label{Propo 3.4}
Assume (A1)--(A4) or ((A1)$'$, (A2)$'$, (A3)--(A4)) with $0< \gamma<1$ and $0 <\gamma' < 1$ when $m>n$, while  when $m\leq n$ we assume $0<\gamma<1$ and $0<\gamma' \leq \gamma/2$. Then FBDSDEJ (\ref{eq:3.1}) has a solution
$\left(y,Y,z,Z,k\right)$  in $\mathbb{H}^{2}.$
\end{theorem}

We divide the proof into three case: $m>n,m<n$ and $m=n.$
\\
\textbf{Case~1:} $m>n$ $\left(  \theta_{1}>0\right)  .$
Consider the following family of FBDSDEJ parameterized by $\alpha\in\left[
0,1\right]  :$%
\begin{eqnarray}\label{eq:3.3}
\left\{
\begin{array}{ll}%
dy_{t}=\left[  \alpha b\left(  t,\upsilon_{t}\right)  +\widetilde{b}%
_{0}\left(  t\right)  \right]  dt+\left[  \alpha \sigma\left(
t,\upsilon_{t}\right)  +\widetilde{\sigma}_{0}\left(  t\right)  \right]
dW_{t}-z_{t}d\overleftarrow{B}_{t}\\ \hspace{5cm}
+\int_{\Theta}\left[  \alpha\varphi\left(  t,\upsilon_{t},\rho\right)
+\varphi_{0}\left(  t,\rho\right)  \right]  \widetilde{N}\left(
d\rho,dt\right)\\ \\
dY_{t}=\left[  \alpha f\left(  t,\upsilon_{t}\right)  -\left(
1-\alpha\right)  \theta_{1}Ry_{t}+\widetilde{f}_{0}\left(  t\right)  \right]
dt+Z_{t}dW_{t}\\  \hspace{2cm}
+\left[  \alpha g\left(  t,\upsilon_{t}\right)  -\left(
1-\alpha\right)  \theta_{1}Rz_{t}+\widetilde{g}_{0}\left(  t\right)  \right]
d\overleftarrow{B}_{t}+\int_{\Theta}k_{t}\left(  \rho\right)  \widetilde
{N}\left(  d\rho,dt\right)  ,
\\ \\
y_{0}=\alpha\Psi\left(  Y_{0}\right)  +\psi,Y_{T}=\alpha h\left(
y_{T}\right)  +\left(  1-\alpha\right)  \theta_{1}Ry_{T}+\phi,
\end{array}
\right.
\end{eqnarray}
where $\upsilon_{t}:=\left(  y_{t},Y_{t},z_{t},Z_{t},k_{t}\right)  ,\left(
\widetilde{b}_{0},\widetilde{f}_{0},\widetilde{\sigma}_{0},\widetilde{g}_{0},\varphi
_{0}\right)  \in\mathbb{H}^{2}$, $\psi\in L^{2}\left(  \Omega,\mathcal{F}%
_{0},\mathbb{P};\mathbb{R}^{n}\right)$ and $\phi\in L^{2}\left(  \Omega
,\mathcal{F}_{T},\mathbb{P};\mathbb{R}^{m}\right)  $\ are given arbitrarily.

\bigskip

Note that when $\alpha=1$ the existence of the solution of (\ref{eq:3.3}) implies clearly
that of (\ref{eq:3.2}) simply by letting $(\widetilde{b}_{0},\widetilde{f}_{0},\widetilde{\sigma}_{0},\widetilde{g}_{0},\varphi_{0})=(0,0,0,0,0)$, while when $\alpha=0$, (\ref{eq:3.3}) becomes a decoupled FBDSDEJ of the form:
\begin{eqnarray}\label{eq:3.4}
\left\{
\begin{array}{ll}%
dy_{t}=\widetilde{b}_{0}\left(  t\right)dt+\widetilde{\sigma}_{0}\left(  t\right)dW_{t}-z_{t}d\overleftarrow{B}_{t}+\int_{\Theta}\varphi_{0}\left(  t,\rho\right)\widetilde{N}\left(
d\rho,dt\right)\\ \\
dY_{t}=- \theta_{1}Ry_{t}dt+Z_{t}dW_{t}+\left[ -\theta_{1}Rz_{t}+\widetilde{g}_{0}\left(  t\right)  \right]
d\overleftarrow{B}_{t}+\int_{\Theta}k_{t}\left(  \rho\right)  \widetilde{N}\left(  d\rho,dt\right),
\\ \\
y_{0}=\psi,Y_{T}= \theta_{1}Ry_{T}+\phi,0\leq t\leq T.
\end{array}
\right.
\end{eqnarray}
Using (\ref{eq:2.4}) (with $W$ replacing $\breve{B}$) and (\ref{eq:2.6}) we can rewrite $$\int_{0}^{t}\widetilde{\sigma}_{0}(s)dW_{s}=-\int_{T-t}^{T}\breve{\widetilde{\sigma}}_{0}(s)\overleftarrow{d\breve{W}}_{s},$$ where $$\breve{\widetilde{\sigma}}_{0}(s)=\widetilde{\sigma}_{0}(T-s),$$ i.e., as a backward It\^{o} integral, and $$\int_{0}^{t}Z_{s}\overleftarrow{dB}_{s}=-\int_{T-t}^{T}\breve{Z}_{s}d\breve{B}_{s},\breve{Z}_{s}=Z_{T-s},$$
as a forward It\^{o} integral, and similar, for Lebesgue integral and integrals with respect to $\widetilde{N}$ in (\ref{eq:3.2}), which then enable us to rewrite the first (forward) equation as a BDSDE with jumps as follows:

\begin{eqnarray}\label{eq:3.5}
&& \hspace{-1cm} \breve{y}_t  =\breve{y}_T+\int_{t}^{T}\breve{\widetilde{b}}_{0}(s)ds-\int_{t}^{T}\breve{\widetilde{\sigma}}_{0}(s)\overleftarrow{dW}_{s}-\int_{t}^{T}\breve{z}_{s}d\breve{B}_{s} \nonumber \\
&& \hspace{2cm}-\int_{t}^{T}\int_{\Theta}\breve{\varphi}_{0}\left(  s,\rho\right)\breve{\widetilde{N}}\left(d\rho,ds\right),0\leq  t\leq T,
\end{eqnarray}
where $\breve{y}_t=y_{T-t},$ (hence $\breve{y}_{T}=y_{0}),\breve{\widetilde{b}}_{0}(s)=\widetilde{b}(T-s),\breve{\varphi}_{0}\left(  s,\rho\right)=\varphi(T-s,\rho).$
So under our assumptions (A3) we deduce from Yin and Situ \cite{YS} that (\ref{eq:3.4}) has a unique solution $(\breve{y},\breve{z})$ in $\mathcal{M}^{2}\left( 0,T;\mathbb{R}^{n}\right)
\times\mathcal{M}^{2}\left(0,T;\mathbb{R}^{n\times l}\right).$

One can simply apply a generalized martingale representation theorem (as in Al-Hussien and Gherbal \cite{AG}) to get an explicit formula for this unique solution $(\breve{y},\breve{z})$ since all integrals here do not depend on $\breve{y}$ or $\breve{z}.$

Consequently, $\{(y_{s},z_{s}):=(\breve{y}_{T-s},\breve{z}_{T-s}),0\leq t\leq T\}$ is the unique solution of the forward equation of (\ref{eq:3.4}). $(y,z)$ in the second equation of (\ref{eq:3.3}), afterwhich it is becomes a BDSDE of type (\ref{eq:3.4}) discussed earlier, so it admits a unique solution $(Y,Z,k).$ Therefore, we derive a unique solution $\upsilon=(y,Y,z,Z,k)$ of (\ref{eq:3.4}) in $\mathbb{H}^{2}.$

The following apriori lemma is a key step in the proof of the method of
continuation. It shows that for a fixed $\alpha=\alpha_{0}\in[
0,1)  ,$ if (\ref{eq:3.3}) is uniquely solvable, then it is also uniquely solvable
for any $\alpha\in\left[  \alpha_{0},\alpha_{0}+\delta_{0}\right]  $, for some
positive constant $\delta_{0}$\ independent of $\alpha_{0}$.
\begin{lemma}\label{Lemma: 3.5}
We assume $m>n$. Under assumption (A1)--(A4), with $0< \gamma<1$ and $0< \gamma' <1 ,$ there exists a positive constant $\delta_{0}$\ such
that if, apriori, for each $\psi\in L^{2}\left(  \Omega,\mathcal{F}%
_{0},\mathbb{P};\mathbb{R}^{n}\right)  ,\phi\in L^{2}\left(  \Omega,\mathcal{F}%
_{T},\mathbb{P};\mathbb{R}^{m}\right)  $\ and $\left(  \widetilde{b}_{0},f%
_{0},\widetilde{\sigma}_{0},\widetilde{g}_{0},\varphi_{0}\right)  \in\mathbb{H}^{2}$,
(\ref{eq:3.3}) is uniquely solvable for some $\alpha_{0}\in[  0,1)  $, then
for each $\alpha\in\left[  \alpha_{0},\alpha_{0}+\delta_{0}\right]  ,\psi\in
L^{2}\left(  \Omega,\mathcal{F}_{0},\mathbb{P};\mathbb{R}^{n}\right)  ,\phi\in
L^{2}\left(  \Omega,\mathcal{F}_{T},\mathbb{P};\mathbb{R}^{m}\right)  $\ and $\left(
\widetilde{b}_{0},\widetilde{f}_{0},\widetilde{\sigma}_{0},\widetilde{g}_{0},\varphi
_{0}\right)  \in\mathbb{H}^{2}$, (\ref{eq:3.3}) is uniquely solvable.
\end{lemma}
\begin{proof}  Assume that, for each $\psi\in L^{2}\left(  \Omega,\mathcal{F}%
_{0},\mathbb{P};\mathbb{R}^{n}\right)  ,\phi\in L^{2}\left(  \Omega,\mathcal{F}%
_{T},\mathbb{P};\mathbb{R}^{m}\right)  $\ and $(  \widetilde{b}_{0},\widetilde{f}
_{0},%
\widetilde{\sigma}_{0},\widetilde{g}_{0},\varphi_{0})  \in\mathbb{H}^{2},$\
there exists a unique solution of (\ref{eq:3.3}) for $\alpha=\alpha_{0}$. Then, for each
$\overline{\upsilon}_{\cdot}:=\left(  \overline{y}_{\cdot},\overline{Y}_{\cdot},\overline{z}_{\cdot},\overline{Z}
_{\cdot},\overline{k}_{\cdot}\right)  \in\mathbb{H}^{2}$, there exists a unique element
$\upsilon_{\cdot}:=\left(  y_{\cdot},Y_{\cdot},z_{\cdot},Z_{\cdot},k_{\cdot}\right)$  of $\mathbb{H}^{2}$\ satisfying the following FBDSDEJ:
\begin{eqnarray}\label{eq:3.6}
\left\{
\begin{array}{ll}%
dy_{t}=\left[  \alpha_{0}b\left(  t,\upsilon_{t}\right)
+\delta b\left(  t,\overline{\upsilon}_{t}\right)  +\widetilde{b}_{0}\left(
t\right)  \right]  dt-z_{t}d\overleftarrow{B}_{t}\\ \hspace{1.5cm}
+\left[  \alpha_{0}\sigma\left(  t,\upsilon_{t}\right)  +\delta
\sigma\left(  t,\overline{\upsilon}_{t}\right)  +\widetilde{\sigma}_{0}\left(
t\right)  \right]  dW_{t}\\ \hspace{1.5cm}
+\int_{\Theta}\left[  \alpha_{0}\varphi\left(  t,\upsilon_{t}%
,\rho\right)  +\delta\varphi\left(  t,\overline{\upsilon}_{t},\rho\right)
+\varphi_{0}\left(  t,\rho\right)  \right]  \widetilde{N}\left(  d\rho,dt\right)\\ \\
dY_{t}=\left[  \alpha_{0}f\left(  t,\upsilon_{t}\right)  -\left(
1-\alpha_{0}\right)  \theta_{1}Ry_{t}+\delta\left(  f\left(
t,\overline{\upsilon}_{t}\right)  +\theta_{1}R \overline{y}_{t}\right)  +\widetilde{f}%
_{0}\left(  t\right)  \right]  dt\\ \hspace{1.5cm}
+\left[  \alpha_{0}g\left(  t,\upsilon_{t}\right)  -\left(
1-\alpha_{0}\right)  \theta_{1}Rz_{t}+\delta\left(  g\left(
t,\overline{\upsilon}_{t}\right)  +\theta_{1}R\overline{z}_{t}\right)  +\widetilde{g}%
_{0}\left(  t\right)  \right]  d\overleftarrow{B}_{t}\\ \hspace{2.5in}
+\, Z_{t}~dW_{t}+\int_{\Theta}k_{t}\left(  \rho\right)  \widetilde{N}\left(
d\rho,dt\right) ,
\\ \\
y_{0}=\alpha_{0}\Psi\left(  Y_{0}\right)  +\delta\Psi\left(  \overline{Y}
_{0}\right)  +\psi ,\\
Y_{T}=\alpha_{0}h\left(  y_{T}\right)  +\left(  1-\alpha_{0}\right)
Ry_{T}+\delta\left(  h\left(  \overline{y}_{T}\right)  -R\overline{y}
_{T}\right)  +\phi,
\end{array}
\right.
\end{eqnarray}
where $\delta\in(  0,1)  $ is a parameter independent of $\alpha_{0}$ which is small enough, and will be determined later in the proof.

We shall prove that the mapping defined by:
\begin{eqnarray*}
&&\hspace{-1cm}  I_{\alpha_{0}+\delta}\left(
\overline{\upsilon}_{\cdot},\overline{y}_{T},\overline{Y}_{0}\right) := \left(  \upsilon_{\cdot},y_{T},Y_{0}\right) :\mathbb{H}^{2}\times
L^{2}\left(  \Omega,\mathcal{F}_{T},\mathbb{P};\mathbb{R}^{n}\right)  \times
L^{2}\left(  \Omega,\mathcal{F}_{0},\mathbb{P};\mathbb{R}^{m}\right)\\
 && \hspace{5cm} \rightarrow\mathbb{H}^{2}\times L^{2}\left(  \Omega,\mathcal{F}%
_{T},\mathbb{P};\mathbb{R}^{n}\right)  \times L^{2}\left(  \Omega,\mathcal{F}%
_{0},\mathbb{P};\mathbb{R}^{m}\right),
\end{eqnarray*}
is contraction.

Let
\[
\overline{\upsilon}_{\cdot} :=\left(  \overline{y}_{\cdot},\overline{Y}
_{\cdot},\overline{z}_{\cdot},\overline{Z}_{\cdot},\overline{k}_{\cdot}\right), \overline{\upsilon}_{\cdot}^{\prime} :=\left(  \overline{y}_{\cdot}^{\prime},\overline{Y}
_{\cdot}^{\prime},\overline{z}_{\cdot}^{\prime},\overline{Z}_{\cdot}^{\prime},\overline{k}_{\cdot}^{\prime
}\right)  \in\mathbb{H}^{2},
\]
and let
\begin{align*}
(\left(  y_{\cdot},Y_{\cdot},z_{\cdot},Z_{\cdot},k_{\cdot}\right), y_{T},Y_{0})= \left({\upsilon}_{\cdot},{y}_{T},{Y}_{0}\right) & := I_{\alpha_{0}+\delta
}\left(  \overline{\upsilon}_{\cdot}, \overline{y}_{T},\overline{Y}_{0}\right)  ,\\
(\left(  y_{\cdot}^{\prime},Y_{\cdot}^{\prime},z_{\cdot}^{\prime},Z_{\cdot}^{\prime},k_{\cdot}^{\prime}\right), y^{\prime}_{T},Y^{\prime}_{0})  = (\upsilon_{\cdot}^{\prime}, y^{\prime}_{T},Y^{\prime}_{0}) & :=I_{\alpha_{0}+\delta
}\left(  \overline{\upsilon}^{\prime}_{\cdot}, \overline{y}^{\prime}_{T},\overline{Y}^{\prime}_{0} \right).
\end{align*}
Set the notation
\begin{align*}
\left(  \widehat{\overline{y}}_{t},\widehat{\overline{Y}}_{t},\widehat{\overline{z}}%
_{t},\widehat{\overline{Z}}_{t},\widehat{\overline{k}}_{t}\right)    & =\widehat
{\overline{\upsilon}}_{t}=\overline{\upsilon}_{t}-\overline{\upsilon}_{t}^{\prime}\\
& =\left(  \overline{y}_{t}-\overline{y}_{t}^{\prime},\overline{Y}_{t}-\overline{Y}_{t}^{\prime
},\overline{z}_{t}-\overline{z}_{t}^{\prime},\overline{Z}_{t}-\overline{Z}_{t}^{\prime},\overline{k}_{t}-\overline{k}_{t}^{\prime}\right)  ,\\
\left(  \widehat{y}_{t},\widehat{Y}_{t},\widehat{z}_{t},\widehat{Z}%
_{t},\widehat{k}_{t}\right)    & =\widehat{\upsilon}_{t}=\upsilon_{t}%
-\upsilon_{t}^{\prime}\\
& =\left(  y_{t}-y_{t}^{\prime},Y_{t}-Y_{t}^{\prime},z_{t}-z_{t}^{\prime
},Z_{t}-Z_{t}^{\prime},k_{t}-k_{t}^{\prime}\right).
\end{align*}
We shall use these definitions and notations in Lemma~\ref{lem:final-lemma} below as well.

\medskip

Applying It\^{o}'s formula to $\left\langle R\widehat{y},\widehat{Y}\right\rangle
$\ on $\left[  0,T\right]  $ yields
\begin{eqnarray*}
&& \hspace{-0.75cm}  d\left\langle R\widehat{y}_{t},\widehat{Y}_{t}\right\rangle  =\left\langle R\widehat{y}_{t},\alpha_{0}\widehat{f}\left(
t,\upsilon_{t}\right)  -\left(  1-\alpha_{0}\right)  \theta_{1}R\widehat{y}%
_{t}+\delta\left(  \widehat{f}\left(  t,\overline{\upsilon}_{t}\right)
+\theta_{1}R\widehat{\overline{y}}_{t}\right)  \right\rangle dt\\
&& \hspace{-0.5cm}  +\left\langle \alpha_{0}R\widehat{b}\left(
t,\upsilon_{t}\right)  +\delta R\widehat{b}\left(  t,\overline{\upsilon}
_{t}\right),\widehat{Y}_{t}  \right\rangle dt+\left\langle \alpha_{0}R\widehat{\sigma
}\left(  t,\upsilon_{t}\right)  +\delta R\widehat{\sigma}\left(
t,\overline{\upsilon}_{t}\right)  ,\widehat{Z}_{t}\right\rangle dt\\
&& \hspace{-0.5cm}+\int_{\Theta
}\left\langle R^{\ast}\widehat{Y}_{t},\alpha_{0}\widehat{\varphi}\left(
t,\upsilon_{t},\rho\right)  +\delta \widehat{\varphi}\left(  t,\overline{\upsilon}
_{t},\rho\right)  \right\rangle
\widetilde{N}\left(  d\rho,dt\right)\\
&& \hspace{-0.5cm}  +\left\langle R\widehat{z}_{t},\alpha_{0}\widehat{g}\left(
t,\upsilon_{t}\right)  -\left(  1-\alpha_{0}\right)  \theta_{1}R\widehat{z}%
_{t}+\delta\left(  \widehat{g}\left(  t,\overline{\upsilon}_{t}\right)
+\theta_{1}R\widehat{\overline{z}}_{t}\right)  \right\rangle dt\\
&& \hspace{-0.5cm}  +\int_{\Theta}\left\langle \alpha_{0}R\widehat{\varphi}\left(
t,\upsilon_{t},\rho\right)  +\delta R\widehat{\varphi}\left(
t,\overline{\upsilon}_{t},\rho\right)  ,\widehat{k}_{t}\left(  \rho\right)
\right\rangle \Pi\left(  d\rho\right)  dt\\
&& \hspace{-0.5cm}  -\left\langle R^{\ast}\widehat{Y}_{t},\widehat{z}_{t}d\overleftarrow{B}%
_{t}\right\rangle +\left\langle R^{\ast}\widehat{Y}_{t},\left(\alpha_{0}\widehat{\sigma}\left(  t,\upsilon
_{t}\right)  +\delta\widehat{\sigma}\left(  t,\overline{\upsilon}
_{t}\right)\right)  dW_{t} \right\rangle\\
&& \hspace{-0.5cm}  +\left\langle R\widehat{y}_{t},\left(\alpha_{0}\widehat{g}\left(
t,\upsilon_{t}\right)  -\left(  1-\alpha_{0}\right)  \theta_{1}R\widehat{z}%
_{t}+\delta\left(  \widehat{g}\left(  t,\overline{\upsilon}_{t}\right)
+\theta_{1}R\widehat{\overline{z}}_{t}\right) \right) d\overleftarrow{B}%
_{t} \right\rangle \\
&& \hspace{-0.5cm}  +\left\langle R\widehat{y}_{t},\widehat{Z}_{t}dW_{t}\right\rangle +\int_{\Theta
}\left\langle R\widehat{y}_{t},\widehat{k}_{t}\left(  \rho\right)  \right\rangle
\widetilde{N}\left(  d\rho,dt\right)  ,
\end{eqnarray*}
where
\begin{align*}
\widehat{\pi}\left(  t,\overline{\upsilon}_{t}\right)    & :=\pi
\left(  t,\overline{\upsilon}_{t}\right)  -\pi \left(  t,\overline{\upsilon}
_{t}^{\prime}\right)  ,\widehat{\pi}\left(  t,\upsilon_{t}\right)
:=\pi\left(  t,\upsilon_{t}\right)  -\pi \left(  t,\upsilon
_{t}^{\prime}\right)
\end{align*}
for $\pi  :=b,\sigma,f,g,$
\[
\widehat{\varphi}\left(  t, \overline{\upsilon}_{t},\cdot\right) :=\varphi
\left(  t, \overline{\upsilon}_{t},\cdot\right)  -\varphi\left(  t, \overline{\upsilon}_{t}^{\prime},\cdot\right), \widehat{\varphi}\left(  t, \upsilon_{t},\cdot\right) :=\varphi
\left(  t, \upsilon_{t},\cdot\right)  -\varphi\left(  t, \upsilon_{t}^{\prime},\cdot\right).
\]
So, by integrating from $0$ to $T$, using the identities $$\widehat{y}_{0}=\alpha_{0}\widehat{\Psi}\left(
Y_{0}\right)  +\delta \widehat{\Psi}\left(  \overline{Y}_{0}\right),$$ and $$\widehat{Y}_{T}=\alpha_{0}\widehat{h}\left(  y_{T}\right)  +\left(
1-\alpha_{0}\right)  R\widehat{y}_{T}+\delta \left(  \widehat{h}\left(  \overline{y}_{T}\right)  -R\widehat{\overline{y}}
_{T}\right),$$ where $$
\widehat{\Psi}\left(  Y_{0}\right)  :=\Psi\left(  Y_{0}\right)-\Psi\left(  Y'_{0}\right),
\widehat{\Psi}\left(  \overline{Y}_{0}\right) :=\Psi\left(  \overline{Y}_{0}\right)-\Psi\left(  \overline{Y}'_{0}\right),$$
$$\widehat{h}\left(  y_{T}\right)  :=h\left(  y_{T}\right)-h\left(  y'_{T}\right),
\widehat{h}\left(  \overline{y}_{T}\right)  :=h\left(  \overline{y}_{T}\right)-h\left(  \overline{y}'_{T}\right),$$ and taking the expectation, it follows that
\begin{eqnarray*}
&& \hspace{-1.25cm} \mathbb{E}\, \left[  \left\langle R\widehat{y}_{T},\alpha_{0}\widehat{h}\left(
y_{T}\right)  +\left(  1-\alpha_{0}\right)  R\widehat{y}_{T}+\delta\left(
\widehat{h}\left(  \overline{y}_{T}\right)  -R\widehat{\overline{y}}_{T}\right)
\right\rangle \right] \\
&& \hspace{-0.5cm}   = \, \mathbb{E}\, \left[  \left\langle \alpha_{0}R\widehat{\Psi}\left(
Y_{0}\right)  +\delta R\widehat{\Psi}\left(  \overline{Y}_{0}\right)  ,\widehat{Y}%
_{0}\right\rangle \right]  -\left(  1-\alpha_{0}\right)  \theta_{1}%
\, \mathbb{E} \, \left[  \int_{0}^{T}\left\vert R\widehat{y}_{t}\right\vert
^{2}dt\right]\\
&&  +\, \delta\, \mathbb{E} \, \left[  \int_{0}^{T}\left\langle R\widehat{y}_{t}%
,\widehat{f}\left(  t,\overline{\upsilon}_{t}\right)  \right\rangle
dt\right]  +\theta_{1}\, \delta\, \mathbb{E} \, \left[  \int_{0}^{T}\left\langle
R\widehat{y}_{t},R\widehat{\overline{y}}_{t}\right\rangle dt\right]\\
&&  +\, \delta\, \mathbb{E} \, \left[  \int_{0}^{T}\left\langle \widehat{Y}_{t}  , R\widehat
{b}\left(  t,\overline{\upsilon}_{t}\right) \right\rangle dt\right]
+\delta\, \mathbb{E} \, \left[  \int_{0}^{T}\left\langle R\widehat{\sigma%
}\left(  t,\overline{\upsilon}_{t}\right)  ,\widehat{Z}_{t}\right\rangle dt\right] \\
&& + \, \delta\, \mathbb{E} \, \left[
\int_{0}^{T}\left\langle R\widehat{z}_{t},\widehat{g}\left(
t,\overline{\upsilon}_{t}\right)  \right\rangle dt\right]-\, \left(  1-\alpha_{0}\right)  \theta_{1}\, \mathbb{E} \, \left[  \int_{0}%
^{T}\left\vert \left\vert R\widehat{z}_{t}\right\vert \right\vert^{2}dt\right]\\
&& +\, \delta \, \theta_{1}\, \mathbb{E} \, \left[  \int_{0}^{T}\left\langle R\widehat{z}%
_{t},R\widehat{\overline{z}}_{t}\right\rangle dt\right]  +\delta\, \mathbb{E} \, \left[
\int_{0}^{T}\left\langle \left\langle R\widehat{\varphi}\left(  t%
,\overline{\upsilon}_{t},\cdot\right)  ,\widehat{k}_{t}\right\rangle \right\rangle
dt\right] \\
&& \hspace{4.25cm}  +\, \alpha_{0}\, \mathbb{E} \, \left[  \int_{0}^{T}\left\langle A\left(
t,\upsilon_{t}\right)  -A\left(  t,\upsilon_{t}^{\prime}\right)
,\widehat{\upsilon}_{t}\right\rangle dt\right] .
\end{eqnarray*}
This implies
\begin{eqnarray}\label{eq:3.7}
&& \hspace{-0.5cm}\mathbb{E}\, \left[  \left\langle R\widehat{y}_{T},\alpha_{0}\widehat{h}\left(
y_{T}\right)  +\left(  1-\alpha_{0}\right)  R\widehat{y}_{T}\right\rangle \right]
+\left(  1-\alpha_{0}\right)  \theta_{1}\, \mathbb{E} \, \left[  \int_{0}^{T}\left(
\left\vert R\widehat{y}_{t}\right\vert ^{2}+\left\Vert R\widehat{z}_{t}\right\Vert
^{2}\right)  dt\right]\nonumber \\
&& \hspace{6.5cm}-\, \alpha_{0} \, \mathbb{E} \, \left[  \int_{0}^{T}\left\langle A\left(
t,\upsilon_{t}\right)  -A\left(  t,\upsilon_{t}^{\prime}\right)
,\widehat{\upsilon}_{t}\right\rangle dt\right]
\nonumber \\
&& =-\, \delta\, \mathbb{E} \, \left[  \left\langle R\widehat{y}_{T},\widehat{h}\left(
\overline{y}_{T}\right)  -R\widehat{\overline{y}}_{T}\right\rangle \right]
+\mathbb{E}\, \left[  \left\langle \alpha_{0}R\widehat{\Psi}\left(  Y_{0}\right)
+\delta R\widehat{\Psi}\left(  \overline{Y}_{0}\right)  ,\widehat{Y}_{0}\right\rangle
\right]  \nonumber \\
&& \hspace{2cm}+\, \delta\, \mathbb{E} \, \left[  \int_{0}^{T}\left\langle R\widehat{y}_{t}%
,\widehat{f}\left(  t,\overline{\upsilon}_{t}\right)  \right\rangle
dt\right]  +\theta_{1}\delta\, \mathbb{E} \, \left[  \int_{0}^{T}\left\langle
R\widehat{y}_{t},R\widehat{\overline{y}}_{t}\right\rangle dt\right]\nonumber \\
&&  \hspace{2cm}+\, \delta\, \mathbb{E} \, \left[  \int_{0}^{T}\left\langle \widehat{Y}_{t},R\widehat
{b}\left(  t,\overline{\upsilon}_{t}\right)  \right\rangle dt\right]
+\, \delta\, \mathbb{E} \, \left[  \int_{0}^{T}\left\langle R\widehat{\sigma%
}\left(  t,\overline{\upsilon}_{t}\right)  ,\widehat{Z}_{t}\right\rangle dt\right]\nonumber\\
&& \hspace{2cm}+\, \delta\, \mathbb{E} \, \left[  \int_{0}^{T}\left\langle R\widehat{z}_{t}%
,\widehat{g}\left(  t,\overline{\upsilon}_{t}\right)  \right\rangle
dt\right]  +\delta\theta_{1}\, \mathbb{E} \, \left[  \int_{0}^{T}\left\langle
R\widehat{z}_{t},R\widehat{\overline{z}}_{t}\right\rangle dt\right]\nonumber \\
&& \hspace{6cm}+ \,\delta\, \mathbb{E} \, \left[  \int_{0}^{T}\left\langle \left\langle
R\widehat{\varphi}\left(  t,\overline{\upsilon}_{t},\cdot\right)  ,\widehat{k}_{t
}\right\rangle \right\rangle dt\right]  .
\end{eqnarray}

From assumption (A2) on $h$\ and (A1) we obtain%
\begin{eqnarray}\label{eq:3.8}
&& \hspace{-1.25cm}\alpha_{0}\, \mathbb{E} \, \left[  \langle R\widehat{y}_{T},\widehat{h}\left(  y_{T}\right)\rangle  \right] =\mathbb{E}\, \left[  \langle R\left(  y_{T}-y_{T}^{\prime}\right)  ,\alpha_{0}\left(
h\left(  y_{T}\right)  -h\left(  y_{T}^{\prime}\right)  \right)\rangle  \right]
\nonumber \\ &&
\hspace{7.5cm} \geq\alpha_{0}\, \beta_{1}\, \mathbb{E} \, \left[  \left\vert R\widehat{y}_{T}\right\vert
^{2}\right]  ,
\end{eqnarray}
and
\begin{eqnarray}\label{eq:3.9}
&& \hspace{-0.5cm}-\alpha_{0}\, \mathbb{E} \, \left[  \int_{0}^{T}\left\langle A\left(
s,\upsilon_{s}\right)  -A\left(  s,\upsilon_{s}^{\prime}\right)
,\widehat{\upsilon}_{s}\right\rangle ds\right]
\geq\alpha_{0}\theta_{1}\, \mathbb{E} \, \left[  \int_{0}^{T}\left(  \left\vert
R\widehat{y}_{t}\right\vert ^{2}+\left\Vert R\widehat{z}_{t}\right\Vert^{2}\right)
dt\right]  \nonumber\\
&& \hspace{2cm}+\, \alpha_{0}\, \theta_{2}\, \mathbb{E} \, \left[  \int_{0}^{T}\left(  \left\vert
R^{\ast}\widehat{Y}_{t}\right\vert ^{2}+\left\Vert R^{\ast}\widehat{Z}_{t}\right\Vert
^{2}+||| R^{\ast}\widehat{k}_{t}|||^{2}\right)  dt\right]  .
\end{eqnarray}
Applying these two inequalities (\ref{eq:3.8}) and (\ref{eq:3.9}) in (\ref{eq:3.7}) gives
\begin{eqnarray}\label{eq:3.10}
&& \hspace{-0.5cm}\alpha_{0}\beta_{1}\, \mathbb{E} \, \left[  \left\vert R\widehat{y}_{T}\right\vert
^{2}\right]  +\left(  1-\alpha_{0}\right)  \, \mathbb{E} \, \left[  \left\vert
R\widehat{y}_{T}\right\vert ^{2}\right]  +\left(  1-\alpha_{0}\right)  \theta
_{1}\, \mathbb{E} \, \left[  \int_{0}^{T}\left(  \left\vert R\widehat{y}_{t}\right\vert
^{2}+\left\Vert R\widehat{z}_{t}\right\Vert^{2}\right)  dt\right]\nonumber \\
&& \hspace{1in}+\, \alpha_{0}\theta_{1}\, \mathbb{E} \, \left[  \int_{0}^{T}\left(  \left\vert
R\widehat{y}_{t}\right\vert ^{2}+\left\Vert R\widehat{z}_{t}\right\Vert^{2}\right)
dt\right]  \nonumber \\
&& \hspace{1in}+\, \alpha_{0}\theta_{2}\, \mathbb{E} \, \left[  \int_{0}^{T}\left(
\left\vert R^{\ast}\widehat{Y}_{t}\right\vert ^{2}+\left\Vert R^{\ast}\widehat{Z}
_{t}\right\Vert^{2}+||| R^{\ast}\widehat{k}_{t}|||^{2}\right)  dt\right]
\nonumber \\
&&  \leq-\delta\, \mathbb{E} \, \left[  \left\langle R\widehat{y}_{T},\widehat{h}\left(
\overline{y}_{T}\right)  -R\widehat{\overline{y}}_{T}\right\rangle \right]
+\mathbb{E}\, \left[  \left\langle \alpha_{0}R\widehat{\Psi}\left(  Y_{0}\right)
+\delta R\widehat{\Psi}\left(  \bar{Y}_{0}\right)  ,\widehat{Y}_{0}\right\rangle
\right]\nonumber \\
&& \hspace{1in}+\, \delta\, \mathbb{E} \, \left[  \int_{0}^{T}\left\langle R\widehat{y}_{t}%
,\widehat{f}\left(  t,\overline{\upsilon}_{t}\right)  \right\rangle
dt\right]  +\theta_{1}\delta\, \mathbb{E} \, \left[  \int_{0}^{T}\left\langle
R\widehat{y}_{t},R\widehat{\overline{y}}_{t}\right\rangle dt\right]\nonumber \\
&&  \hspace{1in}+\, \delta\, \mathbb{E} \, \left[  \int_{0}^{T}\left\langle R\widehat
{b}\left(  t,\overline{\upsilon}_{t}\right) ,\widehat{Y}_{t} \right\rangle dt\right]
+\delta\, \mathbb{E} \, \left[  \int_{0}^{T}\left\langle R\widehat{\sigma%
}\left(  t,\overline{\upsilon}_{t}\right)  ,\widehat{Z}_{t}\right\rangle dt\right]\nonumber\\
&& \hspace{1in}+\, \delta\, \mathbb{E} \, \left[  \int_{0}^{T}\left\langle R\widehat{z}_{t}%
,\widehat{g}\left(  t,\overline{\upsilon}_{t}\right)  \right\rangle
dt\right]  +\delta\theta_{1}\, \mathbb{E} \, \left[  \int_{0}^{T}\left\langle
R\widehat{z}_{t},R\widehat{\overline{z}}_{t}\right\rangle dt\right]\nonumber \\
&&   \hspace{2.5in}+\, \delta\, \mathbb{E} \, \left[  \int_{0}^{T}\left\langle \left\langle
R\widehat{\varphi}\left(  t,\overline{\upsilon}_{t},\cdot\right)  ,\widehat{k}_{t%
}\right\rangle \right\rangle dt\right]  .
\end{eqnarray}

Using (A4) and (A2) we have
\begin{eqnarray*}
&&  \hspace{-1.cm} -\delta\, \mathbb{E} \, \left[  \left\langle R\widehat{y}_{T},\widehat{h}\left(
\overline{y}_{T}\right) -R\widehat{\overline{y}}_{T}\right\rangle \right]   =\,-\delta\, \mathbb{E} \, \left[  \left\langle R\widehat{y}_{T},\widehat{h}\left(
\overline{y}_{T}\right)  \right\rangle \right]  +\delta\, \mathbb{E} \, \left[
\left\langle R\widehat{y}_{T},R\widehat{\overline{y}}_{T}\right\rangle \right]  \\
&&  \leq\frac{\delta}{2}\,\left(  \mathbb{E}\, \left[  \left\vert R\widehat{y}_{T}\right\vert
^{2}\right]  +\mathbb{E}\, \left[  \left\vert \widehat{h}\left(  \overline{y}
_{T}\right)  \right\vert ^{2}\right]  \right)  +\frac{\delta}{2}\,\left(  \mathbb{E}\, %
\left[  \left\vert R\widehat{y}_{T}\right\vert ^{2}\right]  +\mathbb{E}\, \left[
\left\vert R\widehat{\overline{y}}_{T}\right\vert ^{2}\right]  \right)  \\
&&
\leq\,\delta\, \mathbb{E}\, \left[  \left\vert R\widehat{y}_{T}\right\vert
^{2}\right]  +\delta\,\left(\frac{1}{2}+\frac{c^{2}}{2} \right)\,  \mathbb{E}\, \left[ \left\vert R\widehat{\overline{y}
}_{T}\right\vert ^{2}\right] ,
\end{eqnarray*}
and
\begin{eqnarray*}
&&  \hspace{-1.25cm} \alpha_{0}\, \mathbb{E} \, \left[  \left\langle R\widehat{\Psi}\left(
Y_{0}\right)  ,\widehat{Y}_{0} \right\rangle \right]  +\delta\, \mathbb{E} \, \left[
\left\langle R\widehat{\Psi}\left(  \overline{Y}_{0}\right)  ,\widehat{Y}_{0}\right\rangle \right]  \\
&& \hspace{2.75cm}\leq-\alpha_{0}\beta_{2}\, \mathbb{E} \, \left[  \left\vert \widehat{Y}_{0}\right\vert
^{2}\right]  +\delta\, \mathbb{E} \, \left[  \left\vert R\widehat{\Psi}\left(
\overline{Y}_{0}\right)  \right\vert \cdot\left\vert \widehat{Y}_{0}\right\vert
\right]  \\
&& \hspace{2.75cm}\leq- \alpha_{0}\, \beta_{2}\, \mathbb{E} \, \left[  \left\vert \widehat{Y}_{0}\right\vert
^{2}\right]  +\delta \, c\, \mathbb{E} \, \left[  \left\vert \widehat{\overline{Y}}
_{0}\right\vert \cdot\left\vert \widehat{Y}_{0}\right\vert \right]  \\
&& \hspace{2.75cm}\leq - \alpha_{0}\, \beta_{2}\, \mathbb{E} \, \left[  \left\vert \widehat{Y}_{0}\right\vert
^{2}\right]  +\frac{\delta}{2} \, c\,\left(  \mathbb{E}\, \left[  \left\vert \widehat{\overline{Y}
}_{0}\right\vert ^{2}\right]  + \mathbb{E}\, \left[  \left\vert \widehat{Y}%
_{0}\right\vert ^{2}\right]  \right)  .
\end{eqnarray*}
Thus, by making use of (A4) and Remark 3.2, inequality (\ref{eq:3.10}) becomes%
\begin{eqnarray}\label{eq:3.11}
&& \hspace{-0.5cm}\left(  1-\alpha_{0}+\alpha_{0}\beta_{1}\right)  \, \mathbb{E} \, \left[
\left\vert R\widehat{y}_{T}\right\vert ^{2}\right]  +\theta_{1}\, \mathbb{E} \, \int
_{0}^{T}\left(  \left\vert R\widehat{y}_{t}\right\vert ^{2}+\left\Vert R\widehat{z}_{t}\right\Vert^{2}\right)  dt \nonumber
\\
&& + \, \alpha_{0}\, \theta_{2}\, \mathbb{E} \, \left[  \int_{0}^{T}\left(  \left\vert
R^{\ast}\widehat{Y}_{t}\right\vert ^{2}+\left\Vert R^{\ast}\widehat{Z}_{t}\right\Vert
^{2}+||| R^{\ast}\widehat{k}_{t}||| ^{2}\right)  dt\right]+\alpha_{0}\beta_{2}\, \mathbb{E} \, \left[  \left\vert \widehat{Y}_{0}\right\vert
^{2}\right]
\nonumber \\ &&
\leq  \delta \, \left(  \mathbb{E}\, \left[  \left\vert R\widehat{y}_{T}\right\vert
^{2}\right]  + C\, \mathbb{E} \, \left[  \left\vert R\widehat{\overline{y}}_{T}\right\vert
^{2}\right]  \right) + \delta \, C \left(\mathbb{E} \, \left[  \left\vert \widehat{\overline{Y}}%
_{0}\right\vert ^{2}\right]  + \mathbb{E} \, \left[  \left\vert \widehat{Y}_{0}\right\vert ^{2}\right] \right) \nonumber \\
&& \hspace{6cm}+ \, \delta \, C\, \mathbb{E} \, \int_{0}^{T}\left(  \left\Vert \widehat{\upsilon}%
_{t}\right\Vert^{2}+\left\Vert \widehat{\overline{\upsilon}}_{t}\right\Vert^{2}\right)  dt,
\end{eqnarray}
with some constant $C>0,$ which from here on can be different from place to place and depends at most on the Lipschitz constants $c,\gamma,\gamma',\alpha_{0},\delta,\theta_{1},\beta_{1}, \beta_2, R$ and $T.$

Next, we want to estimate $\mathbb{E} \, \left[  \left\vert \widehat{Y}_{0}\right\vert ^{2}\right]$.  Applying It\^{o}'s formula (Proposition~\ref{Lemma 3.3}) to $\left\vert \widehat{Y}_{s}\right\vert ^{2}$ over
$\left[  t,T\right]  ,$ and taking the expectation afterwards give
\begin{align*}
& \mathbb{E}\, \left[  \left\vert \widehat{Y}_{t}\right\vert ^{2}\right]
+\mathbb{E}\, \left[  \int_{t}^{T}\left\Vert \widehat{Z}_{s}\right\Vert
^{2}ds\right]  +\mathbb{E}\, \left[  \int_{t}^{T}|||  \widehat
{k}_{s}|||^{2}ds\right] \\
& =\mathbb{E}\, \left[  \left\vert \alpha_{0}\widehat{h}\left(  y_{T}\right)
+\left(  1-\alpha_{0}\right)  R\widehat{y}_{T}+\delta\left(  \widehat{h}\left(
\overline{y}_{T}\right)  -R\widehat{\overline{y}}_{T}\right)  \right\vert ^{2}\right] \\
& \hspace{0.5cm}-2\, \mathbb{E}\, \left[\int_{t}^{T}\left\langle \widehat{Y}_{s},\alpha_{0}\widehat
{f}\left(  s,v_{s}\right)  -\left(  1-\alpha_{0}\right)  \theta
_{1}R\widehat{y}_{s}+\delta\left(  \widehat{f}\left(  s,\overline{\upsilon}
_{s}\right)  +\theta_{1}R\widehat{\overline{y}}_{s}\right)  \right\rangle ds\right]\\
& \hspace{0.5cm}+\, \mathbb{E}\, \left[\int_{t}^{T}\left\vert \alpha_{0}\widehat{g}\left(
s,\upsilon_{s}\right)  -\left(  1-\alpha_{0}\right)  \theta_{1}R\widehat{z}_{s}%
+\delta\left(  \widehat{g}\left(  s,\overline{\upsilon}_{s}\right)  +\theta
_{1}R\widehat{\overline{z}}_{s}\right)  \right\vert ^{2}ds\right].
\end{align*}
Therefore
\begin{align*}
& \mathbb{E}\, \left[  \left\vert \widehat{Y}_{t}\right\vert ^{2}\right]
+\mathbb{E}\, \left[  \int_{t}^{T}\left\Vert \widehat{Z}_{s}\right\Vert
^{2}ds\right]  +\mathbb{E}\, \left[  \int_{t}^{T}||| \widehat
{k}_{s}||| ^{2}ds\right]  \\
& \leq 4\,\mathbb{E}\, \left[  \alpha_{0}^{2}\left\vert \widehat{h}\left(
y_{T}\right)  \right\vert ^{2}+\left(  1-\alpha_{0}\right)  ^{2}\left\vert
R\widehat{y}_{T}\right\vert ^{2}+\delta^{2}\left\vert \widehat{h}\left(  \overline{y}_{T}\right)  \right\vert ^{2}+\delta^{2}\left\vert R\widehat{\overline{y}}%
_{T}\right\vert ^{2}\right]  \\
& +2\,\mathbb{E}\, \left[  \int_{t}^{T}\left\vert \widehat{Y}_{s}\right\vert
\cdot\left(  \alpha_{0}\left\vert \widehat{f}\left(  s,\upsilon
_{s}\right)  \right\vert +\left(  1-\alpha_{0}\right)  \theta_{1}\left\vert
R\widehat{y}_{s}\right\vert +\delta\left\vert \widehat{f}\left(
s,\overline{\upsilon}_{s}\right)  \right\vert +\delta\theta_{1}\left\vert
R\widehat{\overline{y}}_{s}\right\vert \right)  ds\right]  \\
& +\mathbb{E}\, \left[  \int_{t}^{T}\left(  \left(  1+\varepsilon\right)
\alpha_{0}^{2}\left\vert \widehat{g}\left(  s,\upsilon_{s}\right)
\right\vert ^{2} \right. \right. \\
& \hspace{3.5cm} \left. \left.+\, \left(  1+\frac{1}{\varepsilon}\right)  \left\vert \left(
1-\alpha_{0}\right)  \theta_{1}R\widehat{z}_{s}+\delta\left(  \widehat{g%
}\left(  s,\overline{\upsilon}_{s}\right)  +\theta_{1}R\widehat{\overline{z}}%
_{s}\right)  \right\vert ^{2}\right)  ds\right],
\end{align*}
for any $\varepsilon >0.$ Choose then \ $\varepsilon=\frac{1-\gamma}{2\gamma}$ to get%
\begin{eqnarray}\label{eq:3.12}
&& \hspace{-0,5cm}\mathbb{E}\, \left[  \left\vert \widehat{Y}_{t}\right\vert ^{2}\right]
+\mathbb{E}\, \left[  \int_{t}^{T}\left\Vert \widehat{Z}_{s}\right\Vert
^{2}ds\right]  +\mathbb{E}\, \left[  \int_{t}^{T}||| \widehat
{k}_{s}||| ^{2}ds\right]\nonumber \\ && \hspace{-0,5cm}
\leq 4\, \mathbb{E}\, \left[  \alpha_{0}^{2}\left\vert \widehat{h}\left(
y_{T}\right)  \right\vert ^{2}+\left(  1-\alpha_{0}\right)  ^{2}\left\vert
R\widehat{y}_{T}\right\vert ^{2}+\delta^{2}\left\vert \widehat{h}\left(  \overline{y}_{T}\right)  \right\vert ^{2}+\delta^{2}\left\vert R\widehat{\overline{y}}%
_{T}\right\vert ^{2}\right]
\nonumber \\ && \hspace{-0,5cm}
+ \, 2\, \mathbb{E} \, \left[  \int_{t}^{T}\left\vert \widehat{Y}_{s}\right\vert
\cdot\left(  \alpha_{0}\left\vert \widehat{f}\left(  s,\upsilon
_{s}\right)  \right\vert +\left(  1-\alpha_{0}\right)  \theta_{1}\left\vert
R\widehat{y}_{s}\right\vert +\delta\left\vert \widehat{f}\left(
s,\overline{\upsilon}_{s}\right)  \right\vert + \delta\theta_{1}\left\vert
R\widehat{\overline{y}}_{s}\right\vert \right)  ds\right]  \nonumber\\
&& \hspace{4cm} +\, \mathbb{E} \, \left[  \int_{t}^{T}\left(\frac{1+\gamma}{2\gamma}\right)\alpha_{0}%
^{2}\left\vert \widehat{g}\left(  s,\upsilon_{s}\right)  \right\vert
^{2}ds\right]\nonumber \\ \hspace{-0,5cm}
&& + \, 3\, \mathbb{E} \, \left[  \int_{t}^{T}\left(  \frac{1+\gamma}{1-\gamma
}\right)  \left(  \left(  1-\alpha_{0}\right)  ^{2}\theta_{1}^{2}\left\Vert R\widehat
{z}_{s}\right\Vert^{2}+\delta^{2}\left\vert \widehat{g}\left(
s,\overline{\upsilon}_{s}\right)  \right\vert ^{2}+\delta^{2}\theta_{1}%
^{2}\left\Vert R\widehat{\overline{z}}_{s}\right\Vert^{2}\right)  ds\right]\nonumber \\
&& \hspace{-0,5cm} =: I_{1}+I_{2}+I_{3}+I_{4},
\end{eqnarray}
where $I_{i}$ is the quantity in term $i$ of the right hand side of this inequality for each $i=1,2,3,4.$

It follows by the fact that $h$ is Lipschitz that%
\begin{eqnarray}\label{eq:3.13}
&& \hspace{-0.5cm} I_{1}   =4 \, \mathbb{E}\, \left[  \alpha_{0}^{2}\left\vert \widehat{h}\left(
y_{T}\right)  \right\vert ^{2}+\left(  1-\alpha_{0}\right)  ^{2}\left\vert
R\widehat{y}_{T}\right\vert ^{2}+\delta^{2}\left\vert \widehat{h}\left(  \overline{y}_{T}\right)  \right\vert ^{2}+\delta^{2}\left\vert R\widehat{\overline{y}}%
_{T}\right\vert ^{2}\right]\nonumber \\
&& \leq C\, \mathbb{E} \, \left[  \left\vert \widehat{y}_{T}\right\vert ^{2}%
+\delta\left\vert \widehat{\overline{y}}_{T}\right\vert ^{2}\right].
\end{eqnarray}
Also%
\begin{eqnarray*}
&& \hspace{-0.5cm} I_{2}= 2 \, \mathbb{E}\, \left[  \int_{t}^{T}\left\vert \widehat{Y}_{s}\right\vert
\cdot\left(  \alpha_{0}\left\vert \widehat{f}\left(  s,\upsilon
_{s}\right)  \right\vert +\left(  1-\alpha_{0}\right)  \theta_{1}\left\vert
R\widehat{y}_{s}\right\vert +\delta\left\vert \widehat{f}\left(
s,\overline{\upsilon}_{s}\right)  \right\vert +\delta\theta_{1}\left\vert
R\widehat{\overline{y}}_{s}\right\vert \right)  ds\right]  \\
&& \hspace{-0,5cm} \leq\mathbb{E}\, \left[  \int_{t}^{T}\left(  \left(  \frac{8\;c\;\alpha_{0}^{2}%
}{1-\gamma}\right)  \left\vert \widehat{Y}_{s}\right\vert ^{2}+\left(\frac{1-\gamma}%
{8c}\right)\left\vert \widehat{f}\left(  s,\upsilon_{s}\right)  \right\vert
^{2}\right)  +\left(  1-\alpha_{0}\right)  \theta_{1}\left(  \left\vert
\widehat{Y}_{s}\right\vert ^{2}+\left\vert R\widehat{y}_{s}\right\vert ^{2}\right)
\right.  \\
&& \left.  \hspace{2cm} + \, \delta\left(  \left\vert \widehat{Y}_{s}\right\vert ^{2}+\left\vert
\widehat{f}\left(  s,\overline{\upsilon}_{s}\right)  \right\vert
^{2}\right)  +\delta\theta_{1}\left(  \left\vert \widehat{Y}_{s}\right\vert
^{2}+\left\vert R\widehat{\overline{y}}_{s}\right\vert ^{2}\right)  )ds\right]  \\
&& \hspace{-0,5cm} \leq\mathbb{E}\, \left[  \int_{t}^{T}\left(  \left(  \left(  \frac{8\;c\;\alpha
_{0}^{2}}{1-\gamma}\right)  \left\vert \widehat{Y}_{s}\right\vert ^{2}+\left(\frac{1-\gamma
}{8c}\right)c\left\Vert \widehat{\upsilon}_{s}\right\Vert %
^{2}\right)  +\left(  1-\alpha_{0}\right)  \theta_{1}\left(  \left\vert
\widehat{Y}_{s}\right\vert ^{2}+\left\vert R\widehat{y}_{s}\right\vert ^{2}\right)
\right.  \right.  \\
&& \left.  \left. \hspace{3.5cm} +\, \delta \left(  \left\vert \widehat{Y}_{s}\right\vert
^{2}+c\left\Vert \widehat{\overline{\upsilon}}_{s}\right\Vert %
^{2}\right)  +\delta\theta_{1}\left(  \left\vert \widehat{Y}_{s}\right\vert
^{2}+\left\vert R\widehat{\overline{y}}_{s}\right\vert ^{2}\right)  \right)
ds\right]  .
\end{eqnarray*}
Hence, since $\left\Vert \upsilon_{s}\right\Vert %
^{2}=\left\vert y_{s}\right\vert %
^{2}+\left\vert Y_{s}\right\vert %
^{2}+\left\Vert z_{s}\right\Vert %
^{2}+\left\Vert Z_{s}\right\Vert %
^{2} +||| k_{s}||| %
^{2},$ we have
\begin{eqnarray}\label{eq:3.14}
&& \hspace{-0.5cm} I_{2}   \leq C\, \mathbb{E} \, \left[  \int_{t}^{T}\left\vert \widehat{Y}%
_{s}\right\vert ^{2}ds\right]  +C\mathbb{E}\, \left[  \int_{t}^{T}\left\vert
\widehat{y}_{s}\right\vert ^{2}ds\right]  + C\, \delta\, \mathbb{E} \, \left[  \int_{t}%
^{T}\left\vert \widehat{\overline{y}}_{s}\right\vert ^{2}ds\right]  \nonumber \\ &&
\hspace{1cm} + \, C \, \mathbb{E} \, \left[  \int_{t}^{T}\left\Vert \widehat{z}_{s}\right\Vert^{2}ds\right]\nonumber + \; \left(\frac{1-\gamma}{8}\right)\, \mathbb{E} \, \left[  \int_{t}^{T}\left(  \left\Vert \widehat
{Z}_{s}\right\Vert^{2}+||| \widehat{k}_{s}||| ^{2}\right)  ds\right]  \\ &&
\hspace{7cm}  + \, C\, \delta\, \mathbb{E}\, \left[  \int_{t}%
^{T}\left\Vert \widehat{\overline{\upsilon}}_{s}\right\Vert %
^{2}ds\right]  .
\end{eqnarray}

For $I_{3}$ we apply (A4) to derive
\begin{eqnarray}\label{eq:3.15}
&& \hspace{-0.5cm} I_{3}   =\mathbb{E}\, \left[  \int_{t}^{T}\left(  \frac{1+\gamma}{2\gamma
}\right)  \alpha_{0}^{2}\left\vert \widehat{g}\left(  s,\upsilon
_{s}\right)  \right\vert ^{2}ds\right] \nonumber \\ &&
\leq\mathbb{E}\, \left[  \int_{t}^{T}\left(  c\left(  \frac{1+\gamma}{2\gamma
}\right)  \alpha_{0}^{2}\left\vert \widehat{y}_{s}\right\vert ^{2}+c\left(
\frac{1+\gamma}{2\gamma}\right) \alpha_{0}^{2}\left\vert \widehat{Y}%
_{s}\right\vert ^{2}+c\left(  \frac{1+\gamma}{2\gamma}\right)  \alpha_{0}%
^{2}\left\Vert \widehat{z}_{s}\right\Vert^{2}\right.  \right.  \nonumber\\
&& \hspace{3cm} + \left.  \left.  \left(  \frac{1+\gamma}{2\gamma}\right)  \alpha_{0}%
^{2}\gamma\left\Vert \widehat{Z}_{s}\right\Vert^{2}+\left(  \frac{1+\gamma
}{2\gamma}\right)  \alpha_{0}^{2}\gamma||| \widehat{k}%
_{s}||| ^{2}\right)  ds\right]\nonumber \\
&& \leq C\, \mathbb{E} \, \left[  \int_{t}^{T}\left(  \left\vert \widehat{y}%
_{s}\right\vert ^{2}+\left\vert \widehat{Y}_{s}\right\vert ^{2}+\left\Vert \widehat
{z}_{s}\right\Vert^{2}\right)  ds\right] \nonumber  \\ &&
\hspace{3cm} +\left(
\frac{1+\gamma}{2}\right)  \alpha_{0}^{2}\,\mathbb{E}\, \left[ \int_{t}^{T} \left(  \left\Vert \widehat{Z}%
_{s}\right\Vert^{2}+||| \widehat{k}%
_{s}|||^{2}\right)  ds\right].
\end{eqnarray}

Finally,
\begin{eqnarray}\label{eq:3.16}
&& \hspace{-0.5cm}  I_{4}   =3\,\mathbb{E}\, \left[  \int_{t}^{T}\left(  \frac{1+\gamma
}{1-\gamma}\right) \left(   \left(  1-\alpha_{0}\right)  ^{2}\theta_{1}^{2}\left\Vert
R\widehat{z}_{s}\right\Vert^{2}+\delta^{2}\left\vert \widehat{g}\left(
s,\overline{\upsilon}_{s}\right)  \right\vert ^{2}+\delta^{2}\theta_{1}%
^{2}\left\Vert R\widehat{\overline{z}}_{s}\right\Vert^{2}\right)  ds\right]\nonumber \\
&& \leq3\,\mathbb{E}\, \left[  \int_{t}^{T}\left(  \left(  \frac{1+\gamma}{1-\gamma
}\right)  \left(  1-\alpha_{0}\right)  ^{2}\theta_{1}^{2}\left\Vert R\widehat
{z}_{s}\right\Vert^{2}+\delta^{2}\left(  \frac{1+\gamma}{1-\gamma}\right)
c\left\vert \widehat{\overline{y}}_{s}\right\vert ^{2}  \right.  \right. \nonumber \\ &&
\hspace{1cm} + \, \delta^{2}\left(  \frac{1+\gamma}{1-\gamma}\right)  c\left\vert
\widehat{\overline{Y}}_{s}\right\vert ^{2}
+\, \delta^{2}\left(  \frac{1+\gamma
}{1-\gamma}\right)  c\left\Vert \widehat{\overline{z}}_{s}\right\Vert^{2}%
+\delta^{2}\left(  \frac{1+\gamma}{1-\gamma}\right)  \gamma\left\Vert
\widehat{\overline{Z}}_{s}\right\Vert^{2} \nonumber \\ && \hspace{2.5cm} +\, \delta^{2}\left(  \frac{1+\gamma}{1-\gamma}\right)  \gamma \, ||| \widehat{\overline{k}}_{s}||| ^{2}+\delta^{2}%
\theta_{1}^{2}\left\Vert R\widehat{\overline{z}}_{s}\right\Vert^{2}\nonumber\\
&& \leq C \, \mathbb{E} \, \left[  \int_{t}^{T}\left(  \left\Vert R\widehat{z}%
_{s}\right\Vert^{2}+\delta\left\vert \widehat{\overline{y}}_{s}\right\vert
^{2}+\delta\left\Vert \widehat{\overline{z}}_{s}\right\Vert^{2}+\delta\left\vert
\widehat{\overline{Y}}_{s}\right\vert ^{2}\right)  ds\right]\nonumber\\
&& \hspace{3.5cm} +\, \delta \, C\,\left(  \frac{1+\gamma}{1-\gamma}\right)  \gamma\, \mathbb{E} \, \left[
\int_{t}^{T}\left(  \left\Vert \widehat{\overline{Z}}_{s}\right\Vert^{2}%
+||| \widehat{\overline{k}}_{s}||| ^{2}\right)  ds\right]  .
\end{eqnarray}

Now we substitute (\ref{eq:3.13})-(\ref{eq:3.16}) in (\ref{eq:3.12}) to find that%
\begin{eqnarray*}
&&  \mathbb{E}\, \left[  \left\vert \widehat{Y}_{t}\right\vert ^{2}\right]
+\mathbb{E}\, \left[  \int_{t}^{T}\left(  \left\Vert \widehat{Z}_{s}\right\Vert
^{2}+||| \widehat{k}_{s}||| ^{2}\right)
ds\right]  \\
&&  \leq C\, \mathbb{E} \, \left[  \left\vert \widehat{y}_{T}\right\vert ^{2}%
+\delta^{2}\left\vert \widehat{\overline{y}}_{T}\right\vert ^{2}\right]
+C\, \mathbb{E} \, \left[  \int_{t}^{T}\left\vert \widehat{Y}_{s}\right\vert
^{2}ds\right] +C\, \mathbb{E} \, \left[  \int_{t}^{T}\left\vert \widehat{y}%
_{s}\right\vert ^{2}ds\right]  \\ && +\, C \, \delta\, \mathbb{E} \, \left[  \int_{t}^{T}\left\vert \widehat{\overline{y}}%
_{s}\right\vert ^{2}ds\right] + C\, \mathbb{E} \, \left[  \int_{t}^{T}\left\Vert
\widehat{z}_{s}\right\Vert^{2}ds\right]  +C\delta\, \mathbb{E} \, \left[  \int_{t}%
^{T}\left\Vert \widehat{\overline{\upsilon}}_{s}\right\Vert %
^{2}ds\right]   \\
&&  \hspace{2cm} +\, \left(\frac{1-\gamma}{8}\right)\, \mathbb{E} \, \left[  \int_{t}^{T}\left(  \left\Vert \widehat
{Z}_{s}\right\Vert^{2}+||| \widehat{k}_{s}||| ^{2}\right)  ds\right]
\\ && \hspace{2cm}
 + \, C\, \mathbb{E} \, \left[  \int_{t}%
^{T}\left(  \left\vert \widehat{y}_{s}\right\vert ^{2}+\left\vert \widehat{Y}%
_{s}\right\vert ^{2} +\left\Vert \widehat{z}_{s}\right\Vert^{2}\right)  ds\right]
  \\
&& \hspace{2cm}
+ \left(  \frac{1+\gamma}{2}\right)  \alpha_{0}^{2}\,\mathbb{E}\, \left[  \int
_{t}^{T}\left(  \left\Vert \widehat{Z}_{s}\right\Vert^{2}+||| \widehat{k}_{s}|||^{2}\right)  ds\right]  \\
&& \hspace{2cm}+\, C\, \mathbb{E} \, \left[  \int_{t}^{T}\left(  \left\Vert R\widehat{z}%
_{s}\right\Vert^{2}+\delta\left\vert \widehat{\overline{y}}_{s}\right\vert
^{2}+\delta\left\Vert \widehat{\overline{z}}_{s}\right\Vert^{2}+\delta\left\vert
\widehat{\overline{Y}}_{s}\right\vert ^{2}\right)  ds\right]  \\
&&  \hspace{2cm}+\, \delta\, C\,\left(  \frac{1+\gamma}{1-\gamma}\right)  \gamma~ \mathbb{E} \, \left[
\int_{t}^{T}\left(  \left\Vert \widehat{\overline{Z}}_{s}\right\Vert^{2}%
+||| \widehat{\overline{k}}_{s}||| ^{2}\right)  ds\right]  .
\end{eqnarray*}
This implies
\begin{align*}
& \hspace{-0.5cm} \mathbb{E}\, \left[  \left\vert \widehat{Y}_{t}\right\vert ^{2}\right]
+\mathbb{E}\, \left[  \int_{t}^{T}\left(  \left\Vert \widehat{Z}_{s}\right\Vert
^{2}+|||  \widehat{k}_{s}||| ^{2}\right)
ds\right]  \\
& \leq C\, \mathbb{E} \, \left[  \left\vert \widehat{y}_{T}\right\vert ^{2}\right]
+\delta C\, \mathbb{E} \, \left[  \left\vert \widehat{\overline{y}}_{T}\right\vert
^{2}\right]  +C\, \mathbb{E} \, \left[  \int_{t}^{T}\left\vert \widehat{Y}%
_{s}\right\vert ^{2}ds\right]  +C\delta\, \mathbb{E} \, \left[  \int_{t}%
^{T}\left\Vert \widehat{\overline{\upsilon}}_{s}\right\Vert %
^{2}ds\right]  \\
& \hspace{2cm}+\, C\, \mathbb{E} \, \left[  \int_{t}^{T}\left(  \left\vert \hat{y}%
_{s}\right\vert ^{2}+\left\Vert \hat{z}_{s}\right\Vert^{2}\right)  ds\right]
\\ & \hspace{2cm} +\, \left(  \left(\frac{1-\gamma}{8}\right)+\left(  \frac{1+\gamma}{2}\right)  \alpha_{0}%
^{2}\right)  \mathbb{E}\, \left[  \int_{t}^{T}\left(  \left\Vert \widehat{Z}%
_{s}\right\Vert^{2}+||| \widehat{k}_{s}||| ^{2}\right)  ds\right]  \\
& \leq C\left(  \mathbb{E}\, \left[  \left\vert \widehat{y}_{T}\right\vert
^{2}\right]  +\delta\, \mathbb{E} \, \left[  \left\vert \widehat{\overline{y}}%
_{T}\right\vert ^{2}\right]  \right)  +C\, \mathbb{E} \, \left[  \int_{t}^{T}\left(
\left\vert \widehat{y}_{s}\right\vert ^{2}+\left\Vert \widehat{z}_{s}\right\Vert
^{2}+\left\vert \widehat{Y}_{s}\right\vert ^{2}\right)  ds\right]  \\
& \hspace{1cm} +\, C\delta\, \mathbb{E} \, \left[  \int_{t}^{T}\left\Vert \widehat{\overline{\upsilon}
}_{s}\right\Vert^{2}ds\right]  +\left(  \frac{5+3\gamma}%
{8}\right)  \mathbb{E}\, \left[  \int_{t}^{T}\left(  \left\Vert \widehat{Z}%
_{s}\right\Vert^{2}+||| \widehat{k}_{s}||| ^{2}\right)  ds\right]  ,
\end{align*}
since $\alpha_{0}<1$\ and $\frac{1-\gamma}{8}+\frac{1+\gamma}{2}%
=\frac{5+3\gamma}{8}.$ Hence there exists an universal constant $C>0$ such that
\begin{eqnarray}\label{eq:3.17}
&& \hspace{-1.5cm}\mathbb{E}\, \left[  \left\vert \widehat{Y}_{t}\right\vert ^{2}\right]  +\left(
1-\left(\frac{5+3\gamma}{8}\right)\right)  \mathbb{E}\, \left[  \int_{t}^{T}\left(  \left\Vert
\widehat{Z}_{s}\right\Vert^{2}+||| \widehat{k}_{s}||| ^{2}\right)  ds\right] \nonumber\\ &&
\leq C\, \mathbb{E} \, \left[  \int_{t}^{T}\left\vert \widehat{Y}_{s}\right\vert
^{2}ds\right]  +C\left(  \mathbb{E}\, \left[  \left\vert \widehat{y}_{T}\right\vert
^{2}\right]  +\delta\, \mathbb{E} \, \left[  \left\vert \widehat{\overline{y}}%
_{T}\right\vert ^{2}\right]  \right) \nonumber\\ &&
 \hspace{2cm}+ \, C\, \mathbb{E} \, \left[  \int_{t}^{T}\left(
\left\vert \widehat{y}_{s}\right\vert ^{2}+\left\Vert \widehat{z}_{s}\right\Vert
^{2}+\delta\left\Vert \widehat{\overline{\upsilon}}_{s}\right\Vert ^{2}\right)  ds\right]  .
\end{eqnarray}
By using Gronwall's inequality we deduce then%
\begin{eqnarray}\label{eq:3.18}
&& \hspace{-1.5cm}\mathbb{E}\, \left[  \left\vert \widehat{Y}_{t}\right\vert ^{2}\right]  \leq
Ce^{T-t}\left(  C\left(  \mathbb{E}\, \left[  \left\vert \widehat{y}_{T}\right\vert
^{2}\right]  +\delta\, \mathbb{E} \, \left[  \left\vert \widehat{\overline{y}}%
_{T}\right\vert ^{2}\right]  \right) \right.\nonumber \\ &&
\hspace{2cm}
\left. +\, C\, \mathbb{E} \, \left[  \int_{t}^{T}\left(
\left\vert \widehat{y}_{s}\right\vert ^{2}+\left\Vert \widehat{z}_{s}\right\Vert
^{2}+\delta\left\Vert \widehat{\overline{\upsilon}}_{s}\right\Vert ^{2}\right)  ds\right]  \right)  ,
\end{eqnarray}
for each $0\leq t\leq T$. Consequently
\begin{eqnarray}\label{eq:3.19}
&& \hspace{-1.5cm} \mathbb{E}\, \left[  \left\vert \widehat{Y}_{0}\right\vert ^{2}\right]  \leq
C^{2}e^{T}\left(  C\left(  \mathbb{E}\, \left[  \left\vert \widehat{y}_{T}\right\vert
^{2}\right]  +\delta\, \mathbb{E} \, \left[  \left\vert \widehat{\overline{y}}%
_{T}\right\vert ^{2}\right]  \right) \right.\nonumber \\ &&
\hspace{2cm} \left. +\, C\, \mathbb{E} \, \left[  \int_{0}^{T}\left(
\left\vert \widehat{y}_{s}\right\vert ^{2}+\left\Vert \widehat{z}_{s}\right\Vert
^{2}+\delta\left\Vert \widehat{\overline{\upsilon}}_{s}\right\Vert ^{2}\right)  ds\right]  \right)  .
\end{eqnarray}

\bigskip

Take the supremum over $t$ in (\ref{eq:3.18}) to obtain%
\begin{eqnarray}\label{eq:3.20}
&& \hspace{-1.5cm} \underset{0\leq t\leq T}{\sup}\, \mathbb{E} \, \left[  \left\vert \widehat{Y}%
_{t}\right\vert ^{2}\right]  \leq Ce^{T}\left(  C\left(  \mathbb{E}\, \left[
\left\vert \widehat{y}_{T}\right\vert ^{2}\right]  +\delta\, \mathbb{E} \, \left[
\left\vert \widehat{\overline{y}}_{T}\right\vert ^{2}\right]  \right)\right.
\nonumber \\ &&
\hspace{2cm} \left. +\, C\, \mathbb{E} \, \left[  \int_{0}^{T}\left(  \left\vert \widehat{y}_{s}\right\vert
^{2}+\left\Vert \widehat{z}_{s}\right\Vert^{2}+\delta\left\Vert \widehat
{\overline{\upsilon}}_{s}\right\Vert^{2}\right)  ds\right]
\right)  .
\end{eqnarray}

Applying this in (\ref{eq:3.17}) gives also
\begin{eqnarray}\label{eq:3.21}
&& \mathbb{E}\, \left[  \int_{0}^{T}\left(  \left\Vert \widehat{Z}_{s}\right\Vert
^{2}+ ||| \widehat{k}_{s}||| ^{2}\right)
ds\right]\nonumber \\ &&
\leq C\left(  T-t\right)  \underset{0\leq s\leq T}{\sup}\, \mathbb{E} \, \left[
\left\vert \widehat{Y}_{s}\right\vert ^{2}\right]  +C\left(  \mathbb{E}\, \left[
\left\vert \widehat{y}_{T}\right\vert ^{2}\right]  +\delta\, \mathbb{E} \, \left[
\left\vert \widehat{\overline{y}}_{T}\right\vert ^{2}\right]  \right)
\nonumber \\ &&
\hspace{2cm} +C\, \mathbb{E} \, \left[  \int_{0}^{T}\left(  \left\vert \widehat{y}_{s}\right\vert
^{2}+\left\Vert \widehat{z}_{s}\right\Vert^{2}+\delta\left\Vert \widehat
{\overline{\upsilon}}_{s}\right\Vert^{2}\right)  ds\right] \nonumber  \\
&& \leq C\left(  \mathbb{E}\, \left[  \left\vert \widehat{y}_{T}\right\vert
^{2}\right]  +\delta\, \mathbb{E} \, \left[  \left\vert \widehat{\overline{y}}%
_{T}\right\vert ^{2}\right]  \right)  +C\, \mathbb{E} \, \left[  \int_{0}%
^{T}\left(  \left\vert \widehat{y}_{s}\right\vert ^{2}+\left\Vert \widehat{z}%
_{s}\right\Vert^{2}+\delta\left\Vert \widehat{\overline{\upsilon}}_{s}\right\Vert
^{2}\right)  ds\right]. \nonumber \\
\end{eqnarray}

From (\ref{eq:3.20}) we get%
\begin{eqnarray}\label{eq:3.22}
&& \hspace{-1cm}\mathbb{E}\, \left[  \int_{0}^{T}\left\vert \widehat{Y}_{t}\right\vert ^{2}dt\right]
\leq TCe^{T}\left(  C\left(  \mathbb{E}\, \left[  \left\vert \widehat{y}%
_{T}\right\vert ^{2}\right]  +\delta\, \mathbb{E} \, \left[  \left\vert \widehat
{\overline{y}}_{T}\right\vert ^{2}\right]  \right) \right. \nonumber \\ &&
\hspace{3.5cm} \left. + \, C\, \mathbb{E} \, \left[  \int
_{0}^{T}\left(  \left\vert \widehat{y}_{s}\right\vert ^{2}+\left\Vert \widehat{z}%
_{s}\right\Vert^{2}+\delta\left\Vert \widehat{\overline{\upsilon}}_{s}\right\Vert
^{2}\right)  ds\right]  \right) .
\end{eqnarray}

Finally, add (\ref{eq:3.19}), (\ref{eq:3.22}) and (\ref{eq:3.21}) together to conclude that%
\begin{eqnarray}\label{eq:3.23}
&& \mathbb{E}\, \left[  \left\vert \widehat{Y}_{0}\right\vert ^{2}\right]
+\mathbb{E}\, \left[  \int_{0}^{T}\left(  \left\vert \widehat{Y}_{s}\right\vert
^{2}+\left\Vert \widehat{Z}_{s}\right\Vert^{2}+||| \widehat{k}_{s}||| ^{2}\right)  ds\right] \nonumber \\
&& \leq C\left(  \mathbb{E}\, \left[  \left\vert \widehat{y}%
_{T}\right\vert ^{2}\right]  +\delta\, \mathbb{E} \, \left[  \left\vert \widehat
{\overline{y}}_{T}\right\vert ^{2}\right]  \right)  +C\, \mathbb{E} \, \left[  \int_{0}^{T}\left(  \left\vert \widehat{y}_{s}\right\vert
^{2}+\left\Vert \widehat{z}_{s}\right\Vert^{2}+\delta\left\Vert \widehat
{\overline{\upsilon}}_{s}\right\Vert^{2}\right)  ds\right]. \nonumber \\
\end{eqnarray}

Let us recall that (\ref{eq:3.11}) can be rewritten as
\begin{eqnarray}\label{eq:3.24}
&& \left(  1-\alpha_{0}+\alpha_{0}\beta_{1}\right)  \mathbb{E}\, \left[
\left\vert R\widehat{y}_{T}\right\vert ^{2}\right]  +\theta_{1}\, \mathbb{E} \, \left[
\int_{0}^{T}\left(  \left\vert R\widehat{y}_{s}\right\vert ^{2}+\left\Vert
R\widehat{z}_{s}\right\Vert^{2}\right)  ds\right]\nonumber \\ &&
\hspace{2.5cm} +\, \theta_{2}\, \mathbb{E} \, \left[  \int_{0}^{T}\left(  \left\vert \widehat{Y}%
_{s}\right\vert ^{2}+\left\Vert \widehat{Z}_{s}\right\Vert^{2}+||| \widehat{k}_{s}||| ^{2}\right)  ds\right]
+\beta_{2}\, \mathbb{E} \, \left[  \left\vert \widehat{Y}_{0}\right\vert ^{2}\right]\nonumber\\
&& \hspace{2cm} \leq\delta \, C\,\left(  \mathbb{E}\, \left[  \left\vert \widehat{\overline{Y}}%
_{0}\right\vert ^{2}\right]  +\mathbb{E}\, \left[  \left\vert \widehat{Y}%
_{0}\right\vert ^{2}\right]  +\mathbb{E}\, \left[  \left\vert R\widehat{y}%
_{T}\right\vert ^{2}\right]  +\mathbb{E}\, \left[  \left\vert R\widehat{\overline{y}
}_{T}\right\vert ^{2}\right]  \right)  \nonumber\\
&& \hspace{5cm} +\, \delta \, C\, \mathbb{E} \, \left[  \int_{0}^{T}\left(  \left\Vert \widehat{\upsilon
}_{s}\right\Vert^{2}+\left\Vert \widehat{\overline{\upsilon}}%
_{s}\right\Vert^{2}\right)  ds\right],
\end{eqnarray}
where we have used in (\ref{eq:3.11}) the fact that if two quantities $\xi_{1},\xi_{2}\geq 0$ satisfy $\xi_{1}+\alpha_{0}\xi_{2}\leq K$ then $\xi_{1}+\xi_{2} \leq C K$ for $C=1+\frac{1}{\alpha_{0}}.$

Notice that we want to study the cases when: I) $\theta_{2}>0,\beta_{2}>0,$\ II) $\theta_{2}=0,\beta
_{2}>0,$\, III) $\theta_{2}>0,\beta_{2}=0,$ IV) $\theta_{2}=0,\beta_{2}=0.$

\medskip
I) $\theta_{2}>0,\beta_{2}>0.$\ Let $\Gamma_{1}=\min\left\{  1-\alpha_{0}+\alpha_{0}%
\beta_{1},\theta_{1},\beta_{2},\theta_{2}\right\}.$ Then $\Gamma_{1}>0,$ and inequality (\ref{eq:3.24}) becomes%
\begin{eqnarray*}
&& \hspace{-1cm} \Gamma_{1}\left(  \mathbb{E}\, \left[  \left\vert R\widehat{y}_{T}\right\vert
^{2}\right]  +\mathbb{E}\, \left[  \int_{0}^{T}\left(  \left\vert R\widehat{y}%
_{s}\right\vert ^{2}+\left\Vert R\widehat{z}_{s}\right\Vert^{2}\right)
ds\right]  +\mathbb{E}\, \left[  \left\vert \widehat{Y}_{0}\right\vert ^{2}\right]
\right.  \\
&& \hspace{4cm} +\left.  \mathbb{E}\, \left[  \int_{0}^{T}\left(  \left\vert \widehat{Y}%
_{s}\right\vert ^{2}+\left\Vert \widehat{Z}_{s}\right\Vert^{2}+||| \widehat{k}_{s}||| ^{2}\right)  ds\right]  \right)
\\
&& \hspace{2cm} \leq\delta C\left(  \mathbb{E}\, \left[  \left\vert \widehat{\overline{Y}}%
_{0}\right\vert ^{2}\right]  +\mathbb{E}\, \left[  \left\vert \widehat{Y}%
_{0}\right\vert ^{2}\right]  +\mathbb{E}\, \left[  \left\vert R\widehat{y}%
_{T}\right\vert ^{2}\right]  +\mathbb{E}\, \left[  \left\vert R\widehat{\overline{y}
}_{T}\right\vert ^{2}\right]  \right)  \\
&& \hspace{5cm} +\, \delta \, C\, \mathbb{E} \, \left[  \int_{0}^{T}\left(  \left\Vert \widehat{\upsilon
}_{s}\right\Vert^{2}+\left\Vert \widehat{\overline{\upsilon}}
_{s}\right\Vert^{2}\right)  ds\right],
\end{eqnarray*}
or in particular
\begin{eqnarray*}
&& \hspace{-1cm}\mathbb{E}\, \left[  \left\vert R\widehat{y}_{T}\right\vert ^{2}\right]
+\mathbb{E}\, \left[  \int_{0}^{T}\left(  \left\vert R\widehat{y}_{s}\right\vert
^{2}+\left\Vert R\widehat{z}_{s}\right\Vert^{2}\right)  ds\right]  +\mathbb{E}\, %
\left[  \left\vert \widehat{Y}_{0}\right\vert ^{2}\right]  \\
&& \hspace{4cm}+\mathbb{E}\, \left[  \int_{0}^{T}\left(  \left\vert \widehat{Y}_{s}\right\vert
^{2}+\left\Vert \widehat{Z}_{s}\right\Vert^{2}+||| \widehat{k}_{s}|||^{2}\right)  ds\right]  \\
&& \hspace{0.5cm}\leq\frac{\delta C}{\Gamma_{1}}\left(  \mathbb{E}\, \left[  \left\vert
\widehat{\overline{Y}}_{0}\right\vert ^{2}\right]  +\mathbb{E}\, \left[  \left\vert
\widehat{Y}_{0}\right\vert ^{2}\right]  +\mathbb{E}\, \left[  \left\vert R\widehat{y}%
_{T}\right\vert ^{2}\right]  +\mathbb{E}\, \left[  \left\vert R\widehat{\overline{y}
}_{T}\right\vert ^{2}\right]  \right)  \\
&& \hspace{4cm}+\, \frac{\delta C}{\Gamma_{1}}\, \mathbb{E} \, \left[  \int_{0}^{T}\left(  \left\Vert
\widehat{\upsilon}_{s}\right\Vert^{2}+\left\Vert
\widehat{\overline{\upsilon}}_{s}\right\Vert^{2}\right)
ds\right].
\end{eqnarray*}
Taking $\delta=\frac{\Gamma_{1}}{3C},$ hence $\frac{\delta C}{\Gamma_{1}}%
=\frac{1}{3}$, we obtain (after multiplying the resulting inequality by $\frac{3}{2}$):
\begin{eqnarray*}
&& \mathbb{E}\, \left[  \left\vert R\widehat{y}_{T}\right\vert ^{2}\right]
+\mathbb{E}\, \left[  \left\vert \widehat{Y}_{0}\right\vert ^{2}\right]
+\mathbb{E}\, \left[  \int_{0}^{T}\left\Vert \widehat{\upsilon}_{s}\right\Vert
^{2}ds\right]  \\
&& \hspace{2cm}\leq\frac{1}{2}\left(  \mathbb{E}\, \left[  \left\vert R\widehat{\overline{y}}%
_{T}\right\vert ^{2}\right]  +\mathbb{E}\, \left[  \left\vert \widehat{\overline{Y}
}_{0}\right\vert ^{2}\right]  +\mathbb{E}\, \left[  \int_{0}^{T}\left\Vert
\widehat{\overline{\upsilon}}_{s}\right\Vert^{2}ds\right]
\right)  .
\end{eqnarray*}

II). $\theta_{2}=0,\, \beta_{2}>0.$ Let $\Gamma_{2}=\min\left\{  1-\alpha_{0}+\alpha_{0}%
\beta_{1},\theta_{1},\beta_{2}\right\}.$ Then $\Gamma_{2}>0,$ and inequality (\ref{eq:3.24}) becomes%
\begin{eqnarray}\label{eq:3.25}
&& \hspace{-0.5cm}\Gamma_{2}\left(  \mathbb{E}\, \left[  \left\vert R\widehat{y}_{T}\right\vert
^{2}\right]  +\mathbb{E}\, \left[  \int_{0}^{T}\left(  \left\vert R\widehat{y}%
_{s}\right\vert ^{2}+\left\Vert R\widehat{z}_{s}\right\Vert^{2}\right)
ds\right]  +\mathbb{E}\, \left[  \left\vert \widehat{Y}_{0}\right\vert ^{2}\right]
\right)\nonumber \\ &&
\hspace{2cm}\leq\delta C\left(  \mathbb{E}\, \left[  \left\vert R\widehat{y}_{T}\right\vert
^{2}\right]  +\mathbb{E}\, \left[  \left\vert R\widehat{\overline{y}}_{T}\right\vert
^{2}\right]  +\mathbb{E}\, \left[  \left\vert \widehat{Y}_{0}\right\vert ^{2}\right]
+\mathbb{E}\, \left[  \left\vert \widehat{\overline{Y}}_{0}\right\vert ^{2}\right]
\right)\nonumber\\
&&   \hspace{5cm}+\, \delta \, C\, \mathbb{E} \, \left[  \int_{0}^{T}\left(  \left\Vert \widehat{\upsilon
}_{s}\right\Vert^{2}+\left\Vert \widehat{\overline{\upsilon}}
_{s}\right\Vert^{2}\right)  ds\right] .
\end{eqnarray}

On the other hand, by using this inequality (\ref{eq:3.25}) and (\ref{eq:3.23}), we have%
\begin{eqnarray*}
&&   \hspace{-0.5cm} \Gamma_{2}\left(  \mathbb{E}\, \left[  \left\vert R\widehat{y}_{T}\right\vert
^{2}\right]  +\mathbb{E}\, \left[  \int_{0}^{T}\left(  \left\vert R\widehat{y}%
_{s}\right\vert ^{2}+\left\Vert R\widehat{z}_{s}\right\Vert^{2}\right)
ds\right]  \right.  \\
&& \left.   \hspace{1.5cm}+\, \mathbb{E}\, \left[  \left\vert \widehat{Y}_{0}\right\vert ^{2}\right]
+\mathbb{E}\, \left[  \int_{0}^{T}\left(  \left\vert \widehat{Y}_{s}\right\vert
^{2}+\left\Vert \widehat{Z}_{s}\right\Vert^{2}+||| \widehat{k}_{s}|||^{2}\right)  ds\right]  \right)  \\
&& \leq\delta \, C\,\left(  \mathbb{E}\, \left[  \left\vert R\widehat{y}_{T}\right\vert
^{2}\right]  +\mathbb{E}\, \left[  \left\vert R\widehat{\overline{y}}_{T}\right\vert
^{2}\right]  +\mathbb{E}\, \left[  \left\vert \widehat{Y}_{0}\right\vert ^{2}\right]
+\mathbb{E}\, \left[  \left\vert \widehat{\overline{Y}}_{0}\right\vert ^{2}\right]
\right)  \\
&&  \hspace{6cm}+\, \delta\, C \, \mathbb{E} \, \left[  \int_{0}^{T}\left(  \left\Vert \widehat{\upsilon
}_{s}\right\Vert^{2}+\left\Vert \widehat{\overline{\upsilon}}
_{s}\right\Vert^{2}\right)  ds\right]  \\
&& +\, \Gamma_{2} \, C \,\left(  \mathbb{E}\, \left[  \left\vert \widehat{y}%
_{T}\right\vert ^{2}\right]  +\delta\mathbb{E}\, \left[  \left\vert \widehat
{\overline{y}}_{T}\right\vert ^{2}\right]  \right)  +\Gamma_{2}C\,%
\mathbb{E}\, \left[  \int_{0}^{T}\left(  \left\vert \widehat{y}_{s}\right\vert
^{2}+\left\Vert \widehat{z}_{s}\right\Vert^{2}+\delta\left\Vert \widehat
{\overline{\upsilon}}_{s}\right\Vert^{2}\right)  ds\right]  \\
&& \leq \delta \, C \,\left(  \mathbb{E}\, \left[  \left\vert R\widehat{y}_{T}\right\vert
^{2}\right]  +\mathbb{E}\, \left[  \left\vert R\widehat{\overline{y}}_{T}\right\vert
^{2}\right]  +\mathbb{E}\, \left[  \left\vert \widehat{Y}_{0}\right\vert ^{2}\right]
+\mathbb{E}\, \left[  \left\vert \widehat{\overline{Y}}_{0}\right\vert ^{2}\right]
\right)  \\
&& +\, \delta\, C \, \mathbb{E} \, \left[  \int_{0}^{T}\left(  \left\Vert \widehat{\upsilon
}_{s}\right\Vert^{2}+\left\Vert \widehat{\overline{\upsilon}}
_{s}\right\Vert^{2}\right)  ds\right]  +\Gamma_{2}
C\delta\, \mathbb{E} \, \left[  \left\vert \widehat{\overline{y}}%
_{T}\right\vert ^{2}\right] \\
&&  +\, \Gamma_{2}\, C\, \delta\, \mathbb{E} \, \left[  \int_{0}^{T}\left\Vert \widehat{\overline{\upsilon}}_{s}\right\Vert
^{2}ds\right]  +\, \Gamma_{2}\, C\, \left(  \mathbb{E}\, \left[  \left\vert \widehat{y}%
_{T}\right\vert ^{2}\right]  +\mathbb{E}\, \left[  \int_{0}^{T}\left(  \left\vert
\widehat{y}_{s}\right\vert ^{2}+\left\Vert \widehat{z}_{s}\right\Vert^{2}\right)
ds\right]  \right)  .
\end{eqnarray*}

\medskip
Apply again (\ref{eq:3.25}) to the last term of this inequality and Remark~\ref{Remark 3.2} to get%
\begin{eqnarray*}
&& \hspace{-0.5cm}\Gamma_{2}\left(  \mathbb{E}\, \left[  \left\vert R\widehat{y}_{T}\right\vert
^{2}\right]  +\mathbb{E}\, \left[  \int_{0}^{T}\left(  \left\vert R\widehat{y}%
_{s}\right\vert ^{2}+\left\Vert R\widehat{z}_{s}\right\Vert^{2}\right)
ds\right]  \right.  \\
&& \left.  \hspace{1.5cm}+\mathbb{E}\, \left[  \left\vert \widehat{Y}_{0}\right\vert ^{2}\right]
+\mathbb{E}\, \left[  \int_{0}^{T}\left(  \left\vert \widehat{Y}_{t}\right\vert
^{2}+\left\Vert \widehat{Z}_{s}\right\Vert^{2}+||| \widehat{k}_{s}||| ^{2}\right)  ds\right]  \right)  \\
&& \leq\delta C\left(  \mathbb{E}\, \left[  \left\vert R\widehat{y}_{T}\right\vert
^{2}\right]  +\mathbb{E}\, \left[  \left\vert R\widehat{\overline{y}}_{T}\right\vert
^{2}\right]  +\mathbb{E}\, \left[  \left\vert \widehat{Y}_{0}\right\vert ^{2}\right]
+\mathbb{E}\, \left[  \left\vert \widehat{\overline{Y}}_{0}\right\vert ^{2}\right]
\right)  \\
&& \hspace{1.5cm}+\, \delta \, C\, \mathbb{E} \, \left[  \int_{0}^{T}\left(  \left\Vert \widehat{\upsilon
}_{s}\right\Vert^{2}+\left\Vert \widehat{\overline{\upsilon}}
_{s}\right\Vert^{2}\right)  ds\right]  +\Gamma_{2}\,
C\, \delta\, \mathbb{E} \, \left[  \left\vert \widehat{\overline{y}}%
_{T}\right\vert ^{2}\right]   \\
&& \hspace{1.5cm}+\, C\, \delta \left(  \mathbb{E}\, \left[  \left\vert R\widehat
{y}_{T}\right\vert ^{2}\right]  +\mathbb{E}\, \left[  \left\vert R\widehat
{\overline{y}}_{T}\right\vert ^{2}\right]  +\mathbb{E}\, \left[  \left\vert \widehat
{Y}_{0}\right\vert ^{2}\right]  +\mathbb{E}\, \left[  \left\vert \widehat{\overline{Y}}_{0}\right\vert ^{2}\right]  \right)  \\
&& \hspace{1.5cm}+\, \Gamma_{2}\, C \, \delta\, \mathbb{E} \, %
\left[  \int_{0}^{T}\left\Vert \widehat{\overline{\upsilon}}_{s}\right\Vert
^{2}ds\right] + C \, \delta \, \mathbb{E}\, \left[  \int_{0}^{T}\left(
\left\Vert \widehat{\upsilon}_{s}\right\Vert^{2}+\left\Vert
\widehat{\overline{\upsilon}}_{s}\right\Vert^{2}\right)
ds\right]  \\
&& \leq C\delta \,\left(  \mathbb{E}\, \left[  \left\vert R\widehat{y}%
_{T}\right\vert ^{2}\right]  +\mathbb{E}\, \left[  \left\vert R\widehat{\overline{y}
}_{T}\right\vert ^{2}\right]  +\mathbb{E}\, \left[  \left\vert \widehat{Y}%
_{0}\right\vert ^{2}\right]  +\mathbb{E}\, \left[  \left\vert \widehat{\overline{Y}
}_{0}\right\vert ^{2}\right]  \right)  \\
&& \hspace{6cm}+\, C\, \delta \,\mathbb{E}\, \left[  \int_{0}^{T}\left(  \left\Vert
\widehat{\upsilon}_{s}\right\Vert^{2}+\left\Vert
\widehat{\overline{\upsilon}}_{s}\right\Vert^{2}\right)
ds\right]  .
\end{eqnarray*}
As a result of this and Remark~\ref{Remark 3.2} we have
\begin{eqnarray*}
&& \hspace{-0.5cm}\mathbb{E}\, \left[  \left\vert R\widehat{y}_{T}\right\vert ^{2}\right]
+\mathbb{E}\, \left[  \left\vert \widehat{Y}_{0}\right\vert ^{2}\right]
+\mathbb{E}\, \left[  \int_{0}^{T}\left\Vert \widehat{\upsilon}_{s}\right\Vert
^{2}ds\right]  \\
&& \hspace{2cm}\leq\frac{C\delta}{\Gamma_{2}}\,\left(  \mathbb{E}\, \left[  \left\vert
R\widehat{y}_{T}\right\vert ^{2}\right]  +\mathbb{E}\, \left[  \left\vert
R\widehat{\overline{y}}_{T}\right\vert ^{2}\right]  +\mathbb{E}\, \left[  \left\vert
\widehat{Y}_{0}\right\vert ^{2}\right]  +\mathbb{E}\, \left[  \left\vert
\widehat{\overline{Y}}_{0}\right\vert ^{2}\right]  \right)  \\
&& \hspace{6cm}+\, \frac{C\, \delta}{\Gamma_{2}}~\mathbb{E}\, \left[  \int_{0}^{T}\left(
\left\Vert \widehat{\upsilon}_{s}\right\Vert^{2}+\left\Vert
\widehat{\overline{\upsilon}}_{s}\right\Vert^{2}\right)
ds\right]  .
\end{eqnarray*}

Now take $\delta=\frac{\Gamma_{2}}{3C}$ to get
\begin{eqnarray}\label{eq:3.26}
&& \hspace{-0.5cm} \mathbb{E}\, \left[  \left\vert R\widehat{y}_{T}\right\vert ^{2}\right]
+\mathbb{E}\, \left[  \left\vert \widehat{Y}_{0}\right\vert ^{2}\right]
+\mathbb{E}\, \left[  \int_{0}^{T}\left\Vert \widehat{\upsilon}_{s}\right\Vert
^{2}ds\right]  \nonumber \\
&& \hspace{2cm} \leq\frac{1}{2}\left(  \mathbb{E}\, \left[  \left\vert R\widehat{\overline{y}}%
_{T}\right\vert ^{2}\right]  +\mathbb{E}\, \left[  \left\vert \widehat{\overline{Y}
}_{0}\right\vert ^{2}\right]  +\mathbb{E}\, \left[  \int_{0}^{T}\left\Vert
\widehat{\overline{\upsilon}}_{s}\right\Vert^{2}ds\right]
\right)  .
\end{eqnarray}

III). $\theta_{2}>0,\, \beta_{2}=0.$ Let $\Gamma_{3}=\min\left\{  1-\alpha_{0}+\alpha_{0}\beta_{1},\theta_{1}%
,\theta_{2}\right\}.$ Then $\Gamma_{3}>0,$ and inequality (\ref{eq:3.24}) becomes%
\begin{eqnarray*}
&& \hspace{-0.5cm} \Gamma_{3}\left(  \mathbb{E}\, \left[  \left\vert R\widehat{y}_{T}\right\vert
^{2}\right]  +\mathbb{E}\, \left[  \int_{0}^{T}\left(  \left\vert R\widehat{y}%
_{s}\right\vert ^{2}+\left\Vert R\widehat{z}_{s}\right\Vert^{2}\right)
ds\right]  \right.  \\
&& \hspace{4.5cm} +\left.  \mathbb{E}\, \left[  \int_{0}^{T}\left(  \left\vert \widehat{Y}%
_{s}\right\vert ^{2}+\left\Vert \widehat{Z}_{s}\right\Vert^{2}+||| \widehat{k}_{s}||| ^{2}\right)  ds\right]  \right)
\\
&& \hspace{2cm} \leq\delta C\left(  \mathbb{E}\, \left[  \left\vert R\widehat{y}_{T}\right\vert
^{2}\right]  +\mathbb{E}\, \left[  \left\vert R\widehat{\overline{y}}_{T}\right\vert
^{2}\right]  +\mathbb{E}\, \left[  \left\vert \widehat{Y}_{0}\right\vert ^{2}\right]
+\mathbb{E}\, \left[  \left\vert \widehat{\overline{Y}}_{0}\right\vert ^{2}\right]
\right)  \\
&& \hspace{6cm} +\, \delta \, C\, \mathbb{E} \, \left[  \int_{0}^{T}\left(  \left\Vert \widehat{\upsilon
}_{s}\right\Vert^{2}+\left\Vert \widehat{\overline{\upsilon}}
_{s}\right\Vert^{2}\right)  ds\right]  ,
\end{eqnarray*}
which together with (\ref{eq:3.19}) implies that
\begin{eqnarray}\label{eq:3.27}
&& \hspace{-0.5cm} \Gamma_{3}\left(  \mathbb{E}\, \left[  \left\vert R\widehat{y}_{T}\right\vert
^{2}\right]  +\mathbb{E}\, \left[  \left\vert \widehat{Y}_{0}\right\vert ^{2}\right]
+\mathbb{E}\, \left[  \int_{0}^{T}\left\Vert \widehat{\upsilon}_{s}\right\Vert
^{2}ds\right]  \right)  \nonumber\\
&& \hspace{0.5cm} \leq\delta C\left(  \mathbb{E}\, \left[  \left\vert R\widehat{y}_{T}\right\vert
^{2}\right]  +\mathbb{E}\, \left[  \left\vert R\widehat{\overline{y}}_{T}\right\vert
^{2}\right]  +\mathbb{E}\, \left[  \left\vert \widehat{Y}_{0}\right\vert ^{2}\right]
+\mathbb{E}\, \left[  \left\vert \widehat{\overline{Y}}_{0}\right\vert ^{2}\right]
\right) \nonumber \\
&& \hspace{2cm} +\, \delta \, C\, \mathbb{E} \, \left[  \int_{0}^{T}\left(  \left\Vert \widehat{\upsilon
}_{s}\right\Vert^{2}+\left\Vert \widehat{\overline{\upsilon}}%
_{s}\right\Vert^{2}\right)  ds\right]  +\Gamma_{3}
\, \mathbb{E} \, \left[  \left\vert \widehat{Y}_{0}\right\vert ^{2}\right]  .
\end{eqnarray}

On the other hand, making use of (\ref{eq:3.19}) again and applying (\ref{eq:3.27})  together with Remark~\ref{Remark 3.2} give%
\begin{eqnarray}\label{eq:3.28}
&& \hspace{-0.5cm}\Gamma_{3}\, \mathbb{E}\, \left[  \left\vert \widehat{Y}_{0}\right\vert ^{2}\right]
\leq\Gamma_{3} \, C^{2}\, e^{T}\,\left(  C\left(  \mathbb{E}\, \left[  \left\vert \widehat{y}%
_{T}\right\vert ^{2}\right]  + \delta\, \mathbb{E} \, \left[  \left\vert \widehat
{\overline{y}}_{T}\right\vert ^{2}\right] \right) \right. \nonumber \\
&& \hspace{3.5cm}\left. +\, C\, \mathbb{E} \, \left[  \int
_{0}^{T}\left(  \left\vert \widehat{y}_{s}\right\vert ^{2}+\left\Vert \widehat{z}%
_{s}\right\Vert^{2}+\delta\left\Vert \widehat{\overline{\upsilon}}_{s}\right\Vert
^{2}\right)  ds\right]   \right) \nonumber\\
&& \hspace{2cm} \leq\delta \, C\,\left(  \mathbb{E}\, \left[  \left\vert R\widehat{y}_{T}\right\vert
^{2}\right]  +\mathbb{E}\, \left[  \left\vert R\widehat{\overline{y}}_{T}\right\vert
^{2}\right]  +\mathbb{E}\, \left[  \left\vert \widehat{Y}_{0}\right\vert ^{2}\right]
+\mathbb{E}\, \left[  \left\vert \widehat{\overline{Y}}_{0}\right\vert ^{2}\right]
\right)  \nonumber\\
&& \hspace{5cm} +\, \delta\, C\, \mathbb{E} \, \left[  \int_{0}^{T}\left(  \left\Vert \widehat{\upsilon
}_{s}\right\Vert^{2}+\left\Vert \widehat{\overline{\upsilon}}%
_{s}\right\Vert^{2}\right)  ds\right].
\end{eqnarray}

Therefore, substituting  (\ref{eq:3.28}) in (\ref{eq:3.27}) yields
\begin{eqnarray*}
&& \hspace{-0.5cm} \Gamma_{3} \left(  \mathbb{E}\, \left[  \left\vert R\widehat{y}_{T}\right\vert
^{2}\right]  +\mathbb{E}\, \left[  \left\vert \widehat{Y}_{0}\right\vert ^{2}\right]
+\mathbb{E}\, \left[  \int_{0}^{T}\left\Vert \widehat{\upsilon}_{s}\right\Vert
^{2}ds\right]  \right)  \\
&& \hspace{2cm} \leq\delta \, C\,\left(  \mathbb{E}\, \left[  \left\vert R\widehat{y}%
_{T}\right\vert ^{2}\right]  +\mathbb{E}\, \left[  \left\vert R\widehat{\overline{y}
}_{T}\right\vert ^{2}\right]  +\mathbb{E}\, \left[  \left\vert \widehat{Y}%
_{0}\right\vert ^{2}\right]  +\mathbb{E}\, \left[  \left\vert \widehat{\overline{Y}
}_{0}\right\vert ^{2}\right]  \right)  \\
&& \hspace{6cm} + \, \delta \, C\, \mathbb{E} \, \left[  \int_{0}^{T}\left(  \left\Vert
\widehat{\upsilon}_{s}\right\Vert^{2}+\left\Vert
\widehat{\overline{\upsilon}}_{s}\right\Vert^{2}\right)
ds\right]  .
\end{eqnarray*}
Hence, as in case III) by taking $\delta=\frac{\Gamma_{3}}{3C},$ we deduce that (\ref{eq:3.26}) holds also here.

\bigskip

IV). $\theta_{2}=0, \beta_{2}=0.$  Let $\Gamma_{4}=\min\{1-\alpha_{0}+\alpha_{0}\beta_{1},\theta_{1}\}.$ Then $\Gamma_{4}>0$ and inequality (\ref{eq:3.24}) becomes

\begin{eqnarray}\label{eq:3.29}
&& \hspace{-0.5cm} \Gamma_{4} \left(  \mathbb{E}\, \left[  \left\vert R\widehat{y}_{T}\right\vert
^{2}\right] +\mathbb{E}\, \left[  \int_{0}^{T}\left(\left\vert R\widehat{y}_{s}\right\vert
^{2}+\left\Vert R\widehat{z}_{s}\right\Vert
^{2}\right)ds\right]  \right)  \nonumber \\
&& \hspace{2cm} \leq\delta \, C\,\left(  \mathbb{E}\, \left[  \left\vert \widehat{\overline{Y}}%
_{0}\right\vert ^{2}\right]  +\mathbb{E}\, \left[  \left\vert \widehat{Y
}_{0}\right\vert ^{2}\right]  +\mathbb{E}\, \left[  \left\vert R\widehat{y}%
_{T}\right\vert ^{2}\right]  +\mathbb{E}\, \left[  \left\vert R\widehat{\overline{y}
}_{T}\right\vert ^{2}\right]  \right) \nonumber  \\
&& \hspace{4cm} +\, \delta \, C\, \mathbb{E} \, \left[  \int_{0}^{T}\left(  \left\Vert
\widehat{\upsilon}_{s}\right\Vert^{2}+\left\Vert
\widehat{\overline{\upsilon}}_{s}\right\Vert^{2}\right)
ds\right]  .
\end{eqnarray}
Thus by using (\ref{eq:3.23}) we deduce
\begin{eqnarray*}
&&   \hspace{-0.5cm} \Gamma_{4} \left(  \mathbb{E}\, \left[  \left\vert R\widehat{y}_{T}\right\vert
^{2}\right]  +\mathbb{E}\, \left[  \int_{0}^{T}\left(  \left\vert R\widehat{y}%
_{s}\right\vert ^{2}+\left\Vert R\widehat{z}_{s}\right\Vert^{2}\right)
ds\right]  \right.  \\
&& \left.  \hspace{2cm}  +\, \mathbb{E}\, \left[  \left\vert \widehat{Y}_{0}\right\vert ^{2}\right]
+\mathbb{E}\, \left[  \int_{0}^{T}\left(  \left\vert \widehat{Y}_{s}\right\vert
^{2}+\left\Vert \widehat{Z}_{s}\right\Vert^{2}+||| \widehat{k}_{s}||| ^{2}\right)  ds\right]  \right)  \\
&& \leq\delta C\,\left(  \mathbb{E}\, \left[  \left\vert R\widehat{y}_{T}\right\vert
^{2}\right]  +\mathbb{E}\, \left[  \left\vert R\widehat{\overline{y}}_{T}\right\vert
^{2}\right]  +\mathbb{E}\, \left[  \left\vert \widehat{Y}_{0}\right\vert ^{2}\right]
+\mathbb{E}\, \left[  \left\vert \widehat{\overline{Y}}_{0}\right\vert ^{2}\right]
\right)  \\
&&  \hspace{6cm}+\, \delta \, C \, \mathbb{E} \, \left[  \int_{0}^{T}\left(  \left\Vert \widehat{\upsilon
}_{s}\right\Vert^{2}+\left\Vert \widehat{\overline{\upsilon}}
_{s}\right\Vert^{2}\right)  ds\right]  \\
&& +\, \Gamma_{4}\, C\, \left(  \mathbb{E}\, \left[  \left\vert \widehat{y}%
_{T}\right\vert ^{2}\right]  +\delta\mathbb{E}\, \left[  \left\vert \widehat
{\overline{y}}_{T}\right\vert ^{2}\right]  \right)  +\Gamma_{4}C\,%
\mathbb{E}\, \left[  \int_{0}^{T}\left(  \left\vert \widehat{y}_{s}\right\vert
^{2}+\left\Vert \widehat{z}_{s}\right\Vert^{2}+\delta\,\left\Vert \widehat
{\overline{\upsilon}}_{s}\right\Vert^{2}\right)  ds\right]  \\
&& \leq \delta C \,\left(  \mathbb{E}\, \left[  \left\vert R\widehat{y}_{T}\right\vert
^{2}\right]  +\mathbb{E}\, \left[  \left\vert R\widehat{\overline{y}}_{T}\right\vert
^{2}\right]  +\mathbb{E}\, \left[  \left\vert \widehat{Y}_{0}\right\vert ^{2}\right]
+\mathbb{E}\, \left[  \left\vert \widehat{\overline{Y}}_{0}\right\vert ^{2}\right]
\right)  \\
&& +\, \delta \, C \, \mathbb{E} \, \left[  \int_{0}^{T}\left(  \left\Vert \widehat{\upsilon
}_{s}\right\Vert^{2}+\left\Vert \widehat{\overline{\upsilon}}
_{s}\right\Vert^{2}\right)  ds\right]  +\Gamma_{4}
C\delta\, \mathbb{E} \, \left[  \left\vert \widehat{\overline{y}}%
_{T}\right\vert ^{2}\right]  +\Gamma_{4}C\delta\, \mathbb{E} \, %
\left[  \int_{0}^{T}\left\Vert \widehat{\overline{\upsilon}}_{s}\right\Vert
^{2}ds\right]  \\
&&  \hspace{3.5cm}+\, \Gamma_{4}\, C \,\left(  \mathbb{E}\, \left[  \left\vert \widehat{y}%
_{T}\right\vert ^{2}\right]  +\mathbb{E}\, \left[  \int_{0}^{T}\left(  \left\vert
\widehat{y}_{s}\right\vert ^{2}+\left\Vert \widehat{z}_{s}\right\Vert^{2}\right)
ds\right]  \right)  .
\end{eqnarray*}
Apply again (\ref{eq:3.29}) and Remark~\ref{Remark 3.2} to get%
\begin{eqnarray*}
&& \hspace{-0.5cm}\Gamma_{4} \left(  \mathbb{E}\, \left[  \left\vert R\widehat{y}_{T}\right\vert
^{2}\right]  +\mathbb{E}\, \left[  \int_{0}^{T}\left(  \left\vert R\widehat{y}%
_{s}\right\vert ^{2}+\left\Vert R\widehat{z}_{s}\right\Vert^{2}\right)
ds\right]  \right.  \\
&& \left.  \hspace{1cm} +\, \mathbb{E}\, \left[  \left\vert \widehat{Y}_{0}\right\vert ^{2}\right]
+\mathbb{E}\, \left[  \int_{0}^{T}\left(  \left\vert \widehat{Y}_{t}\right\vert
^{2}+\left\Vert \widehat{Z}_{s}\right\Vert^{2}+||| \widehat{k}_{s}||| ^{2}\right)  ds\right]  \right)  \\
&& \leq\delta C\left(  \mathbb{E}\, \left[  \left\vert R\widehat{y}_{T}\right\vert
^{2}\right]  +\mathbb{E}\, \left[  \left\vert R\widehat{\overline{y}}_{T}\right\vert
^{2}\right]  +\mathbb{E}\, \left[  \left\vert \widehat{Y}_{0}\right\vert ^{2}\right]
+\mathbb{E}\, \left[  \left\vert \widehat{\overline{Y}}_{0}\right\vert ^{2}\right]
\right)  \\
&& +\, \delta \, C\, \mathbb{E} \, \left[  \int_{0}^{T}\left(  \left\Vert \widehat{\upsilon
}_{s}\right\Vert^{2}+\left\Vert \widehat{\overline{\upsilon}}
_{s}\right\Vert^{2}\right)  ds\right]  +\Gamma_{4}
C\delta\, \mathbb{E} \, \left[  \left\vert \widehat{\overline{y}}%
_{T}\right\vert ^{2}\right]  +\Gamma_{4}C \delta\, \mathbb{E} \, %
\left[  \int_{0}^{T}\left\Vert \widehat{\overline{\upsilon}}_{s}\right\Vert
^{2}ds\right]  \\
&& \hspace{1.5cm}+\, C\, \delta \left(  \mathbb{E}\, \left[  \left\vert R\widehat
{y}_{T}\right\vert ^{2}\right]  +\mathbb{E}\, \left[  \left\vert R\widehat
{\overline{y}}_{T}\right\vert ^{2}\right]  +\mathbb{E}\, \left[  \left\vert \widehat
{Y}_{0}\right\vert ^{2}\right]  + \mathbb{E}\, \left[  \left\vert \widehat{\overline{Y}}_{0}\right\vert ^{2}\right]  \right)  \\
&& \hspace{7cm}+ \,C\delta \, \mathbb{E}\, \left[  \int_{0}^{T}\left(
\left\Vert \widehat{\upsilon}_{s}\right\Vert^{2}+\left\Vert
\widehat{\overline{\upsilon}}_{s}\right\Vert^{2}\right)
ds\right]  \\
&& \leq C\, \delta \left(  \mathbb{E}\, \left[  \left\vert R\widehat{y}%
_{T}\right\vert ^{2}\right]  +\mathbb{E}\, \left[  \left\vert R\widehat{\overline{y}
}_{T}\right\vert ^{2}\right]  +\mathbb{E}\, \left[  \left\vert \widehat{Y}%
_{0}\right\vert ^{2}\right]  +\mathbb{E}\, \left[  \left\vert \widehat{\overline{Y}
}_{0}\right\vert ^{2}\right]  \right)  \\
&& \hspace{7cm}+ \, C\, \delta \, \mathbb{E}\, \left[  \int_{0}^{T}\left(  \left\Vert
\widehat{\upsilon}_{s}\right\Vert^{2}+\left\Vert
\widehat{\overline{\upsilon}}_{s}\right\Vert^{2}\right)
ds\right]  .
\end{eqnarray*}
It follows that
\begin{eqnarray*}
&& \hspace{-0.5cm}\mathbb{E}\, \left[  \left\vert R\widehat{y}_{T}\right\vert ^{2}\right]
+\mathbb{E}\, \left[  \left\vert \widehat{Y}_{0}\right\vert ^{2}\right]
+\mathbb{E}\, \left[  \int_{0}^{T}\left\Vert \widehat{\upsilon}_{s}\right\Vert
^{2}ds\right]  \\
&& \hspace{2cm}\leq\frac{C\delta}{\Gamma_{4}}\left(  \mathbb{E}\, \left[  \left\vert
R\widehat{y}_{T}\right\vert ^{2}\right]  +\mathbb{E}\, \left[  \left\vert
R\widehat{\overline{y}}_{T}\right\vert ^{2}\right]  +\mathbb{E}\, \left[  \left\vert
\widehat{Y}_{0}\right\vert ^{2}\right]  +\mathbb{E}\, \left[  \left\vert
\widehat{\overline{Y}}_{0}\right\vert ^{2}\right]  \right)  \\
&& \hspace{5cm}+ \, \frac{C\, \delta}{\Gamma_{4}} \, \mathbb{E}\, \left[  \int_{0}^{T}\left(
\left\Vert \widehat{\upsilon}_{s}\right\Vert^{2}+\left\Vert
\widehat{\overline{\upsilon}}_{s}\right\Vert^{2}\right)
ds\right]  .
\end{eqnarray*}
then we take $\delta=\frac{\Gamma_{2}}{3C}$ to get the same inequality as in (\ref{eq:3.26}).

From the preceding four cases we conclude that the mapping $I_{\alpha_{0}+\delta}$\ is contraction, in the sense
that%
\begin{eqnarray*}
&& \hspace{-0.5cm} \mathbb{E}\, \left[  \int_{0}^{T}\left\Vert \widehat{\upsilon}_{s}\right\Vert
^{2}ds\right]  +\mathbb{E}\, \left[  \left\vert \widehat{y}%
_{T}\right\vert ^{2}\right]  +\mathbb{E}\, \left[  \left\vert \widehat{Y}%
_{0}\right\vert ^{2}\right]  \\
&& \hspace{2cm} \leq\frac{1}{2}\left(  \mathbb{E}\, \left[  \int_{0}^{T}\left\Vert \widehat
{\overline{\upsilon}}_{s}\right\Vert^{2}ds\right]  +\mathbb{E}\, %
\left[  \left\vert \widehat{\overline{y}}_{T}\right\vert ^{2}\right]
+\mathbb{E}\, \left[  \left\vert \widehat{\overline{Y}}_{0}\right\vert ^{2}\right]
\right)  .
\end{eqnarray*}
Therefore this mapping has a unique fixed point $\upsilon=\left(
y,Y,z,Z,k\right)$  in $\mathbb{H}^{2},$ which can easily be seen to be the unique solution of (\ref{eq:3.3}) for
$\alpha=\alpha_{0}+\delta,\delta\in\left[  0,\delta_{0}\right]  .$
\end{proof}

\bigskip
\textbf{Case~2.} $m<n.$ If $m<n$, then $\theta_{2}>0$ and $\beta_{2}>0$. We consider the following system:
\begin{eqnarray}\label{eq:3.30}
\left\{
\begin{array}{ll}
dy_{t}=\left[  \alpha b\left(  t,\upsilon_{t}\right)  -\left(
1-\alpha\right)  \theta_{2}R^{\ast}Y_{t}+\widetilde{b}_{0}\left(  t\right)
\right]  dt-z_{t}d\overleftarrow{B}_{t}
 \\  \hspace{1cm}
+\left[  \alpha \sigma\left(  t,\upsilon_{t}\right)  -\left(
1-\alpha\right)  \theta_{2}R^{\ast}Z_{t}+\widetilde{\sigma}_{0}\left(  t\right)
\right]  dW_{t}\\  \hspace{1cm}
+\int_{\Theta}\left[  \alpha\varphi\left(  t,\upsilon_{t},\rho\right)
-\left(  1-\alpha\right)  \theta_{2}R^{\ast}k_{t}\left(  \rho\right)
+\varphi_{0}\left(  t,\rho\right)  \right]  \widetilde{N}\left(
d\rho,dt\right), \\ \\
dY_{t}=\left[  \alpha f\left(  t,\upsilon_{t}\right)  +\widetilde{f}%
_{0}\left(  t\right)  \right]  dt+\left[  \alpha g\left(
t,\upsilon_{t}\right)  +\widetilde{g}_{0}\left(  t\right)  \right]
d\overleftarrow{B}_{t}+Z_{t}dW_{t}  \\ \hspace{1cm}
+\int_{\Theta}k_{t}\left(  \rho\right)  \widetilde{N}\left(  d\rho,dt\right), \\ \\
y_{0}=\alpha\Psi\left(  Y_{0}\right)  -\left(  1-\alpha\right)  R^{\ast}%
Y_{0}+\psi,Y_{T}=\alpha h\left(  y_{T}\right)  +\phi.
\end{array}
\right.
\end{eqnarray}

When $\alpha=0$ equation (\ref{eq:3.30}) is uniquely solvable as shown in Case~1. Since $m<n,$ when $\alpha=0,$ system (\ref{eq:3.30}) becomes:
\begin{eqnarray}\label{eq:3.31}
\left\{
\begin{array}{ll}
dy_{t}=\left[  -\theta_{2}R^{\ast}Y_{t}+\widetilde{b}_{0}\left(  t\right)
\right]  dt-z_{t}d\overleftarrow{B}_{t}
+\left[  - \theta_{2}R^{\ast}Z_{t}+\widetilde{\sigma}_{0}\left(  t\right)
\right]  dW_{t}\\  \hspace{1cm}
+\int_{\Theta}\left[ - \theta_{2}R^{\ast}k_{t}\left(  \rho\right)
+\varphi_{0}\left(  t,\rho\right)  \right]  \widetilde{N}\left(
d\rho,dt\right), \\ \\
dY_{t}=\widetilde{f}%
_{0}\left(  t\right)  dt+\widetilde{g}_{0}\left(  t\right)
d\overleftarrow{B}_{t}+Z_{t}dW_{t}
+\int_{\Theta}k_{t}\left(  \rho\right)  \widetilde{N}\left(  d\rho,dt\right), \\ \\
y_{0}=  -R^{\ast}%
Y_{0}+\psi,Y_{T}=\phi.
\end{array}
\right.
\end{eqnarray}
So, by following the same procedure done for (\ref{eq:3.4}), we derive a unique solution $\upsilon=\left(
y,Y,z,Z,k\right)$  in $\mathbb{H}^{2}$\ for (\ref{eq:3.31}).

When $\alpha=1$\ the
existence of the solution of (\ref{eq:3.31}) implies clearly that of (\ref{eq:3.2}). By the same
techniques used for proving Lemma~\ref{Lemma: 3.5}, one can prove the following lemma.
\begin{lemma}\label{lem:final-lemma}
Assume $m<n$. Under assumptions (A1)-(A4),  with $0< \gamma<1$ and $0< \gamma' \leq \gamma/2 ,$ there exists
a positive constant $\delta_{0}$\ such that if, apriori, for each $\psi\in
L^{2}\left(  \Omega,\mathcal{F}_{0},\mathbb{P};\mathbb{R}^{n}\right)  ,\phi\in
L^{2}\left(  \Omega,\mathcal{F}_{T},\mathbb{P};\mathbb{R}^{m}\right)  $\ and $\left(
\widetilde{b}_{0},\widetilde{f}_{0},\widetilde{\sigma}_{0},\widetilde{g}_{0},\varphi
_{0}\right)  \in\mathbb{H}^{2}$, (\ref{eq:3.30}) is uniquely solvable for some
$\alpha_{0}\in\left[  0,1\right)  $, then for each $\alpha\in\left[
\alpha_{0},\alpha_{0}+\delta_{0}\right]  $\ and $\psi\in L^{2}\left(
\Omega,\mathcal{F}_{0},\mathbb{P};\mathbb{R}^{n}\right)  ,\phi\in L^{2}\left(
\Omega,\mathcal{F}_{T},\mathbb{P};\mathbb{R}^{m}\right)  ,$ $\left(  \widetilde{b}%
_{0},\widetilde{f}_{0},\widetilde{\sigma}_{0},\widetilde{g}_{0},\varphi_{0}\right)
\in\mathbb{H}^{2}$, (\ref{eq:3.30}) is also uniquely solvable in $\mathbb{H}^{2}$.
\end{lemma}
\begin{proof}  Assume that, for each $\psi\in L^{2}\left(  \Omega,\mathcal{F}%
_{0},\mathbb{P};\mathbb{R}^{n}\right)  ,\phi\in L^{2}\left(  \Omega,\mathcal{F}%
_{T},\mathbb{P};\mathbb{R}^{m}\right)  $\ and $(  \widetilde{b}_{0},\widetilde{f}
_{0},%
\widetilde{\sigma}_{0},\widetilde{g}_{0},\varphi_{0})  \in\mathbb{H}^{2},$\
there exists a unique solution of (\ref{eq:3.3}) for $\alpha=\alpha_{0}$. Then for each
$\overline{\upsilon}_{t}:=\left(  \overline{y}_{t},\overline{Y}_{t},\overline{z}_{t},\overline{Z}
_{t},\overline{k}_{t}\right)  \in\mathbb{H}^{2}$, there exists a unique element
$\upsilon:=\left(  y,Y,z,Z,k\right)$  of $\mathbb{H}%
^{2}$\ satisfying the following FBDSDEJ:
\begin{eqnarray}\label{eq:3.32}
\left\{
\begin{array}{ll}%
dy_{t}=[  \alpha_{0}b\left(  t,\upsilon_{t}\right) -(1-\alpha_{0})\theta_{2}R^{\ast}Y_{t}\\ \hspace{4cm}
+\delta\left( b(  t,\overline{\upsilon}_{t})+\theta_{2}R^{\ast} \overline{Y}_{t}\right) +\widetilde{b}_{0}\left(
t\right)  ]  dt-z_{t}d\overleftarrow{B}_{t}\\ \hspace{1.5cm}
+[  \alpha_{0}\sigma\left(  t,\upsilon_{t}\right) -(1-\alpha_{0})\theta_{2}R^{\ast}Z_{t}\\ \hspace{4cm}
+\delta\left(  \sigma(  t,\overline{\upsilon}_{t})+\theta_{2}R^{\ast} \overline{Z}_{t} \right)   +\widetilde{\sigma}_{0}\left(
t\right)  ]  dW_{t}\\ \hspace{1.5cm}
+\int_{\Theta}[  \alpha_{0}\varphi\left(  t,\upsilon_{t}%
,\rho\right)  -(1-\alpha_{0})\theta_{2}R^{\ast}k_{t}(\rho)\\ \hspace{4cm}
+\delta\left(  \varphi(  t,\overline{\upsilon}_{t},\rho)+\theta_{2}R^{\ast} \overline{k}_{t}(\rho)\right)
+\varphi_{0}\left(  t,\rho\right)  ]  \widetilde{N}\left(  d\rho,dt\right)\\ \\
dY_{t}=\left[  \alpha_{0}f\left(  t,\upsilon_{t}\right) +\delta  f\left(
t,\overline{\upsilon}_{t}\right)   +\widetilde{f}%
_{0}\left(  t\right)  \right]  dt\\ \hspace{1.5cm}
+\left[  \alpha_{0}g\left(  t,\upsilon_{t}\right) +\delta g\left(
t,\overline{\upsilon}_{t}\right)   +\widetilde{g}%
_{0}\left(  t\right)  \right]  d\overleftarrow{B}_{t}
+Z_{t}dW_{t}+\int_{\Theta}k_{t}\left(  \rho\right)  \widetilde{N}\left(
d\rho,dt\right) ,
\\ \\
y_{0}=\alpha_{0}\Psi\left(  Y_{0}\right)  -\left(  1-\alpha_{0}\right)
R^{\ast}Y_{0}+\delta\left(  \Psi\left(  \overline{Y}_{0}\right)  -R^{\ast}\overline{Y}
_{0}\right)  +\psi,\\
Y_{T}=\alpha_{0}h \left(  y_{T}\right)  +\delta h\left(  \overline{y}
_{T}\right)  +\phi ,
\end{array}
\right.
\end{eqnarray}
where $\delta\in (0,1)$ is a parameter independent of $\alpha_{0}$ which is small enough, and will be determined later in the proof.

Let us consider the mapping $I_{\alpha_{0}+\delta}$,  defined in the proof of Lemma~\ref{Lemma: 3.5}, and let
\begin{align*}
& \overline{\upsilon}_{\cdot} :=\left(  \overline{y}_{\cdot},\overline{Y}
_{\cdot},\overline{z}_{\cdot},\overline{Z}_{\cdot},\overline{k}_{\cdot}\right), \overline{\upsilon}_{\cdot}^{\prime} :=\left(  \overline{y}_{\cdot}^{\prime},\overline{Y}
_{\cdot}^{\prime},\overline{z}_{\cdot}^{\prime},\overline{Z}_{\cdot}^{\prime},\overline{k}_{\cdot}^{\prime
}\right)  \in\mathbb{H}^{2}, \\
& (\left(  y_{\cdot},Y_{\cdot},z_{\cdot},Z_{\cdot},k_{\cdot}\right), y_{T},Y_{0})= \left({\upsilon}_{\cdot},{y}_{T},{Y}_{0}\right) := I_{\alpha_{0}+\delta
}\left(  \overline{\upsilon}_{\cdot}, \overline{y}_{T},\overline{Y}_{0}\right)  ,\\
& (\left(  y_{\cdot}^{\prime},Y_{\cdot}^{\prime},z_{\cdot}^{\prime},Z_{\cdot}^{\prime},k_{\cdot}^{\prime}\right), y^{\prime}_{T},Y^{\prime}_{0})  = (\upsilon_{\cdot}^{\prime}, y^{\prime}_{T},Y^{\prime}_{0}) :=I_{\alpha_{0}+\delta
}\left(  \overline{\upsilon}^{\prime}_{\cdot}, \overline{y}^{\prime}_{T},\overline{Y}^{\prime}_{0} \right).
\end{align*}
We keep also the same notations such as $\widehat{\upsilon}_{\cdot},...$etc., which is set after system~(\ref{eq:3.6}) in the proof of
Lemma~\ref{Lemma: 3.5}.

Applying It\^{o}'s formula to $\left\langle R\widehat{y},\widehat{Y}\right\rangle
$\ on $\left[  0,T\right]  $ yields
\begin{eqnarray*}
&& \hspace{-0.75cm}  d\left\langle \widehat{y}_{t},R^{\ast}\widehat{Y}_{t}\right\rangle  =\left\langle \widehat{y}_{t},\alpha_{0}R^{\ast}\widehat{f}\left(
t,\upsilon_{t}\right) +\delta R^{\ast} \widehat{f}\left(  t,\overline{\upsilon}_{t}\right)
\right\rangle dt\\
&& \hspace{-0.5cm}  +\left\langle R^{\ast}\widehat{Y}_{t} ,\alpha_{0}\widehat{b}\left(
t,\upsilon_{t}\right)  -\left(  1-\alpha_{0}\right)  \theta_{2}R^{\ast}\widehat{Y}_{t}+\delta\left(  \widehat{b}\left(  t,\overline{\upsilon}_{t}\right)
+\theta_{2}R^{\ast}\widehat{\overline{Y}}_{t}\right)  \right\rangle dt
\\
&& \hspace{-0.5cm} + \, \left\langle \alpha_{0}\widehat{\sigma
}\left(  t,\upsilon_{t}\right)   -\left(  1-\alpha_{0}\right)  \theta_{2}R^{\ast}\widehat{Z}_{t}+\delta\left(  \widehat{\sigma}\left(  t,\overline{\upsilon}_{t}\right)
+\theta_{2}R^{\ast}\widehat{\overline{Z}}_{t}\right) ,R^{\ast}\widehat{Z}_{t}\right\rangle dt\\
&& \hspace{-0.5cm}  + \, \left\langle \widehat{z}_{t},\alpha_{0}R^{\ast}\widehat{g}\left(
t,\upsilon_{t}\right)  +\delta R^{\ast} \widehat{g}\left(  t,\overline{\upsilon}_{t}\right)
\right\rangle dt-\left\langle R^{\ast}\widehat{Y}_{t},\widehat{z}_{t}d\overleftarrow{B}_{t}\right\rangle \\
&& \hspace{-0.5cm}  + \, \int_{\Theta}\left\langle \alpha_{0}\widehat{\varphi}\left(
t,\upsilon_{t},\rho\right) -\left(  1-\alpha_{0}\right)  \theta_{2} R^{\ast} \widehat{k}_{t}(\rho) \right.
\\ && \hspace{3cm} \left.
 + \, \delta\left(  \widehat{\varphi}\left(  t,\overline{\upsilon}_{t},\rho\right)
+\theta_{2}R^{\ast} \widehat{\overline{k}}_{t}(\rho)\right)  , R^{\ast}\widehat{k}_{t}\left(  \rho\right)
\right\rangle \Pi\left(  d\rho\right)  dt\\
&& \hspace{-0.5cm}  + \, \left\langle R^{\ast}\widehat{Y}_{t},\int_{\Theta}\left( \alpha_{0}\widehat{\varphi}\left(
t,\upsilon_{t},\rho\right) -\left(  1-\alpha_{0}\right)  \theta_{2} R^{\ast} \widehat{k}_{t}(\rho) \right. \right.
\\ && \hspace{2cm}
\left. \left.
 +\, \delta\left(  \widehat{\varphi}\left(  t,\overline{\upsilon}_{t},\rho\right)
+\theta_{2}R^{\ast} \widehat{\overline{k}}_{t}(\rho)\right)\right)\widetilde{N}\left(  d\rho,dt\right)
\right\rangle \\
&& \hspace{-0.5cm}  + \, \left\langle R^{\ast}\widehat{Y}_{t},\left(\alpha_{0}\widehat{\sigma}\left(  t,\upsilon
_{t}\right)  -\left(  1-\alpha_{0}\right)  \theta_{2}R^{\ast}\widehat{Z}_{t}+\delta\left(  \widehat{\sigma}\left(  t,\overline{\upsilon}_{t}\right)
+\theta_{2}R^{\ast}\widehat{\overline{Z}}_{t}\right) \right)dW_{t}\right\rangle \\
&& \hspace{-0.5cm}  + \, \left\langle \widehat{y}_{t},\left(\alpha_{0}R^{\ast}\widehat{g}\left(
t,\upsilon_{t}\right)  + \delta R^{\ast}\widehat{g}\left(  t,\overline{\upsilon}_{t}\right)\right)
d\overleftarrow{B}_{t}\right\rangle \\
&& \hspace{3cm}  + \, \left\langle \widehat{y}_{t},R^{\ast}\widehat{Z}_{t}dW_{t}\right\rangle +\int_{\Theta
}\left\langle \widehat{y}_{t},R^{\ast}\widehat{k}_{t}\left(  \rho\right) \widetilde{N}\left(  d\rho,dt\right)  \right\rangle
 ,
\end{eqnarray*}
with
\[
\widehat{y}_{0}  = \alpha_{0}\widehat{\Psi}\left(  Y_{0}\right)  -\left(  1-\alpha_{0}\right)
R^{\ast}\widehat{Y}_{0}+\delta\left(  \widehat{\Psi}\left(  \overline{Y}_{0}\right)  -R^{\ast}\widehat{\overline{Y}}
_{0}\right),\]
and
\[
R^{\ast}\widehat{Y}_{T} =\alpha_{0}R^{\ast}\widehat{h} \left(  y_{T}\right)  +\delta R^{\ast} \widehat{h}\left(  \overline{y}
_{T}\right) .
\]

\bigskip

Integrating from $0$ to $T,$ taking expectation and using assumptions (A3) give
\begin{eqnarray*}
&& \hspace{-1.25cm} \mathbb{E}\, \left[  \left\langle \widehat{y}_{T},R^{\ast}\widehat{Y}_{T}\right\rangle \right] \\
&& \hspace{-0.5cm}   =\mathbb{E}\, \left[  \left\langle \widehat{y}_{0},R^{\ast}\widehat{Y}_{0}\right\rangle \right] +\mathbb{E}\, \left[  \int_{0}^{T}\left\langle \widehat{y}%
_{s}, \alpha_{0}R^{\ast}\widehat{f}\left(
s,\upsilon_{s}\right) +\delta R^{\ast}\widehat{f}\left( s, \overline{\upsilon}_{s}\right) \right\rangle ds\right] \\
&& \hspace{-0.5cm}
+\, \mathbb{E}\, \left[  \int_{0}^{T}\left\langle R^{\ast}\widehat{Y}_{s},\alpha_{0}\widehat{b}\left(
s,\upsilon_{s}\right)-\left(  1-\alpha_{0}\right)  \theta_{2}R^{\ast}\widehat{Y}_{s} +\delta \left( \widehat{b%
}\left(  s,\overline{\upsilon}_{s}\right)+\theta_{2}R^{\ast}\widehat{\overline{Y}}_{s} \right) \right\rangle ds\right] \\
&& \hspace{-0.5cm}    +\, \mathbb{E}\, \left[
\int_{0}^{T}\left\langle \widehat{z}_{s},\alpha_{0}R^{\ast}\widehat{g}\left(
s,\upsilon_{s}\right) +\delta R^{\ast}\widehat{g}\left(
s,\overline{\upsilon}_{s}\right) \right\rangle ds\right]\\
&& \hspace{-0.5cm}  +\, \mathbb{E}\, \left[  \int_{0}^{T}\left\langle \alpha_{0}\widehat{\sigma}\left(
s,\upsilon_{s}\right)-\left(  1-\alpha_{0}\right)  \theta_{2}R^{\ast}\widehat{Z}_{s} +\delta \left( \widehat{\sigma%
}\left(  s,\overline{\upsilon}_{s}\right)+\theta_{2}R^{\ast}\widehat{\overline{Z}}_{s} \right),R^{\ast}\widehat{Z}_{s} \right\rangle dt\right]\\
&& \hspace{-0.5cm}  +\, \mathbb{E}\, \big[  \int_{0}^{T}\int_{\Theta
}\left\langle \alpha_{0}\widehat{\varphi}\left(
s,\upsilon_{s},\rho\right)-\left(  1-\alpha_{0}\right)  \theta_{2}R^{\ast}\widehat{k}_{s}(\rho) +\delta \left( \widehat{\varphi%
}\left(  s,\overline{\upsilon}_{s},\rho\right), \right. \right. \\
&&\hspace{2.5in} \left. \left. +\, \theta_{2}R^{\ast}\widehat{\overline{k}}_{s}(\rho) \right) R^{\ast}\widehat{k}_{s}(\rho) \right\rangle \Pi(d\rho)ds\big]  ,
\end{eqnarray*}
Using notation preceding (A1) this equality reads as
\begin{eqnarray}\label{eq:3.33}
&& \hspace{-0.5cm}\mathbb{E}\, \left[  \left\langle \widehat{y}_{T},\alpha_{0}R^{\ast}\widehat{h}\left(
y_{T}\right)+\delta R^{\ast} \widehat{h}\left(
\overline{y}_{T}\right)\right\rangle \right]\nonumber \\
&& =\mathbb{E}\, \left[  \left\langle \alpha_{0}\widehat{\Psi}\left(
Y_{0}\right) -\left(  1-\alpha_{0}\right)R^{\ast}\widehat{Y}_{0} +\delta \left( \widehat{\Psi}\left(  \overline{Y}_{0}\right)- R^{\ast}\widehat{\overline{Y}}_{0}\right) ,R^{\ast}\widehat{Y}%
_{0}\right\rangle \right] \nonumber\\
&& +\, \alpha_{0}\mathbb{E}\, \left[  \int_{0}^{T}\left\langle \widehat{A}\left(
s,\upsilon_{s}\right)
,\widehat{\upsilon}_{s}\right\rangle ds\right]
 +\delta\mathbb{E}\, \left[  \int_{0}^{T}\left\langle \widehat{y}_{s}%
,R^{\ast}\widehat{f}\left(  s,\overline{\upsilon}_{s}\right)  \right\rangle
ds\right]\nonumber \\
&& +\, \mathbb{E}\, \left[  \int_{0}^{T}\left\langle R^{\ast}\widehat{Y}_{s},-\left(  1-\alpha_{0}\right)  \theta_{2}R^{\ast}\widehat{Y}_{s} +\delta \left( \widehat{b%
}\left(  s,\overline{\upsilon}_{s}\right)+\theta_{2}R^{\ast}\widehat{\overline{Y}}_{s} \right) \right\rangle ds\right]\nonumber \\
&& +\, \delta \, \mathbb{E} \, \left[
\int_{0}^{T}\left\langle \widehat{z}_{s},R^{\ast}\widehat{g}\left(
s,\overline{\upsilon}_{s}\right) \right\rangle ds\right]\nonumber \\
&&  +\, \mathbb{E}\, \left[  \int_{0}^{T}\left\langle -\left(  1-\alpha_{0}\right)  \theta_{2}R^{\ast}\widehat{Z}_{s} +\delta \left( \widehat{\sigma%
}\left(  s,\overline{\upsilon}_{s}\right)+\theta_{2}R^{\ast}\widehat{\overline{Z}}_{s} \right),R^{\ast}\widehat{Z}_{s} \right\rangle ds\right] \nonumber \\
&&  +\, \mathbb{E}\, \left[  \int_{0}^{T}\int_{\Theta
}\left\langle -\left(  1-\alpha_{0}\right)  \theta_{2}R^{\ast}\widehat{k}_{s}(\rho) \right. \right. \nonumber \\ && \hspace{2cm} \left. \left. + \, \delta \left( \widehat{\varphi%
}\left(  s,\overline{\upsilon}_{s},\rho\right)+\theta_{2}R^{\ast}\widehat{\overline{k}}_{s}(\rho) \right),R^{\ast}\widehat{k}_{s}(\rho) \right\rangle \Pi(d\rho)ds\right].
\end{eqnarray}

We know from (A2) that $\left\langle \alpha_{0}\widehat{\Psi}(Y_{0}),R^{\ast}\widehat{Y}_{0}\right\rangle \leq -\alpha_{0} \beta_{2}|R^{\ast}\widehat{Y}_{0}|^{2},$ which implies with the help of Cauchy-Schwartz inequality that
\begin{eqnarray}\label{eq:3.36}
&& \hspace{-1.5cm} \mathbb{E}\, \left[  \left\langle \alpha_{0}\widehat{\Psi}\left(
Y_{0}\right) -\left(  1-\alpha_{0}\right)R^{\ast}\widehat{Y}_{0} +\delta \left( \widehat{\Psi}\left(  \overline{Y}_{0}\right)- R^{\ast}\widehat{\overline{Y}}_{0}\right) ,R^{\ast}\widehat{Y}%
_{0}\right\rangle \right]
\nonumber\\
&& \leq -\, \alpha_{0} \, \beta_{2}\, \mathbb{E} \, \left[|R^{\ast}\widehat{Y}_{0}|^{2}\right]-\left(  1-\alpha_{0}\right)\, \mathbb{E} \, \left[|R^{\ast}\widehat{Y}_{0}|^{2}\right]+\delta c\, \mathbb{E} \, \left[  |\widehat{\overline{Y}}_{0}|^{2}\right]\nonumber\\
&& \hspace{0.9in}+\, \delta \, c\, \mathbb{E} \, \left[  |R^{\ast}\widehat{Y}%
_{0}|^{2}\right]+\delta \, \mathbb{E} \, \left[  | R^{\ast}\widehat{\overline{Y}}_{0}|^{2}\right]+\delta \, \mathbb{E} \, \left[  |R^{\ast}\widehat{Y}%
_{0}|^{2}\right].
\end{eqnarray}

On the other hand, apply (A1) on $A$ and (A2) on $\widehat{h}$, and use (A3) for the term $ -\delta\, \mathbb{E}\, \left[\left\langle \widehat{y}_{T},R^{\ast}\widehat{h}(\overline{y}_{T})\right\rangle\right]$ appearing in (\ref{eq:3.33}) to get
\begin{eqnarray*}
&& \hspace{-0.5cm}\alpha_{0} \beta_{1}\, \mathbb{E} \, \left[  |R\widehat{y}_{T}|^{2} \right]+\alpha_{0}\theta_{1}\, \mathbb{E} \, \left[\int_{0}^{T}\left(  |R\widehat{y}_{s}|^{2}+\|R\widehat{z}_{s}\|^{2}\right)ds \right] \\
&& \hspace{0.7in}
+ \, \alpha_{0}\theta_{2}\, \mathbb{E} \, \left[\int_{0}^{T}\left(  |R^{\ast}\widehat{Y}_{s}|^{2}+\|R^{\ast}\widehat{Z}_{s}\|^{2}+
|||R^{\ast} \widehat{k}_{s}|||\right)ds \right] \\
&& \leq \delta\left(  \mathbb{E}\, \left[|R\widehat{y}_{T}|^{2}\right]+c\, \mathbb{E} \, \left[|\widehat{\overline{y}}_{T}|^{2}\right]\right)-\alpha_{0} \beta_{2}\, \mathbb{E} \, \left[|R^{\ast}\widehat{Y}_{0}|^{2}\right]-\left(  1-\alpha_{0}\right)\, \mathbb{E} \, \left[|R^{\ast}\widehat{Y}_{0}|^{2}\right]\\
&& + \, \delta c\, \mathbb{E} \, \left[  |\widehat{\overline{Y}}_{0}|^{2}\right]+\delta c\, \mathbb{E} \, \left[  |R^{\ast}\widehat{Y}_{0}|^{2}\right]+\delta \, \mathbb{E} \, \left[  | R^{\ast}\widehat{\overline{Y}}_{0}|^{2}\right]+\delta \, \mathbb{E} \, \left[  |R^{\ast}\widehat{Y}_{0}|^{2}\right] \\
&& + \, \delta \, \mathbb{E} \, \left[ \int_{0}^{T} |R\widehat{y}_{s}|^{2}ds\right] +\delta\, \mathbb{E} \, \left[ \int_{0}^{T} |\widehat{f}(s,\overline{\upsilon}_{s})|^{2}ds\right]\\
&& -\left(  1-\alpha_{0}\right)  \theta_{2}\, \mathbb{E} \, \left[ \int_{0}^{T} |R^{\ast}\widehat{Y}_{s}|^{2}ds\right] +\delta \, \mathbb{E} \, \left[ \int_{0}^{T} |R^{\ast}\widehat{Y}_{s}|^{2}ds\right] + \, \delta\, \mathbb{E} \, \left[ \int_{0}^{T} |\widehat{b}(s,\overline{\upsilon}_{s})|^{2}ds\right] \\
&& + \, \delta \theta_{2}\, \mathbb{E} \, \left[ \int_{0}^{T} |R^{\ast}\widehat{Y}_{s}|^{2}ds\right] +\delta \theta_{2}\, \mathbb{E} \, \left[ \int_{0}^{T} |R^{\ast}\widehat{\overline{Y}}_{s}|^{2}ds\right]\\
&& + \, \delta \, \mathbb{E} \, \left[ \int_{0}^{T} \|R\widehat{z}_{s}\|^{2}ds\right] +\delta\, \mathbb{E} \, \left[ \int_{0}^{T} |\widehat{g}(s,\overline{\upsilon}_{s})|^{2}ds\right]\\
&&  -\left(  1-\alpha_{0}\right)  \theta_{2}\, \mathbb{E} \, \left[ \int_{0}^{T} \|R^{\ast}\widehat{Z}_{s}\|^{2}ds\right] +\delta \, \mathbb{E} \, \left[ \int_{0}^{T} \|R^{\ast}\widehat{Z}_{s}\|^{2}ds\right] \\ && + \,  \,\delta\, \mathbb{E} \, \left[ \int_{0}^{T} |\widehat{\sigma}(s,\overline{\upsilon}_{s})|^{2}ds\right] +\delta \theta_{2}\, \mathbb{E} \, \left[ \int_{0}^{T} \|R^{\ast}\widehat{\overline{Z}}_{s}\|^{2}ds\right] +\delta \theta_{2}\, \mathbb{E} \, \left[ \int_{0}^{T} \|R^{\ast}\widehat{Z}_{s}\|^{2}ds\right]\\
&&  -\left(  1-\alpha_{0}\right)  \theta_{2}\, \mathbb{E} \, \left[ \int_{0}^{T} |||R^{\ast} \widehat{k}_{s}|||^{2}ds\right] +\delta \, \mathbb{E} \, \left[ \int_{0}^{T} |||R^{\ast} \widehat{k}_{s}|||^{2}ds\right]\\
&& + \, \delta\, \theta_{2}\, \mathbb{E} \, \left[ \int_{0}^{T} ||| R^{\ast}\widehat{\overline{k}}_{s}|||^{2}ds\right]+\delta \, \theta_{2}\, \mathbb{E} \, \left[ \int_{0}^{T} |||R^{\ast} \widehat{k}_{s}|||^{2}ds\right] \\ && \hspace{3in} + \, \delta\, \mathbb{E} \, \left[ \int_{0}^{T} \||\widehat{\varphi}(s,\overline{\upsilon}_{s}, \cdot)\||^{2}ds\right].
\end{eqnarray*}
Rearranging this inequality and applying (A3) yield
\begin{eqnarray*}
&& \hspace{-0.5cm} \alpha_{0}\beta_{2}\, \mathbb{E} \, \left[  \left\vert R^{\ast}\widehat{Y}_{0}\right\vert
^{2}\right] +(1-\alpha_{0})\theta_{2}\, \mathbb{E} \, \left[  \int_{0}^{T}\left(
\left\vert R^{\ast}\widehat{Y}_{s}\right\vert ^{2}+\left\Vert R^{\ast}\widehat{Z}_{s}\right\Vert^{2}+|||R^{\ast} \widehat{k}_{s}||| ^{2}\right)  ds\right]    \\
&& \hspace{0.5in}+ \,\alpha_{0}\, \theta_{2}\, \mathbb{E} \, \left[  \int_{0}^{T}\left(
\left\vert R^{\ast}\widehat{Y}_{s}\right\vert ^{2}+\left\Vert R^{\ast}\widehat{Z}_{s}\right\Vert^{2}+|||R^{\ast} \widehat{k}_{s}||| ^{2}\right)  ds\right]\\
&& \hspace{0.5in}+\, \alpha_{0} \, \theta
_{1}\, \mathbb{E} \, \left[  \int_{0}^{T}\left(  \left\vert R\widehat{y}_{s}\right\vert ^{2}+\left\Vert R\widehat{z}%
_{s}\right\Vert^{2}\right)  ds\right] +(1-\alpha_{0})\, \mathbb{E} \, \left[  \left\vert R^{\ast}\widehat{Y}_{0}\right\vert
^{2}\right]\\
&& \leq -\alpha_{0}\beta_{1}\, \mathbb{E} \, \left[  |R\widehat{y}_{T}|^{2} \right]+\delta\left(\mathbb{E}\, \left[  |R\widehat{y}_{T}|^{2}\right] +c\, \mathbb{E} \, \left[  |\widehat{\overline{y}}_{T}|^{2}\right]\right)+\delta c\, \mathbb{E} \, \left[  |\widehat{\overline{Y}}_{0}|^{2}\right]\\
&&\hspace{0.1in}
+\, \delta \, c\, \mathbb{E} \, \left[  |R^{\ast}\widehat{Y}_{0}|^{2}\right]+\delta \, \mathbb{E} \, \left[  |R^{\ast}\widehat{\overline{Y}}_{0}|^{2}\right]+\delta \, \mathbb{E} \, \left[  |R^{\ast}\widehat{Y}_{0}|^{2}\right]+\delta\, \mathbb{E} \, \left[  \int_{0}^{T}\left\vert R\widehat{y}_{s}\right\vert ^{2} ds\right]  \\
&& \hspace{0.1in}+ \, 3\, \delta \, C\, \mathbb{E} \, \left[  \int_{0}^{T}\| \widehat{\overline{\upsilon}}_{s}\|^{2}ds \right] +\delta\, \mathbb{E} \, \left[  \int_{0}^{T}\left\vert R^{\ast}\widehat{Y}_{s}\right\vert ^{2} ds\right]+\delta \theta_{2}\, \mathbb{E} \, \left[  \int_{0}^{T}\left\vert R^{\ast}\widehat{Y}_{s}\right\vert ^{2} ds\right]\\
&& \hspace{0.1in}
+\, \delta \, \theta_{2}\, \mathbb{E} \, \left[  \int_{0}^{T}\left\vert R^{\ast}\widehat{\overline{Y}}_{s}\right\vert ^{2} ds\right]+\delta \, \mathbb{E} \, \left[  \int_{0}^{T}\| R\widehat{z}_{s}\| ^{2} ds\right]+\delta \theta_{2}\, \mathbb{E} \, \left[  \int_{0}^{T}\| R^{\ast}\widehat{\overline{Z}}_{s}\| ^{2} ds\right]\\
&& \hspace{0.1in}+\, \delta \, \theta_{2}\, \mathbb{E} \, \left[  \int_{0}^{T}\| R^{\ast}\widehat{Z}_{s}\| ^{2} ds\right] +\delta \, \mathbb{E} \, \left[  \int_{0}^{T}|||R^{\ast} \widehat{k}_{s}||| ^{2} ds\right]\\
&& \hspace{1.5in}+ \, \delta \, \theta_{2}\, \mathbb{E} \, \left[  \int_{0}^{T}||| R^{\ast}\widehat{\overline{k}}_{s}||| ^{2} ds\right] + \delta \, \theta_{2}\, \mathbb{E} \, \left[  \int_{0}^{T}|||R^{\ast} \widehat{k}_{s}||| ^{2} ds\right] .
\end{eqnarray*}

We conclude
\begin{eqnarray}\label{eq:3.37}
&& \hspace{-0.5cm}(1-\alpha_{0}+\alpha_{0}\beta_{2})\, \mathbb{E} \, \left[  \left\vert R^{\ast}\widehat{Y}_{0}\right\vert
^{2}\right] +\theta_{2}\, \mathbb{E} \, \left[\int
_{0}^{T}\left(  \left\vert R^{\ast}\widehat{Y}_{s}\right\vert ^{2}+\left\Vert R^{\ast}\widehat{Z}_{s}\right\Vert^{2}+|||R^{\ast} \widehat{k}_{s}||| ^{2}\right)  ds \right]
\nonumber \\
&& \hspace{0.75cm} +\alpha_{0}\theta_{1}\, \mathbb{E} \, \left[  \int_{0}^{T}\left(  \left\vert
R\widehat{y}_{s}\right\vert ^{2}+\left\Vert R\widehat{z}_{s}\right\Vert
^{2}\right)  ds\right]+\alpha_{0}\beta_{1}  \, \mathbb{E} \, \left[
\left\vert R\widehat{y}_{T}\right\vert ^{2}\right]
\nonumber \\ &&
\leq \delta C\, \mathbb{E} \, \left[  \left\vert \widehat{Y}%
_{0}\right\vert ^{2}\right] +\delta C\, \mathbb{E} \, \left[  \left\vert R^{\ast}\widehat{\overline{Y}}%
_{0}\right\vert ^{2}\right] +\delta  \, \mathbb{E} \, \left[  \left\vert R\widehat{y}%
_{T}\right\vert ^{2}\right]  +\delta C\, \mathbb{E} \, \left[  \left\vert \widehat{\overline{y}
}_{T}\right\vert ^{2}\right]   \nonumber \\
&& \hspace{6cm}+\, \delta \, C \,\mathbb{E}\, \left[\int_{0}^{T}\left(  \left\Vert \widehat{\upsilon}%
_{s}\right\Vert^{2}+\left\Vert \widehat{\overline{\upsilon}}%
_{s}\right\Vert^{2}\right)  ds\right],
\end{eqnarray}
where $C$ is a universal constant depending on $R^{\ast},R,\theta_{2}$ and $c.$

\medskip

Now we want to find an estimate for $\mathbb{E} \left[  \left\vert \widehat{y}%
_{T}\right\vert ^{2}\right] $ or $\mathbb{E} \left[  \left\vert R\widehat{y}%
_{T}\right\vert ^{2}\right] $ by using It\^{o}'s formula (Proposition~\ref{Lemma 3.3}) to $\left\vert R\widehat{y}_{s}\right\vert ^{2}$ over
$\left[  0,t\right], $
\begin{eqnarray*}
&& \hspace{-0.5cm}\mathbb{E}\, \left[  \left\vert \widehat{y}_{t}\right\vert ^{2}\right]
+\mathbb{E}\, \left[  \int_{0}^{t}\left\Vert \widehat{z}_{s}\right\Vert
^{2}ds\right] \\
&&  \hspace{-0.25cm}= \, \mathbb{E}\, \left[  \left\vert \alpha_{0}\widehat{\Psi}\left(
Y_{0}\right) -\left(  1-\alpha_{0}\right)R^{\ast}\widehat{Y}_{0} +\delta \left( \widehat{\Psi}\left(  \overline{Y}_{0}\right)- R^{\ast}\widehat{\overline{Y}}_{0}\right)\right\vert ^{2}\right]\\
&&  \hspace{-0.25cm}+\, 2\, \mathbb{E}\, \left[  \int_{0}^{t} \left\langle \widehat{y}_{s},
\alpha_{0} \widehat{b}\left(  s,\upsilon
_{s}\right) -\left(  1-\alpha_{0}\right)  \theta_{2}R^{\ast}\widehat{Y}_{s} +\delta\left( \widehat{b}\left(
s,\overline{\upsilon}_{s}\right) +\theta_{2}R^{\ast}\widehat{\overline{Y}}_{s} \right) \right\rangle ds\right] \\
&& \hspace{-0.25cm}+\, \mathbb{E} \, \left[\int_{0}^{t}\left\vert \alpha_{0}\widehat{\sigma}\left(
s,\upsilon_{s}\right)  -\left(  1-\alpha_{0}\right)  \theta_{2}R^{\ast}\widehat{Z}_{s}%
+\delta\left(  \widehat{\sigma}\left(  s,\overline{\upsilon}_{s}\right)  +\theta
_{2}R^{\ast}\widehat{\overline{Z}}_{s}\right)  \right\vert ^{2}ds\right]\\
&& \hspace{-0.25cm}+\, \mathbb{E} \, \left[\int_{0}^{t}\int_{\Theta}\left\vert \alpha_{0}\widehat{\varphi}\left(
s,\upsilon_{s},\rho\right)  -\left(  1-\alpha_{0}\right)  \theta_{2}R^{\ast}\widehat{k}_{s}(\rho) +\delta\left(  \widehat{\varphi}\left(  s,\overline{\upsilon}_{s},\rho\right) \right. \right. \right. \\ && \hspace{3.2in} \left. \left. \left. + \, \theta
_{2}R^{\ast}\widehat{\overline{k}}_{s}(\rho) \right)  \right\vert ^{2}\Pi_{s}(d\rho)ds\right].
\end{eqnarray*}
Therefore
\begin{eqnarray*}
&& \hspace{-0.5cm}\mathbb{E}\, \left[  \left\vert \widehat{y}_{t}\right\vert ^{2}\right]
+\mathbb{E}\, \left[  \int_{0}^{t}\left\Vert \widehat{z}_{s}\right\Vert
^{2}ds\right]\\
&&  \hspace{-0.25cm} \leq  4 \, \mathbb{E}\, \left[ \alpha_{0}^{2}\left\vert \widehat{\Psi}\left(
Y_{0}\right)\right\vert ^{2} +\left(  1-\alpha_{0}\right)^{2}\left\vert R^{\ast}\widehat{Y}_{0}\right\vert ^{2} +\delta ^{2} \left\vert \widehat{\Psi}\left(  \overline{Y}_{0}\right)\right\vert ^{2}+\delta ^{2} \left\vert R^{\ast}\widehat{\overline{Y}}_{0}\right\vert ^{2}\right]\\
&&  \hspace{-0.25cm}+ \, 2\, \mathbb{E}\, \left[  \int_{0}^{t}\left\vert \widehat{y}_{s}\right\vert
\cdot\left(  \alpha_{0}\left\vert \widehat{b}\left(  s,\upsilon
_{s}\right)  \right\vert +\left(  1-\alpha_{0}\right)  \theta_{2}\left\vert
R^{\ast}\widehat{Y}_{s}\right\vert +\delta\left\vert \widehat{b}\left(
s,\overline{\upsilon}_{s}\right)  \right\vert +\delta\theta_{2}\left\vert
R^{\ast}\widehat{\overline{Y}}_{s}\right\vert \right)  ds\right] \\
&& \hspace{1.5in}+\, \mathbb{E} \, \left[\int_{0}^{t}(1+\varepsilon)\alpha_{0}^{2}\left\vert \widehat{\sigma}\left(
s,\upsilon_{s}\right)  \right\vert ^{2}ds\right]\\
&& \hspace{-0.25cm}+\, \mathbb{E} \, \left[\int_{0}^{t}(1+\frac{1}{\varepsilon})\left\vert \left(  1-\alpha_{0}\right)  \theta_{2}R^{\ast}\widehat{Z}_{s}%
+\delta\left(  \widehat{\sigma}\left(  s,\overline{\upsilon}_{s}\right)  +\theta
_{2}R^{\ast}\widehat{\overline{Z}}_{s}\right)  \right\vert ^{2}ds\right] \\
&& \hspace{1.5in}+\, \mathbb{E} \, \left[\int_{0}^{t}\int_{\Theta}(1+\varepsilon')\alpha_{0}^{2}\left\vert \widehat{\varphi}\left(
s,\upsilon_{s},\rho\right)  \right\vert ^{2}\Pi_{s}(d\rho)ds\right]\\
&& \hspace{-0.25cm}+\, \mathbb{E} \, \left[\int_{0}^{t}\int_{\Theta}(1+\frac{1}{\varepsilon'})\left\vert \left(  1-\alpha_{0}\right)  \theta_{2}R^{\ast}\widehat{k}_{s}(\rho)%
+\delta\left(  \widehat{\varphi}\left(  s,\overline{\upsilon}_{s},\rho\right)  +\theta
_{2}R^{\ast}\widehat{\overline{k}}_{s}(\rho)\right)  \right\vert ^{2}\Pi_{s}(d\rho)ds\right],
\end{eqnarray*}
for any $\varepsilon >0$ and $\varepsilon' >0.$ Choose then \ $\varepsilon'=\varepsilon=\frac{1-\gamma}{2\gamma}$ to get
\begin{eqnarray}\label{eq:3.38}
&& \hspace{-0,5cm}\mathbb{E}\, \left[  \left\vert \widehat{y}_{t}\right\vert ^{2}\right]
+\mathbb{E}\, \left[  \int_{0}^{t}\left\Vert \widehat{z}_{s}\right\Vert
^{2}ds\right] \nonumber \\
&& \leq 4\, \mathbb{E}\, \left[ \alpha_{0}^{2}\left\vert \widehat{\Psi}\left(
Y_{0}\right)\right\vert ^{2} +\left(  1-\alpha_{0}\right)^{2}\left\vert R^{\ast}\widehat{Y}_{0}\right\vert ^{2} +\delta ^{2} \left\vert \widehat{\Psi}\left(  \overline{Y}_{0}\right)\right\vert ^{2}+\delta ^{2} \left\vert R^{\ast}\widehat{\overline{Y}}_{0}\right\vert ^{2}\right]\nonumber \\
&&
+\, 2\, \mathbb{E}\, \left[  \int_{0}^{t}\left\vert \widehat{y}_{s}\right\vert
\cdot\left(  \alpha_{0}\left\vert \widehat{b}\left(  s,\upsilon
_{s}\right)  \right\vert +\left(  1-\alpha_{0}\right)  \theta_{2}\left\vert
R^{\ast}\widehat{Y}_{s}\right\vert +\delta\left\vert \widehat{b}\left(
s,\overline{\upsilon}_{s}\right)  \right\vert +\delta\theta_{2}\left\vert
R^{\ast}\widehat{\overline{Y}}_{s}\right\vert \right)  ds\right]
\nonumber \\
&& \hspace{1.5in} +\, \left(\frac{1+\gamma}{2\gamma}\right) \, \alpha_{0}^{2}\, \mathbb{E}\, \left[  \int_{0}^{t}\left\vert \widehat{\sigma}\left(  s,\upsilon_{s}\right)  \right\vert
^{2}ds\right] \nonumber\\
 && +\, 3\, \left(  \frac{1+\gamma}{1-\gamma}\right) \, \mathbb{E}\, \left[  \int_{0}^{t} \int_{\Theta}\left(  \left(  1-\alpha_{0}\right)  ^{2}\theta_{2}^{2}| R^{\ast}\widehat{k}_{s}(\rho)|^{2}+\delta^{2}\left\vert \widehat{\varphi}\left(
s,\overline{\upsilon}_{s},\rho\right)  \right\vert ^{2}\right. \right. \nonumber \\
&& \hspace{3.5in}+\, \left.  \left. \delta^{2}\theta_{2}%
^{2}| R^{\ast}\widehat{\overline{k}}_{s}(\rho)| ^{2}\right)\Pi_{s}(d\rho)  ds\right] \nonumber \\
&& \hspace{1.5in} +\, \left(\frac{1+\gamma}{2\gamma}\right) \alpha_{0}^{2} \, \mathbb{E}\, \left[  \int_{0}^{t} \int_{\Theta}\left\vert \widehat{\varphi}\left(  s,\upsilon_{s},\rho\right)  \right\vert
^{2}\Pi_{s}(d\rho)ds\right]
\nonumber\\
&& +\, 3\, \left(  \frac{1+\gamma}{1-\gamma
}\right) \mathbb{E}\, \left[  \int_{0}^{t} \left(  \left(  1-\alpha_{0}\right)  ^{2}\theta_{2}^{2}\left\Vert R^{\ast}\widehat
{Z}_{s}\right\Vert^{2}+\delta^{2}\left\vert \widehat{\sigma}\left(
s,\overline{\upsilon}_{s}\right)  \right\vert ^{2}+\delta^{2}\theta_{2}^{2}\left\Vert R^{\ast}\widehat{\overline{Z}}_{s}\right\Vert^{2}\right)  ds\right] \nonumber \\
&& =:I'_{1}+I'_{2}+I'_{3}+I'_{4}+I'_{5}+I'_{6}.
\end{eqnarray}
We have
\begin{eqnarray}\label{eq:3.39}
&& \hspace{-0.5cm} I'_{1}  =  4\, \mathbb{E}\, \left[ \alpha_{0}^{2}\left\vert \widehat{\Psi}\left(
Y_{0}\right)\right\vert^{2} +\left(  1-\alpha_{0}\right)^{2}\left\vert R^{\ast}\widehat{Y}_{0}\right\vert ^{2} +\delta ^{2} \left\vert \widehat{\Psi}\left(  \overline{Y}_{0}\right)\right\vert ^{2}+\delta ^{2} \left\vert R^{\ast}\widehat{\overline{Y}}_{0}\right\vert ^{2}\right] \nonumber \\
&& \leq 4\, \mathbb{E}\, \left[ \alpha_{0}^{2}c\left\vert \widehat{Y}_{0}\right\vert ^{2} +\left(  1-\alpha_{0}\right)^{2}\left\vert R^{\ast}\widehat{Y}_{0}\right\vert ^{2} +\delta ^{2}c \left\vert \widehat{\overline{Y}}_{0}\right\vert ^{2}+\delta ^{2} \left\vert R^{\ast}\widehat{\overline{Y}}_{0}\right\vert ^{2}\right] \nonumber \\
&& \leq C~ \mathbb{E}\, \left[ \left\vert \widehat{Y}_{0}\right\vert ^{2} +\delta \left\vert \widehat{\overline{Y}}_{0}\right\vert ^{2}\right],
\end{eqnarray}
It follows from the fact that $b$ is Lipschitz
\begin{eqnarray*}
&& \hspace{-0.5cm} I'_{2}=2\, \mathbb{E}\, \left[  \int_{0}^{t}\left\vert \widehat{y}_{s}\right\vert
\cdot\left(  \alpha_{0}\left\vert \widehat{b}\left(  s,\upsilon
_{s}\right)  \right\vert +\left(  1-\alpha_{0}\right)  \theta_{2}\left\vert
R^{\ast}\widehat{Y}_{s}\right\vert +\delta\left\vert \widehat{b}\left(
s,\overline{\upsilon}_{s}\right)  \right\vert +\delta\theta_{2}\left\vert
R^{\ast}\widehat{\overline{Y}}_{s}\right\vert \right)  ds\right]  \\
&& \leq\mathbb{E}\, \left[  \int_{0}^{t}\left(  \left(  \frac{8c\;\alpha_{0}}{1-\gamma}\right)  \left\vert \widehat{y}_{s}\right\vert ^{2}+\frac{1-\gamma}{8c}\, c\left\Vert \widehat{\upsilon}_{s}\right\Vert^{2} +\left(  1-\alpha_{0}\right)  \theta_{2}\left(  \left\vert
\widehat{y}_{s}\right\vert ^{2}+\left\vert R^{\ast}\widehat{Y}_{s}\right\vert ^{2}\right)
\right. \right. \\
&& \left. \left. \hspace{3.5cm}+ \, \delta \left(  \left\vert \widehat{y}_{s}\right\vert^{2}+c\left\Vert
\widehat{\overline{\upsilon}}_{s} \right\Vert^{2}\right)  +\delta\theta_{2}\left(  \left\vert \widehat{y}_{s}\right\vert
^{2}+\left\vert R^{\ast}\widehat{\overline{Y}}_{s}\right\vert^{2}\right)  \right) ds\right]  \\
&& \leq\mathbb{E}\, \left[  \int_{0}^{t}\left\vert \widehat{y}_{s}\right\vert ^{2}\left(
\frac{8c\;\alpha_{0}}{1-\gamma}+\frac{1-\gamma}{8}+\left(  1-\alpha_{0}\right)
\theta_{2}+\delta+\delta\theta_{2}\right)  ds\right]  \\
&& +\, \mathbb{E} \, \left[  \int_{0}^{t}\left\vert \widehat{Y}_{s}\right\vert ^{2}\left(
\frac{1-\gamma}{8}\right)  ds\right]  +\mathbb{E}\, \left[  \int_{0}^{t}\left\vert R^{\ast}\widehat{Y}_{s}\right\vert ^{2}\left(  \left(  1-\alpha
_{0}\right)  \theta_{2}\right)  ds\right]   \\
&&  +\, \mathbb{E} \, \left[  \int_{0}^{t}||| \widehat{k}_{s}|||^{2}\left(  \frac{1-\gamma}{8}\right)  ds\right]+\mathbb{E}\, \left[  \int_{0}^{t}\left\Vert \widehat{z}_{s}\right\Vert^{2}\left(
\frac{1-\gamma}{8}\right)  ds\right] \\
&&  +\, \mathbb{E}\, \left[  \int_{0}%
^{t}\left\Vert \widehat{Z}_{s}\right\Vert^{2}\left(  \frac{1-\gamma}{8}\right)
ds\right] +\, \mathbb{E} \, \left[  \int_{0}%
^{t}\left\vert R^{\ast}\widehat{\overline{Y}}_{s}\right\vert ^{2}\left(  \delta\theta
_{2}\right)  ds\right] +
\mathbb{E}\, \left[  \int_{0}^{t}\left\Vert \widehat{\overline{\upsilon}}_{s}\right\Vert^{2}(c\, \delta) ds\right]  .
\end{eqnarray*}
Therefore
\begin{eqnarray}\label{eq:3.40}
&& \hspace{-0.5cm} I'_{2}   \leq C\, \mathbb{E} \, \left[  \int_{0}^{t}\left(\left\vert \widehat{y}%
_{s}\right\vert ^{2}+\left\vert
\widehat{Y}_{s}\right\vert ^{2}ds\right)\right] \nonumber \\
&&  + \left(\frac{1-\gamma}{8}\right) \, \mathbb{E} \, \left[  \int_{0}^{t}\left(  \left\Vert \widehat
{z}_{s}\right\Vert^{2}+\left\Vert \widehat{Z}_{s}\right\Vert^{2}+||| \widehat{k}_{s}||| ^{2}\right)  ds\right] + \, \delta C\, \mathbb{E} \, \left[  \int_{0}%
^{t}\left\Vert \widehat{\overline{\upsilon}}_{s} \right\Vert^{2}ds\right]. \nonumber \\
\end{eqnarray}
Similarly,
\begin{eqnarray}\label{eq:3.41}
&& \hspace{-0.5cm} I'_{3}   =\left(  \frac{1+\gamma}{2\gamma
}\right)  \alpha_{0}^{2}\, \mathbb{E} \, \left[  \int_{0}^{t}\left\vert \widehat{\sigma}\left(  s,\upsilon
_{s}\right)  \right\vert ^{2}ds\right] \nonumber \\ &&
\leq \left(
\frac{1+\gamma}{2\gamma}\right)  \alpha_{0}^{2}\, \mathbb{E} \, \left[  \int_{0}^{t}\left(  c\left\vert \widehat{y}_{s}\right\vert ^{2}+c\left\vert \widehat{Y}%
_{s}\right\vert ^{2}+\frac{\gamma}{2}\left\Vert \widehat{z}_{s}\right\Vert^{2}+c\left\Vert \widehat{Z}_{s}\right\Vert^{2}+c||| \widehat{k}_{s}||| ^{2}\right)  ds\right] \nonumber \\
&& \leq C\, \mathbb{E} \, \left[  \int_{0}^{t}\left(  \left\vert \widehat{y}%
_{s}\right\vert ^{2}+\left\vert \widehat{Y}_{s}\right\vert ^{2}+\left\Vert \widehat
{Z}_{s}\right\Vert^{2}+||| \widehat{k}_{s}|||^{2}\right)  ds\right] \nonumber \\
&& \hspace{6cm}  + \left(
\frac{1+\gamma}{4}\right)  \alpha_{0}^{2} \, \mathbb{E} \, \left[ \int_{0}^{t} \left\Vert \widehat{z}%
_{s}\right\Vert^{2}  ds\right],
\end{eqnarray}

\medskip

\begin{eqnarray}\label{eq:3.42}
&& \hspace{-0,5cm}I'_{4}=3\left(  \frac{1+\gamma}{1-\gamma}\right)\, \mathbb{E} \, \left[  \int_{0}^{t}\int_{\Theta}\left(  \left(
1-\alpha_{0}\right)^{2}  \theta_{2}^{2}\left\vert R^{\ast}\widehat{k}_{s}(\rho)\right\vert^{2}+\delta^{2} \left\vert\widehat{\varphi%
}\left(  s,\overline{\upsilon}_{s},\rho\right)\right\vert^{2}  \right.  \right.\nonumber\\
&& \hspace{3in} \left.  \left.
+\, \delta^{2}\theta_{2}^{2}\left\vert R^{\ast}\widehat{\overline{k}}%
_{s}(\rho) \right\vert ^{2}\right) \Pi_{s}(d\rho) ds \right]\nonumber\\
&& \hspace{0.5cm}\leq C\, \mathbb{E} \, \left[  \int_{0}%
^{t}\left(||| \widehat{k}_{s}|||
^{2} +\delta\left\Vert \widehat{\overline{\upsilon}}_{s}\right\Vert^{2}\right)ds\right] ,
\end{eqnarray}
and
\bigskip
\begin{eqnarray}\label{eq:3.43}
&& \hspace{-0.5cm} I'_{5}   =\left(  \frac{1+\gamma}{2\gamma
}\right)  \alpha_{0}^{2}\; \mathbb{E} \, \left[  \int_{0}^{t}\int_{\Theta}\left\vert \widehat{\varphi}\left(  s,\upsilon
_{s},\rho\right)  \right\vert ^{2}\Pi_{s}(d\rho)ds\right] \nonumber \\ &&
\leq \left(
\frac{1+\gamma}{2\gamma}\right)  \alpha_{0}^{2}\; \mathbb{E} \, \left[  \int_{0}^{t}\left(  c\left\vert \widehat{y}_{s}\right\vert ^{2}+c\left\vert \widehat{Y}%
_{s}\right\vert ^{2}+\frac{\gamma}{2}\left\Vert \widehat{z}_{s}\right\Vert^{2}+c\left\Vert \widehat{Z}_{s}\right\Vert^{2}+c||| \widehat{k}_{s}||| ^{2}\right)  ds\right] \nonumber \\
&& \leq C\, \mathbb{E} \, \left[  \int_{0}^{t}\left(  \left\vert \widehat{y}%
_{s}\right\vert ^{2}+\left\vert \widehat{Y}_{s}\right\vert ^{2}+\left\Vert \widehat
{Z}_{s}\right\Vert^{2}+||| \widehat{k}_{s}||| ^{2}\right)  ds\right]  \nonumber \\
&& \hspace{6cm} + \left(
\frac{1+\gamma}{4}\right)  \alpha_{0}^{2} \; \mathbb{E} \, \left[ \int_{0}^{t} \left\Vert \widehat{z}%
_{s}\right\Vert^{2}  ds\right].
\end{eqnarray}
Also,
\begin{eqnarray}\label{eq:3.44}
&& \hspace{-0.5cm}  I'_{6}   =3\, \left(  \frac{1+\gamma
}{1-\gamma}\right)\, \mathbb{E} \, \left[  \int_{0}^{t} \left(   \left(  1-\alpha_{0}\right)^{2}\theta_{2}^{2}\left\Vert
R^{\ast}\widehat{Z}_{s}\right\Vert^{2}+\delta^{2}\left\vert \widehat{\sigma}\left(
s,\overline{\upsilon}_{s}\right)\right\vert ^{2}+\delta^{2}\theta_{2}%
^{2}\left\Vert R^{\ast}\widehat{\overline{Z}}_{s}\right\Vert^{2} \right)  ds\right]\nonumber \\
&& \leq 3 \, \left(  \frac{1+\gamma}{1-\gamma
}\right)\, \mathbb{E} \, \left[  \int_{0}^{t}  \left(\left(  1-\alpha_{0}\right)^{2}\theta_{2}^{2}\left\Vert R^{\ast}\widehat
{Z}_{s}\right\Vert^{2}+\delta^{2}\theta_{2}^{2}\left\Vert R^{\ast}\widehat{\overline{Z}}_{s}\right\Vert^{2} \right) ds\right]
\nonumber\\
&& \hspace{3cm} +\, 3\, \delta^{2}\left(  \frac{1+\gamma
}{1-\gamma}\right) \, c \, \mathbb{E} \, \left[  \int_{0}^{t}  \left( \left\vert \widehat{\overline{y}}_{s}\right\vert ^{2}%
+\left\vert
\widehat{\overline{Y}}_{s}\right\vert ^{2} + ||| \widehat{\overline{k}}_{s}|||^{2}+\left\Vert \widehat{\overline{Z}}_{s}\right\Vert^{2} \right) ds\right]\nonumber\\
&& \hspace{3.5cm}+\, 3\, \delta^{2}\, \left(  \frac{1+\gamma
}{1-\gamma}\right) \frac{\gamma}{2} \, \mathbb{E} \, \left[  \int_{0}^{t} \left\Vert \widehat{\overline{z}}_{s}\right\Vert^{2}ds\right]\nonumber\\
&& \leq C \, \mathbb{E} \, \left[  \int_{0}^{t}\left(  \left\Vert R^{\ast}\widehat{Z}%
_{s}\right\Vert^{2}+\delta\left\Vert R^{\ast}\widehat{\overline{Z}}_{s}\right\Vert
^{2}+\delta \gamma\left\Vert \widehat{\overline{z}}_{s}\right\Vert^{2}\right) ds\right]\nonumber\\
&& \hspace{3.5cm}
+\, \delta C\, \mathbb{E} \, \left[  \int_{0}^{t}  \left(\left\vert
\widehat{\overline{y}}_{s}\right\vert ^{2}+\left\vert \widehat{\overline{Y}}_{s}\right\vert ^{2}+\left\Vert \widehat{\overline{Z}}_{s}\right\Vert^{2}+||| \widehat{\overline{k}}_{s}||| ^{2}\right)  ds\right] \nonumber\\
&& \leq C \, \mathbb{E} \, \left[  \int_{0}^{t}\left(  \left\Vert R^{\ast}\widehat{Z}%
_{s}\right\Vert^{2}+\delta\left\Vert\widehat{\overline{\upsilon}}_{s}\right\Vert^{2}\right) ds\right].
\end{eqnarray}

Now Substitute (\ref{eq:3.39})-(\ref{eq:3.44}) in (\ref{eq:3.38}) to get%
\begin{eqnarray*}
&& \hspace{-0.5cm} \mathbb{E}\, \left[  \left\vert \widehat{y}_{t}\right\vert ^{2}\right]
+\mathbb{E}\, \left[  \int_{0}^{t}  \left\Vert \widehat{z}_{s}\right\Vert
^{2}ds\right] \leq C\, \mathbb{E} \, \left[ \left\vert \widehat{Y}_{0}\right\vert ^{2} +\delta \left\vert \widehat{\overline{Y}}_{0}\right\vert ^{2}\right]\nonumber \\
&& \hspace{2in} + \, C\, \mathbb{E} \, \left[  \int_{0}^{t}\left(\left\vert \widehat{y}%
_{s}\right\vert ^{2}+\left\vert
\widehat{Y}_{s}\right\vert ^{2}\right)ds\right] \nonumber \\
&& \hspace{1in} + \; \left(\frac{1-\gamma}{8}\right)\, \mathbb{E} \, \left[  \int_{0}^{t}\left(  \left\Vert \widehat
{z}_{s}\right\Vert^{2}+\left\Vert \widehat{Z}_{s}\right\Vert^{2}+\||| \widehat{k}_{s}||| ^{2}\right)  ds\right]\nonumber \\ &&
\hspace{3in}   + \, C\delta\, \mathbb{E} \, \left[  \int_{0}%
^{t}\left\Vert \widehat{\overline{\upsilon}}_{s}\right\Vert^{2}ds\right]  \nonumber \\
&& \hspace{-0.5cm} + \, C\, \mathbb{E} \, \left[  \int_{0}^{t}\left(  \left\vert \widehat{y}%
_{s}\right\vert^{2}+\left\vert \widehat{Y}_{s}\right\vert ^{2}+\left\Vert \widehat
{Z}_{s}\right\Vert^{2}+||| \widehat{k}_{s}|||
^{2}\right)  ds\right]  + \left(
\frac{1+\gamma}{4}\right)  \alpha_{0}^{2} \, \mathbb{E} \, \left[ \int_{0}^{t} \left\Vert \widehat{z}%
_{s}\right\Vert^{2}  ds\right]  \nonumber\\
&&   \hspace{2in} + \, C\, \mathbb{E} \, \left[  \int_{0}%
^{t}\left(||| \widehat{k}_{s}|||
^{2} +\delta\left\Vert \widehat{\overline{\upsilon}}_{s}\right\Vert^{2}\right)ds\right] \nonumber \\
&& \hspace{-0.5cm} + \,C\, \mathbb{E} \, \left[  \int_{0}^{t}\left(  \left\vert \widehat{y}%
_{s}\right\vert^{2}+\left\vert \widehat{Y}_{s}\right\vert ^{2}+\left\Vert \widehat
{Z}_{s}\right\Vert^{2}+||| \widehat{k}_{s}|||
^{2}\right)  ds\right]  + \left(
\frac{1+\gamma}{4}\right)  \alpha_{0}^{2} \, \mathbb{E} \, \left[ \int_{0}^{t} \left\Vert \widehat{z}%
_{s}\right\Vert^{2}  ds\right]  \nonumber\\
&& \hspace{2.5in} + \,C \, \mathbb{E} \, \left[  \int_{0}^{t}\left(  \left\Vert R^{\ast}\widehat{Z}%
_{s}\right\Vert^{2}+\delta\left\Vert\widehat{\overline{\upsilon}}_{s}\right\Vert^{2}\right) ds\right]  ,
\end{eqnarray*}
which implies%
\begin{eqnarray*}\label{eq:3.45}
&& \hspace{-1cm} \mathbb{E}\, \left[  \left\vert \widehat{y}_{t}\right\vert^{2}\right]
+\mathbb{E}\, \left[  \int_{0}^{t}  \left\Vert \widehat{z}_{s}\right\Vert^{2}
ds\right] \nonumber \\
&&\leq \, C\, \mathbb{E} \, \left[ \left\vert \widehat{Y}_{0}\right\vert^{2} +\delta \left\vert \widehat{\overline{Y}}_{0}\right\vert ^{2}\right]
+ C\, \mathbb{E} \, \left[  \int_{0}^{t}\left\vert \widehat{y}_{s}\right\vert ^{2}ds\right]  +C\delta\, \mathbb{E} \, \left[  \int_{0}%
^{t}\left\Vert \widehat{\overline{\upsilon}}_{s}\right\Vert^{2}ds\right]  \nonumber\\
&& \hspace{2cm} +\, C\, \mathbb{E} \, \left[  \int_{0}^{t}\left(  \left\vert \widehat{Y}_{s}\right\vert ^{2}+\left\Vert \widehat{Z}_{s}\right\Vert^{2}+||| \widehat{k}_{s}|||^{2}\right)  ds\right]
 \nonumber \\
&& \hspace{2cm}  +\, \left(  \left(\frac{1-\gamma}{8}\right) +\left(  \frac{1+\gamma}{4}\right) + \left(  \frac{1+\gamma}{4}\right)\right)  \, \mathbb{E} \, \left[  \int_{0}^{t} \left\Vert \widehat{z}_{s}\right\Vert^{2} ds\right].
\end{eqnarray*}
In particular, we have
\begin{eqnarray}\label{eq:3.46}
&& \hspace{-1.5cm} \mathbb{E}\, \left[  \left\vert \widehat{y}_{t}\right\vert^{2}\right]
+\left(1-\frac{5+3\gamma}{8}\right) \, \mathbb{E}\, \left[  \int_{0}^{t}  \left\Vert \widehat{z}_{s}\right\Vert
^{2}
ds\right] \nonumber \\
&&\leq C\, \mathbb{E} \, \left[ \left\vert \widehat{Y}_{0}\right\vert^{2} +\delta \left\vert \widehat{\overline{Y}}_{0}\right\vert ^{2}\right]
+C\, \mathbb{E}\, \left[  \int_{0}^{t}\left\vert \widehat{y}%
_{s}\right\vert ^{2}ds\right]  + C \, \delta\, \mathbb{E}\, \left[  \int_{0}^{t}\left\Vert \widehat{\overline{\upsilon}}_{s}\right\Vert^{2}ds\right]  \nonumber\\
&& \hspace{4cm}+ \, C\, \mathbb{E}\, \left[  \int_{0}^{t}\left(  \left\vert \widehat{Y}%
_{s}\right\vert ^{2}+\left\Vert \widehat{Z}_{s}\right\Vert^{2}+||| \widehat{k}_{s}||| ^{2}\right)  ds\right]
.
\end{eqnarray}

Now apply Gronwall's inequality to get
\begin{eqnarray*}\label{eq:3.47}
&& \hspace{-2cm}\mathbb{E} \left[  \left\vert \widehat{y}_{t}\right\vert ^{2}\right]  \leq
C^{2}e^{t}\left(   \mathbb{E}\, \left[  \left\vert \widehat{Y}_{0}\right\vert
^{2}\right]  +\delta\, \mathbb{E} \, \left[  \left\vert \widehat{\overline{Y}}%
_{0}\right\vert ^{2}\right]  \right.\nonumber \\ &&
\hspace{1.5cm}
\left. + \, \mathbb{E}\, \left[  \int_{0}^{t}\left(
\left\vert \widehat{Y}_{s}\right\vert ^{2}+\left\Vert \widehat{Z}_{s}\right\Vert
^{2}+||| \widehat{k}_{s}|||
^{2}+\delta\left\Vert \widehat{\overline{\upsilon}}_{s}\right\Vert^{2}\right)  ds\right]  \right)  ,
\end{eqnarray*}
for all $0\leq t\leq T,$ and
\begin{eqnarray*}\label{eq:3.48}
&& \hspace{-2cm}\mathbb{E}\, \left[  \int_{0}^{t}  \left\Vert \widehat{z}_{s}\right\Vert^{2}
ds\right]  \leq C \left(    \mathbb{E}\, \left[  \left\vert \widehat{Y}_{0}\right\vert^{2}\right]  + \delta\, \mathbb{E} \, \left[  \left\vert \widehat{\overline{Y}}_{0}\right\vert^{2}\right] \right) + C t\, \mathbb{E} \, \left[  \left\vert \widehat{y}_{T}\right\vert
^{2}\right] \nonumber \\ &&
\hspace{1.5cm}
+ \, C\, \mathbb{E} \, \left[  \int_{0}^{T}\left(\left\vert \widehat{Y}_{s}\right\vert^{2}+\left\Vert \widehat{Z}_{s}\right\Vert^{2}+||| \widehat{k}_{s}|||^{2}+\delta\left\Vert \widehat{\overline{\upsilon}}_{s}\right\Vert^{2}\right)  ds\right].
\end{eqnarray*}
Consequently
\begin{eqnarray}\label{eq:3.49}
&& \hspace{-2cm}\mathbb{E}\, \left[  \left\vert R\widehat{y}_{T}\right\vert^{2}\right] +\mathbb{E}\, \left[  \int_{0}^{t}  \left\Vert \widehat{z}_{s}\right\Vert
^{2}
ds\right]  \leq
C\left(   \mathbb{E}\, \left[  \left\vert \widehat{Y}_{0}\right\vert
^{2}\right]  +\delta\, \mathbb{E} \, \left[  \left\vert \widehat{\overline{Y}}%
_{0}\right\vert ^{2}\right] \right) \nonumber \\ &&
\hspace{1.5cm}
 + \, C\, \mathbb{E} \, \left[  \int_{0}^{t}\left(
\left\vert \widehat{Y}_{s}\right\vert ^{2}+\left\Vert \widehat{Z}_{s}\right\Vert
^{2}+||| \widehat{k}_{s}|||
^{2}+\delta\left\Vert \widehat{\overline{\upsilon}}_{s}\right\Vert ^{2}\right)  ds\right]  ,
\end{eqnarray}
with the help of Remark~\ref{Remark 3.2}.

Now we can argue as in the proof of Lemma~\ref{Lemma: 3.5} and in particular in parallel to the part following inequality (\ref{eq:3.19}) to conclude (with the help of Remark~\ref{Remark 3.2}) that
\begin{eqnarray*}
&& \hspace{-0.5cm} \mathbb{E}\, \left[  \left\vert \widehat{y}%
_{T}\right\vert ^{2}\right]  +\mathbb{E}\, \left[  \left\vert \widehat{Y}_{0}\right\vert ^{2}\right] + \mathbb{E}\, \left[  \int_{0}^{T}\left\Vert \widehat{\upsilon}_{s}\right\Vert^{2}ds\right] \\
&& \hspace{2cm} \leq\frac{1}{2}\left(  \mathbb{E}\, \left[  \left\vert \widehat{\overline{y}}_{T}\right\vert ^{2}\right]
+\mathbb{E}\, \left[  \left\vert \widehat{\overline{Y}}_{0}\right\vert^{2}\right] + \mathbb{E}\, \left[  \int_{0}^{T}\left\Vert \widehat
{\overline{\upsilon}}_{s}\right\Vert^{2}ds\right] \right)  ,
\end{eqnarray*}
which shows that the mapping $I_{\alpha_{0}+\delta}$,  defined in the proof of Lemma~\ref{Lemma: 3.5}, is contraction for some small $\delta \in[0,\delta_{0}]$ on  $$\mathbb{H}^{2}\times L^{2}(\Omega,\mathcal{F}_{T},\mathbb{P};\mathbb{R}^{n})\times L^{2}(\Omega,\mathcal{F}_{0},\mathbb{P};\mathbb{R}^{m}),$$ and so it attains a unique fixed point $\upsilon=\left(
y,Y,z,Z,k\right)$  in $\mathbb{H}^{2},$\ which can be seen easily to be the unique solution of (\ref{eq:3.32}) for
$\alpha=\alpha_{0}+\delta.$
\end{proof}

\bigskip
\textbf{Case~3.} $m=n.$ We assume (A1)-(A4) with $0< \gamma<1$ and $0< \gamma' \leq \gamma/2 .$ From (A1) and (A2) we only need to consider two cases:\\
1) If $\theta_{1}>0,\theta_{2}\geq0,\beta_{1}>0$, and $\beta_{2}\geq0$, we can
have the same result as Lemma~\ref{Lemma: 3.5}.
\\
2) If $\theta_{1}\geq0,\theta_{2}>0,\beta_{1}\geq0$, and $\beta_{2}>0$, we can
have the same result as Lemma~\ref{lem:final-lemma}.

\bigskip

We are now ready to complete the proof of Theorem~\ref{Propo 3.4}.

\medskip

\noindent\begin{proof}[Proof completion of Theorem~\ref{Propo 3.4}]
For \textbf{Case~1}, we know that, for each $$\psi\in L^{2}\left(
\Omega,\mathcal{F}_{0},\mathbb{P};\mathbb{R}^{n}\right)  ,\phi\in L^{2}\left(
\Omega,\mathcal{F}_{T},\mathbb{P};\mathbb{R}^{m}\right)  ,\left(  \widetilde{b}
_{0},\widetilde{f}_{0},\widetilde{\sigma}_{0},\widetilde{g}_{0},\varphi_{0}\right)
\in\mathbb{H}^{2},$$ FBDSDEJ~(\ref{eq:3.3}) has a unique solution as $\alpha=0.$ It follows from
Lemma~\ref{Lemma: 3.5} that there exists a positive constant $\delta_{0}=\delta_{0}\left(
c,\gamma,\beta_{1},\theta_{1},R,T\right)  $\ such that for any $\delta
\in\left[  0,\delta_{0}\right]  $\ and $\psi\in L^{2}\left(  \Omega
,\mathcal{F}_{0},\mathbb{P};\mathbb{R}^{n}\right)  ,\phi\in L^{2}\left(
\Omega,\mathcal{F}_{T},\mathbb{P};\mathbb{R}^{m}\right)  ,$ $\left(  \widetilde{b}
_{0},\widetilde{f}_{0},\widetilde{\sigma}_{0},\widetilde{g}_{0},\varphi_{0}\right)
\in\mathbb{H}^{2}$, (\ref{eq:3.3}) has a unique solution for $\alpha=\delta$. Since
$\delta_{0}$\ depends only on $c,\gamma,\beta_{1},\theta_{1},R$\ and $T$, we
can repeat this process $N$\ times with $1\leq N\delta_{0}<1+\delta_{0}$. In
particular, for $\alpha=1$\ with $\left(  \widetilde{b}_{0},\widetilde{f}_{0}
,\widetilde{\sigma}_{0},\widetilde{g}_{0},\varphi_{0}\right)  \equiv0,\phi\equiv
0,\psi\equiv0$, FBDSDEJ~(\ref{eq:3.1}) has a unique solution in $\mathbb{H}^{2}.$

\bigskip

For \textbf{Case~2}, we know that, for each $$\psi\in L^{2}\left(
\Omega,\mathcal{F}_{0},\mathbb{P};\mathbb{R}^{n}\right)  ,\phi\in L^{2}\left(
\Omega,\mathcal{F}_{T},\mathbb{P};\mathbb{R}^{m}\right)  ,\left(  \tilde{b}
_{0},\widetilde{f}_{0},\widetilde{\sigma}_{0},\widetilde{g}_{0},\varphi_{0}\right)
\in\mathbb{H}^{2},$$ FBDSDEJ~(\ref{eq:3.30}) has a unique solution as $\alpha=0$. It follows from
Lemma~\ref{lem:final-lemma} that there exists a positive constant $\delta_{0}=\delta_{0}\left(
c,\gamma,\beta_{2},\theta_{2},R,T\right)  $\ such that for any $\delta
\in\left[  0,\delta_{0}\right]  $\ and $\psi\in L^{2}\left(  \Omega
,\mathcal{F}_{0},\mathbb{P};\mathbb{R}^{n}\right)  ,\phi\in L^{2}\left(
\Omega,\mathcal{F}_{T},\mathbb{P};\mathbb{R}^{m}\right)  ,$ $\left(  \widetilde{b}
_{0},\widetilde{f}_{0},\widetilde{\sigma}_{0},\widetilde{g}_{0},\varphi_{0}\right)
\in\mathbb{H}^{2}$, (\ref{eq:3.30}) has a unique solution for $\alpha=\delta$.
Since $\delta_{0}$ depends only on $c,\gamma,\beta_{2},\theta_{2},R$\ and $T$,
we can repeat this process for $N$\ times with $1\leq N\delta_{0}<1+\delta
_{0}$, and then deduce as in the preceding case that FBDSDEJ~(\ref{eq:3.1}) has a unique solution in $\mathbb{H}^{2}.$

\bigskip

Similar to these cases, the desired result can be obtained in \textbf{Case~3}.
\end{proof}

\fussy


\begin{thebibliography}{99}
\bibitem {AG} A. Al-Hussein and B. Gherbal, Stochastic maximum principle for
Hilbert space valued forward-backward doubly SDEs with Poisson
jumps. 26th IFIP TC 7 Conference, CSMO 2013, Klagenfurt, Austria,
September 9--13, 2013. System modeling and optimization, (2014), 1--10.

\bibitem {Al-G-relaxed} A. Al-Hussein and B. Gherbal, Maximum principle for relaxed and strict control problems of forward-backward doubly SDEs with jumps under full and  partial information, 2018, Submitted.

\bibitem {Ant93} F. Antonelli, Backward-forward stochastic differential
equations, Ann. Appl. Probab., 3 (1993), 777--793.

\bibitem {BBP} G. Barles, R. Buckdahn and E. Pardoux, Backward stochastic differential equations and integral-partial
differential equations, Stoch. Stoch. Rep., 60, 1--2 (1997), 57--83.

\bibitem{Bao} F. Bao, Y. Cao and X. Han. Forward backward doubly stochastic differential equations
and the optimal filtering of diffusion processes, arXiv:1509.06352v3 [math.PR], 2017.

\bibitem{Bour} N. Bourbaki, \'El\'ements de Math´ematiques. I, Livre III: Topologie G\'en\'erale. Chapitre 9: Utilisations des Nombres
R\'eels en Topologie G\'en\'erale, Hermann, Paris, 1958.

\bibitem {DM} F. Delarue and S. Menozzi, A forward-backward stochastic algorithm for quasi-linear PDEs,
Ann. Appl. Probab., 16 (2006), 140--184.

\bibitem {HP} Y. Hu and S. Peng, Solution of forward-backward stochastic differential equations, Prob. Th. \& Rel. Fields, 103 (1995), 273--283.

\bibitem {MPY} J. Ma, P. Protter and J. Yong, Solving forward-backward stochastic differential equations
explicitly–a four step scheme, Prob. Th. \& Rel. Fields, 98 (1994), 339--359.

\bibitem {MY} J. Ma and J. Yong, Forward-backward stochastic differential equations and their
applications. Number 1702. Springer Science Business Media, (1999).

\bibitem {MSZ} J. Ma, J. Shen and Y. Zhao, On numerical approximations of forward-backward
stochastic differential equations, Siam J. Numer. Anal, vol. 46, no. 5 (2008), 2636--2661.

\bibitem{Parth} K.R. Parthasarathy, Probability measures on metric spaces, Academic press, 1967.

\bibitem {PP} E. Pardoux and S. Peng, Backward doubly stochastic differential
equations and system of quasilinear SPDEs, Prob. Th. \& Rel. Fields, 98, 2 (1994), 209--227.

\bibitem {PT} E. Pardoux and S. Tang, Forward-backward stochastic differential equations and quasilinear
parabolic PDEs, Prob. Th. \& Rel. Fields, 114 (1999), 123--150.

\bibitem {PW} S. Peng and Z. Wu, Fully coupled forward-backward stochastic
differential equations and applications to optimal control, SIAM J. Control
Optim, 37 (1999), 825--843.

\bibitem {PS} S. Peng and Y. Shi, A type-symmetric forward-backward stochastic
differential equations, C. R. Acad. Sci. Paris Ser. I, 336, 1 (2003), 773--778.

\bibitem {S} R. Situ, On solution of backward stochastic differential equations with jumps and applications, Stoch.
Process. Appl., 66, 2 (1997), 209–--36.

\bibitem {SL} X. Sun and Y. Lu, The property for solutions of the multi-dimensional backward doubly stochastic
differential equations with jumps, Chin. J. Appl. Probab. Stat., 24 (2008), 73--82.

\bibitem {TL} S.J. Tang and X.J. Li, Necessary conditions for optimal control of stochastic systems with random
jumps, SIAM J. Control Optim, vol. 32, no. 5, (1994), 1447--1475.

\bibitem{Watanabe} S. Watanabe, It\^{o}'s theory of excursion point processes and its developments,
Stochastic Processes Appl. 120, no. 5 (2010), 653--677.

\bibitem {YM} J. Yin and X. Mao, The adapted solution and comparison theorem for backward stochastic differential
equations with Poisson jumps and applications, J. Math. Anal. Appl., 346, 2 (2008), 345--358.

\bibitem {YS} J. Yin and R. Situ, On solutions of forward-backward stochastic
differential equations with Poisson jumps, Stoch. Anal. Appl., 21 (2003), 1419--1448.

\bibitem {Y} J. Yong, Finding adapted solutions of forward-backward stochastic differential equations--method of continuation, Prob. Th. \& Rel. Fields,
 107 (1997), 537--572.

\bibitem {ZSG} Q. Zhu, Y. Shi and X.J. Gong, Solutions to general forward-backward doubly stochastic differential equations,
Applied mathematics and mechanics, vol. 30, no. 4 (2009), 517--526.

\bibitem {ZS} Q. Zhu and Y. Shi, Forward-backward doubly stochastic differential equations and related stochastic partial differential equations, Science China Mathematics, vol. 55, no. 12 (2012), 2517--2534.

\end{thebibliography}
\end{document}